\documentclass[a4paper,11pt]{article}
\usepackage{epsfig}
\usepackage{amsmath}
\usepackage{amssymb}
\usepackage{amscd}

\newcommand{\ud}{\mathrm{d}}

\begin{document}
\title{ON THE GEOMETRY OF THREE-DIMENSIONAL BOL ALGEBRAS WITH SOLVABLE ENVELOPING LIE ALGEBRAS OF SMALL DIMENSION}
\author{THOMAS BOUETOU BOUETOU}
\maketitle
\newpage
\tableofcontents
\newpage
\begin{center}

\section{INTRODUCTION}

\end{center}

The notion of quasigroup and the tangential structure appearing to it
have been intensively developing for the 40 past years. The quasigroup as
a nonassociative algebraical structure is naturally the generalization
of the notion of group. It first appeared in the work of R. Moufang
(1935) \cite{mou}. She obtained some identities. The smooth local loops were
first introduced in the work of Malcev A.I. \cite{malc}, in connection with
the generalization of Lie group.

Left binary algebras were later on called Malcev algebras. It is well
known that modern differential geometry and nonassociative algebras
are interacting on one another. The notion of binary-ternary
operation tangent to the given geodesic loop connected with an
arbitrary point in the affine connected space was introduced in Akivis
works see \cite{ak1,ak2}. Several mathematicians worked in the development of
differential geometry and the study of different classes of
quasigroups and loops \cite{loos,sabmik4,fed2,shele1,shele2,bar,bru,kik2}.

In 1925 E. Cartan in his research lay the beginning of investigation
of symmetric spaces. Today the given structure plays an important role,
in differential geometry and its application. The question arising
here is that of, the description and the classification of symmetric
space naturally leading to the classification of the corresponding
algebraical structure.

The survey on geometry of fiber space stimulates the interest for
special types of 3-webs, in particular Bol 3-Webs and the tangential
structures appearing to it: Bol algebras \cite{mike2,mike6}. In
connection with it
the idea of the description of collections of three-dimensional Bol
algebras elaborates. For separate classes of Bol algebras see
\cite{alhs,fedo}. In
this investigation our approach on Bol algebras will be base on the
classification of solvable Lie triple systems \cite{boue1}. As we know Bol
algebra can be seen as a Lie triple system equipped with an additional
bilinear anti-commutative operation which verify a series of
supplementary conditions.

Bol algebras appear under the infinitesimal description of the class
of local smooth Bol loop, in the work of L.V. Sabinin and
P.O. Mikheev. The interest of studying smooth Bol loop is connected to
the fact that, the geodisical loop see(Sabinin \cite{sab4}) of local
 symmetric
affine connected space, verifies the left Bol identity and automorphic
inverse identity. Its an exact algebraic analog of construction of
symmetric spaces. In particular, the velocity space in the theory of
relativity (STR), is a Bol loop relatively to the additional law of
velocity \cite{sabmik5}. Thus, relying on what is stated above we
are ready to
formulate the purpose of this work: Classification of Bol algebras of
dimension 3 with solvable Lie algebras of dimension $\leq
5$, with accuracy to isomorphism and with accuracy to isotopic which
include:

\begin{enumerate}
\item the  description of three-dimensional Lie triple systems and
their corresponding Lie algebras with invomorphisme.
\item The description of three-dimensional Bol algebras linked with
the distinguished Lie triple systems above.
\item The construction of the object describing the isotopy of Bol algebras.
\item The description of Bol 3-Webs, connected with the selected Bol algebras.
\end{enumerate}
The methods used in these investigations are the methods of geometry
of fiber bound and that of non associative algebras. This work
bears a theoretic character. The obtained results can find their use
in differential geometry, theory of quasigroup and loop, same as in
various applications in physics and mechanics.

The obtained results are new and published in the work of Bouetou
Bouetou T. \cite{boue1,boue2,boue3,boue4,boue5,boue6,boue7,boue8,boue9}, and
also in Bouetou Bouetou T. Mikheev P.O. \cite{boumik}.

This thesis constitutes: an introduction, three chapters and their
references. Each chapter is divided into paragraphs, and the paragraph
are divided into subsections. The enumeration of formulas is connected
with the paragraph.

Chapter I constitutes of seven paragraphs.

In \S 1. {\bf Bol loops} a short discussion about the notion and
terminology in quasigroup theory and Bol loop is given.

In \S 2. {\bf Isotopy of loops} here the definition of isotopy is given. In
short form, and proposition, it states the fundament of isotopy of
loops.

In \S 3. {\bf Local analytic Bol} loops here the notion of local analytic
Bol loops is discussed. In a definitive manner the notion of Bol
algebras as W-algebras \cite{ak2} verifying a system of identities is
introduced. Then the discussion of the imbedding of local analytic Bol
loop into a local Lie group and the imbedding of Bol algebra into Lie
algebra is dealt with. It shows how to calculate the operations $ \xi
\cdot \eta$ and $ (\xi, \eta, \chi)$ of Bol algebra $\mathfrak{B}$ in
the term of enveloping Lie algebra $\mathfrak{G}$, his subalgebra
$\mathfrak{h}$ and subspace $\mathfrak{B}$.

In \S. 4. {\bf 3-Webs and coordinates loop of 3-Webs} we see through
definitions and short propositions the fundamentals of 3-Webs theory
and the coordinates loop of 3-Webs being stated.

In \S .5. {\bf Isotopy of Bol algebras}. In this paragraph a correct
generalization of the notion of isotopy of (global loops to the case
of local analytic Bol loops is given. In connection with this, the
definition of isotopy of Bol algebras is given. The following theorem
is stated and proven:

Let $ B(x)$ and $ \widetilde{B(0)}$ be global analytic Bol loops, and let
 their tangent Bol algebras be isotopic, then  $ \widetilde{B(0)}$ is locally
 isomorphic, to an analytic Bol loop analytically isotopic to $ B(x) $.

In \S .6. {bf About the classification of Bol algebras}. In this paragraph
we state
 in a short way the method used in this investigation.

In \S .7. {\bf Isocline Bol algebras}. This given class of Bol algebras is a
 particular case of Bol algebras. Any Bol algebra is said to be isocline if
 and only if it verifies the plane axiom.

Chapter II. This chapter is devoted to the classification of solvable Lie
 triple systems of dimension 3. It consists of 4 paragraphs.

\S .1. {\bf Some information about lie triple systems}. Following the work
 \cite{loos} and \cite{trof}
 the direct and inverse construction of imbedding of Lie triple system
into a Lie algebra is presented.

\S .2. {\bf Solvable and simple Lie triple systems} On the basis of the
analogue theorem of
 Levi-Malcev theorem \cite{koba}, which states that if $ \mathfrak{M}$ is
 a Lie triple
 system and $ \mathfrak{G}=\mathfrak{M} \dotplus \mathfrak{h} $ his canonical
 enveloping Lie algebra and $r$-the radical of Lie algebra$ \mathfrak{G}$, then
 in $ \mathfrak{G}$ there exist a semisimple subalgebra $ p $, complementary
 to $ r $, such that;
$$
 \mathfrak{M}=\mathfrak{M'}\dotplus \mathfrak{M''}
$$
where   $ \mathfrak{M'}= \mathfrak{M}  \cap r$,  $ \mathfrak{M'}= \mathfrak{M}  \cap r $ and $ \mathfrak{M''}=\mathfrak{M}  \cap p $, $ \mathfrak{h}=\mathfrak{h'}\dotplus \mathfrak{h''} $,  $ \mathfrak{h'}= \mathfrak{M}  \cap r $,
 $ \mathfrak{h''}=\mathfrak{M}  \cap p$ .

One can select 3 cases of 3-dimensional Lie triple systems:
\begin{enumerate}
\item Semi-simple.
\item Splitting.
\item Solvable.
\end{enumerate}

\S .3. {\bf Classification of solvable Lie triple systems}. In this paragraph
the
 following theorem is stated and proven: with accuracy to isomorphism there
 exists only 7 different types of Lie triple systems of dimension 3.

\S .4. {\bf Some example of Bol algebras with solvable trilinear operation}.
 Here
 two examples are given and one remark. The first example is of Bol algebras,
 obtained from the classification of 3-dimensional Lie algebra (Bianchi
 classification \cite{dou}). The second is of Bol algebras, obtained from right
 alternative algebras \cite{mike1}. The remark is about the loops obtained
from the work \cite{fedo}.

Chapter III. This chapter consists of 7 paragraphs. Each paragraph is based on
 the result obtained after the proof of theorem II.2.

\S .1. {\bf Bol algebras with trivial trilinear operation of type I}. Here it
is
 proven that, with accuracy to isomorphism there exists 6 Bol algebras of type I
, and their corresponding 3-webs are described.

\S .2. {\bf Bol algebras with trivial trilinear operation of type II}. Here
 Bol algebras of type  with enveloping Lie algebras of dimension 4  and 5
 are considered. It's shown that with accuracy to isomorphism there exists
 4 Bol algebras with 4 dimensional enveloping Lie algebra, two of them are
 isotopics. Also, with accuracy to isomorphism there exists 3 families of Bol
 algebras with 5 dimensional canonical Lie algebra. Their corresponding
 3-Webs are described.

\S .3. {\bf Bol algebras with trilinear operation of type III}. For a better
 investigation of such algebras the case has been divided into type $ III^-$
 and type $III^+$. In the limit of each of them, Bol algebras with 4
dimensional
 enveloping Lie algebras are considered. It's shown that in type $III^+$ and
 type $III^-$ with accuracy to isomorphism there exists two families and one
 exceptional Bol algebra. And the identification of Bol algebras obtained in
 \cite{mike1} is made. The corresponding 3-Webs are described.

\S .4. {\bf Bol algebras with trilinear operation of type IV}. this examination
 is
 divided into two cases type $IV^+$ and type $IV^-$. in the limit of each of
 them,
 the examination of Bol algebras with a 4-dimensional canonical Lie algebra
 is
 given; therefore it's shown that there exist with accuracy to isomorphism 2
 families of Bol algebras. Their corresponding 3-Webs are described.

\S .5. {\bf Bol algebras with trilinear operation of type V}. Here also the
 consideration is divided into type $ V^+$ and type $V^-$. In the limit of each
 of them the examination of Bol algebras with a 4- dimensional canonical Lie
 algebras is given. And, the corresponding 3-webs are output.

\S .6. {\bf Bol algebras with trilinear operation of type VI}. The canonical
 enveloping Lie algebra of this type of Bol algebras is of dimension 5,
 therefore it's shown that with accuracy to isomorphism there exist 3
 families of Bol algebras.

\S .7. {\bf Bol algebras with trilinear operation of type VII}. The canonical
 enveloping Lie algebra of this type of Bol algebras is of dimension 5,
 therefore it's shown that with accuracy to isomorphism there exist 6
 families and 4 exceptional Bol algebras. And the 3-Webs corresponding
 to each one of them are output.

\newpage


\section{CHAPTER I}


\subsection{BOL LOOP}


{\bf Definition I.1.1.} A set $ Q $ with a fixed element $ e$ and binary
 operations $ (\cdot)$ and $(\setminus)$, satisfying the following conditions

\begin{equation}
\forall a \in Q  \; e\cdot a=a \cdot e= a
\end{equation}
\begin{equation}
\forall a, b \in Q  a\cdot (a \setminus b)=b, a \setminus(a \cdot b)=b
\end{equation}
is called a left loop with two-sided identity element $e$.

{\bf Definition I.1.2.} A left loop with two-sided identity element
 $ (Q, \cdot, \setminus, e) $ is called a Bol loop if the left Bol identity
$$
a \cdot(b \cdot(a \cdot c))=(a \cdot(b \cdot a))\cdot c  \; \; (3)
$$
holds for all
$a$,$ b$, $c$  $ \in Q $.

The following properties are known from the theory of Quasigroups and loops.

\begin{enumerate}
\item $ \forall a \in Q $ there is the unique element $ a^{-1} \in Q $ such
that:
$$
a^{-1}\cdot a= a \cdot a^{-1}=e
$$
\item $ \forall a, b \in Q $ the solution of equation $ a \cdot x=b$ is:
$$
x=a\setminus b=a^{-1}b,
$$
the solution of the equation $ x \cdot a=b $ also exists is uniquely determined
 and is of the form $ x=a^{-1} \cdot ((a \cdot b) \cdot a^{-1}) $.
In particular in any Bol loop one can define the operation of right division by
 $ a \setminus b = a^{-1} \cdot((a \cdot b) \cdot a^{-1}) $, then
 $ (a \cdot b) \setminus b= a $ and $ (a \setminus b) \cdot b= a $
\item $ \forall a\in Q $, $ \forall m \in \mathbb{N} $ defines $ a^m $
recursively by $ a^0=e $, $ a^{m}=a^{m-1} \cdot a $ and $ a^{-m}=(a^{-1})^{m} $
, the left mono-alternative property:
$$
a^{m} \cdot (a^{r} \cdot b)=a^{m+r}\cdot b \; \; \; (4)
$$
holds $ \forall$ $ a$,$b$ $ \in Q $ and$\forall m$, $r$
$\in \mathbb{Z}$, in particular, Bol loops are power associative
 (mono-associative)
$$
\forall a \in Q \;and\; \forall m, r \in \mathbb{Z},\; a{^m} \cdot a^{r}=a^{m+r}\; (5)
$$
\item Besides the left Bol identity (3) one can consider the right Bol identity
$$
((c \cdot a) \cdot b) \cdot a=c \cdot((a \cdot b) \cdot a) \; \; (6)
$$
A loop $ ( M, \cdot, \setminus, /, e ) $ satisfying the properties (3) and (6)
 simultaneously, is called a Moufang loop. Moufang loops are diassociative,
 that means, any two elements of a Moufang loop generate a subgroup. Note that
 diassociativity implies the following properties:
$$
\forall a, b \in M \; \; (a \cdot b)^{-1}=b^{-1} \cdot a^{-1}
$$

$$
a \cdot (b \cdot a)=(a \cdot b) \cdot a \; \; (elasticity) \; \; \; \;  (7)
$$

\item The Moufang condition (3)+(6) is equivalent to each of the identities:

$$
a \cdot (b \cdot(a \cdot c))=((a \cdot b) \cdot a) \cdot c,
$$
$$
 a \cdot (b \cdot (c \cdot b))=((a \cdot b) \cdot c) \cdot b \; \; \; (8)
$$
\end{enumerate}

The Bol condition (3) was first clearly distinguished from stronger
 Moufang condition (8) in the work \cite{bar}. We may say without
 exaggeration that the Bol loop construction, its various particular cases,
 modification and generalizations(Moufang loops, M-loops, Bruck loops,
 semi-Bol loop, etc \ldots)
are the heart of the modern theory of quasigroup and loops. Partly it is
explained by isotopic invariance of the Bol property (3) that is, each loop
isotopic to a Bol loop is a Bol loop \cite{belou,bru} the same result is
true for
 Moufang loops. It is known that the category of left loop is equivalent
to the category of equipped homogeneous spaces \cite{sab1,sab2}. This fact
 is probably the original
cause of the theoretical physicist arising interest in the flexible and
economic construction of loops of transformations. It seems to us that the
 loops, in particular Bol loops, will be a tool of the newest natural
 sciences. (An amazing fact should be pointed out that the symmetrical
space is practically some Bol loop!).


\subsection{IMBEDDING OF LOOPS IN GROUPS}


Let $ (Q, \cdot, \setminus, e) $ be a left loop and $ L_{Q}$- group, generated
 with the set of left translations
$$
L_{x}:L_{x}y=x \cdot y, \forall x,y \in Q, L_{x^{-1}}y=x\setminus y
$$
The group $as_{l}(Q)$, generated with the set $ {l(x,y)} $$\forall x, y \in Q$
 where
$$ l_{(x,y)}=L^{-1}_{(x \cdot y)} \cdot L_{x} \cdot L_{y} $$

is called associant. In fact, it will be called the left associant of the left
quasigroup.  $as_{l}(Q)$ is the subgroup of
the group of  permutation $ \sigma_{Q} $. For $h\in  \sigma_{Q}, q\in Q $ we
will denote the action of $ \sigma_{Q} $ to $Q$ by $hq$. It's clair that if
$(Q, \cdot, \setminus)$ is associative, then $as_{l}(Q)=\{Id_Q\}$
 One can verify that $ asl(Q) \subset L_Q \subset \sigma_Q$.
 In the work \cite{sab2}  the following application is defined
$ m: Q\times \sigma_Q \longrightarrow \sigma_Q$
 $$m_{q}(h)=m(q,h) =L^{-1}_{hq}.h.L_q. h^{-1},  \forall q \in Q, \forall h\in \sigma_Q $$
 it's clear that $m_Q (h)={id_Q}$ if and only if $ h$ -is an
 authomorphis of loop $ Q(\cdot)$. If $m_Q (as_l (Q))=\{id_Q\}$, then $Q$ -is a
 special loop. A group $ H\subset \sigma_Q$ is called left transassociant if a
 left quasigroup $ Q(\cdot)$, is such that $as_l (Q)\subset H$ and
$m_Q (H)\subset H$ and left
 transassociant of loop $ Q(\cdot)$, if $ H$ conserves the right neutral
 element $ e$ of a left loop$ Q(\cdot): he=e, \forall h \in H, As_l (Q)$ -the
minimal
 left transassociant of a left loop, obtained from $ as_l (Q)$ by an
 indefinite extension with the help of the operator $m_Q$ and the
 operation generating the subgroup from the subset. If in the left
 loop $ e$ is a right neutral element, then $ as_l (Q)e={e}, he=e$ and
 consequently $ m_Q (h)e={e}$ and,$ As_l (Q)\subset L_Q$. One can prove that for
 a left loop $ Q(\cdot), m_q (as_l (Q)) \subset as_l (Q)$, $m_q (L_Q)\subset L_Q$
and, hence, $as_l (Q)=As_l (Q$) and $L_Q$ transassociant for $Q (\cdot)$.

\underline{{\bf Definition 1.1.3.}} [Sabinin L.V. \cite{sab2}] Let $Q(\times)$ be a left loop,
 and $ H$ his transassociant, the set $ Q\times H$, equipped with the internal
 composition law
 $$(q_1 ,h_1 ) (q_2 ,h_2 )= (q_1 (h_1 q_2 ),\phi(q_1 ,h_1, q_2, h_2) )$$

where $$\phi(q_1 ,h_1, q_2, h_2)=l_{q_1 \cdot h_1,q_2} \circ m_{q_2}h_1 \circ h_1 \circ h_2$$
 is called a semi-direct product of $ Q$ and $ H$ he is not by $ Q \boxtimes H$.

\underline{{\bf Proposition 1.1.1.}} [Sabinin L.V., \cite{sab2}] The semi-direct product
 of $Q\boxtimes H$ of left loop $Q$ and his transassociant is a group.

 About the general properties of Bol loops see [7,66,84,90]

\begin{center}
\subsection{ ISOTOPIC LOOPS}
\end{center}

 Loops $(Q,\cdot,e)$ and $(Q',\circ,e)$ are called isotopic if the maps
$ U,V$ and $W: Q \longrightarrow Q'$ exist and bijective such that
 $$xU \circ yV=(x\cdot y)W \; \; \forall  x,y \in Q$$

 the triplet $(U,V,W)$ is called isotopy  of $Q$ and $Q'$. If the loops
$(Q,\cdot,e)$
and $(Q', \circ) $ coincide and $W=Id_Q$ then the isotopy $(U,V,Id_Q)$ is
called principal.

The following proposition holds.

\underline{{\bf Proposition I.2.1}} [Bruck \cite{bru}] If the  $(Q', \circ, \epsilon) $
is isotopic to the loop
 $(Q,\cdot,e)$ then it's isomorphic to the principal isotopy $(Q,\bot,b.a)$ where
the operation $\bot$ is define by:
$$x \bot y=(x/a) \cdot (b\setminus y) $$

Let  $(Q,\cdot,e)$-be an abstract loop with the left Bol identity. it's known
see [Robinson D.A.\cite{rob2}] that each loop isotopic to $Q$ will also verify
the left
Bol identity. In addition we have the following proposition:

\underline{{\bf Proposition I.2.2}} [Bruck \cite{bru}] If the
 $(Q', \bot f) $ is isotopic to the loop
 $(Q,\cdot,e)$ with the left Bol identity then it's isomorphic to the principal isotopy $(Q,\bot,f^2)$ where
the operation $\bot$ is define by:
$$x \bot y=(x/f) \cdot (f\setminus y) $$

\subsection{Local analytic Bol loops}

Let $Q$ be a set where is defined the operation $\times$(multiplication),
the left division $\setminus$ and the unity element $e$.

We will say that $(Q,\times,\setminus,e)$ is a left local loop with a double
sided unity if $Q$ is a topological space with a fixed element $e$ and for a
certain neighborhood $\mathcal{U}$ of the element $e$ are defined continuous
map
$$\mathcal{U}\times \mathcal{U}\longrightarrow Q:(x,y)\longrightarrow x\times y$$
and
$$\mathcal{U}\times \mathcal{U}\longrightarrow Q:(x,y)\longrightarrow x\setminus y$$
verifying the condition:
\begin{enumerate}
\item if $x \in \mathcal{U}$, then $e \times x=x \times e=x$
\item if $x,y, x\times y \in \mathcal{U}$ then $x\setminus (x\times y)=y$
\item if $x,y,x\setminus y \in \mathcal{U}$ then $x\times (x\setminus y)=y $
\end{enumerate}

If $Q$ is a smooth manifold of class $C^K (0\leq k \leq \omega)$ and the maps
$((x,y)\longrightarrow x\times y),((x,y)\longrightarrow x\setminus y)$ also of
the  class $C^K$ then $Q$ is called left local loop of  class $C^K $. If at the
place of $\mathcal{U}$ consider all the space $Q$, then $Q$ is called left
topological loop or correspondently left loop of smoothness $C^k$.

\underline{{\bf Definition I.3.1.}} We will say that $Q$ is a local Bol loop of
 class
 $C^K (0\leq k \leq \omega)$, if $Q$ is a manifold of class  $C^K$ with a
fixed element $e$ and for any neighborhood $\mathcal{U}$ of the element $e$
 are defined the maps of class $C^k$ $ \mathcal{U}\ni x\longrightarrow x^{-1} $
and  $ \mathcal{U}\times \mathcal{U}\ni (x,y)\longrightarrow x\cdot y\in Q $
verifying the conditions
\begin{itemize}
\item if $x\in \mathcal{U}$ then $e\cdot x=x\cdot e=x$
\item if $x,y,x\cdot y \in \mathcal{U}$ then $x^{-1}\cdot x=x\cdot x^{-1}=x, x^{-1}\cdot (x\cdot y)=y$
\item if $x,y,z,x\cdot z, x\cdot y, x\cdot yx,y\cdot xz \in \mathcal{U}$ then
the left Bol identity $(x\cdot yx)z=x(y\cdot xz)$ is verify.
\end{itemize}

\underline{{\bf Definition I.3.2.}} Let $V$ be a finite dimensional vector
space where, it's
defined a bilinear and a trilinear operation $x\cdot y$ and $<x,y,z>$. We will
say that $V$ is called a $W$-algebra if the following identities are verify:
\begin{enumerate}
\item 4. $x\cdot x=0$
\item 5. $<x,x,x>=0$
\item 6. $xy\cdot z+yz\cdot x+zx\cdot y=<x,y,z>+<y,z,x>+<z,x,y>-<y,x,z>-<z,y,x>-<x,z,y>$.
\end{enumerate}

The operation of composition of any local loop $(Q,\times,e)$ of class $C^3$
has in the neighborhood of the unity element $e$ the following coordinates
expression(with accuracy up to the 3 order)
$$
(x\times y)^i=x^i+y^i+\tau^i_{jk}x^j y^k +\mu^i_{jkl}x^j x^k y^l +\nu^i_{jkl}x^j y^k y^l+......
$$

from the basic tensor of the local loop
$$ \alpha^i_{jk}=\tau^i_{[jk]}$$
$$
\beta^i_{jkl}=2\mu^i_{jkl}-2\nu^i_{jkl}+\alpha^m_{jk} \alpha^i_{mk}-\alpha^i_{jm}\alpha^m_{ik}
$$

equipping the tangent space $ V=T_{e}(Q) $ with the composition law
$$
[x,y]^{i}=2\alpha^{i}_{jk}x^{j}y^{k},
$$
$$
(x,y,z)^{i}=\beta^{i}_{jkl}x^{j}y^{k}z^{l}.
$$

If $ \alpha(t)$, $\beta(t)$, $ \gamma(t) $ be smooth curves in loop $ Q $ at
 $ t=0 $ going through the point $e $ with tangent vectors $ x$, $y$, $z$
accordingly then

$$
(\beta(t) \times \alpha(t)) \setminus(\alpha(t) \times \beta(t))=t^{2}[x,y]+0(t^{2})
$$
$$
[\alpha(t) \times(\beta(t) \times \gamma(t))] \setminus[(\alpha(t) \times \beta(t)) \times \gamma(t)]=t^{3}<x,y,z>+0(t^{3})
$$
The tensors $ \alpha^{i}_{jk} $ and $ \beta^{i}_{jkl} $ ( or what equivalent
 the operation $[,], <, , >$) are defining the $ C^{w} $ Bol loop with local
isomorphism accuracy under the fulfillment of conditions:

4'. $ [x,x]=0 $\\
5'. $ <x,x,y>=0 $ \\
6'. $ [[x,y],z]+[[y,z],x]+[[z,x],y]=2<x,y,z>+2<y,z,x>+2<z,x,y> $

7.  $ [<x,y,z>,\eta]-[<x,y,\eta>,z]+<[x,y],\eta,z>-<[x,y],z,\eta>+<x,y,[\eta,z]>=0 $,

8.  $ (x,y,(z,\xi,\eta))=((x,y,z),\xi,\eta)+(z,(x,y,\xi),\eta)+(z,\xi,(x,y.\eta)) $, \\

where $ (x,y,z)=2<x,y,z>-[[x,y],z] $.

\underline{{\bf Definition I.3.3.}} \cite{sabmik1} $  W $ -algebra $V$, the
 basic operations for which the following conditions are satisfied:

5''. $ <x,x,y>=0 $

6''. $ xy\cdot z+yz\cdot x+zx \cdot y=2<x,y,z>+2<y,z,x>+2<z,x,y> $

7''. $ <x,y,z>\cdot \xi-<x,y,\xi>\cdot z+<x\cdot y, \xi,z>-<x \cdot y, z, \xi>+<x,y, \xi \cdot z>=0 $

8.  $ (x,y,(z,\xi,\eta))=((x,y,z),\xi,\eta)+(z,(x,y,\xi),\eta)+(z,\xi,(x,y.\eta)) $, \\

where $ (x,y,z)=-2<x,y,z>+xy \cdot z $.

Will be call a Bol algebras.

The system of the identities 5'', 6'', 7'' ' can be rewritten in the equivalent
 form:

9. $ (x,x,y)=0 $

10. $ (x,y,z)+(y,z,x)+(z,x,y)=0 $

11. $ (x,y,z) \cdot \xi-(x,y, \xi) \cdot z+(z,\xi,x \cdot y)-(x,y,z \cdot \xi)+xy \cdot z \xi=0 $

here, one can consider an arbitrary binary-ternary algebra, of finite dimension
 with main operation for which $  x \cdot y $ and $ (x,y,z) $ satisfy the
 identities (8)-(10) as a Bol algebra. This point of view is more preferable
 because in this case one can regard the conditions (8)-(11) as the identities
 defining a Lie triple systems  with the composition law $ (x,y,z) $.

Let $<Q,\cdot, e>$ be a loop with the left Bol identity and $G=<LQ>$ his
canonical enveloping group for the Bol loop $Q$ and the subgroup $H$ which
correspond to the Lie group $\mathfrak{G}$ and his subalgebra $\mathfrak{h}$.
In the point of the local left coset section of the class $ G \bmod H$.
$Q=\exp \mathcal{U}$ where $\mathcal{U}$ is sufficiently small neighborhood
of the $0$ in $\mathfrak{B}$ (Bol algebra). let's introduce the law of
composition $\times$: $x\times y=(x/e)y$ in result we obtain a loop
$<Q,\times,e>$ with left Bol identity and the $W$-algebra is isomorphic to the
to the original Bol algebra.

Let the local section $B=\exp \mathcal{U}$ of left coset class $G \bmod H$
verifying the conditions stated above $x \in B$ and $\widetilde{H}=x\triangle H\triangle x^{-1} $
for $x$ sufficiently near to $e$, we have $\widetilde{H}\cap B=\{e\}$ that is
way in $B$ one can define a new law of composition
$a\top b=\widetilde{\prod}_B (a\triangle b)$  where $\widetilde{\prod}_B$ is
the projection on $B$ parallel to $H$. Since $Q$ is the section $G/H$ and
$f=x^2$ then $(Q,\top)$ is obtain by the rotation $x\cdot H\cdot x^-1 $. For
$y$ relatively near to $e$ let's consider the operation:
$$a\bot b=(a/y)\times(y\times b).$$  According to Robinson \cite{rob3}, any
arbitrary loop isotopic to the loop $(Q,\times,e)$ is isomorphic to the loop
$(Q,\bot,e)$ for any $y$ see [79-80].

\underline{{\bf Theorem I.31.1}} Local analytic Bol loops $B$ and $B'$ are
local isomorphic
if and only if their correspondent Bol algebras $V$ and $V'$ are isomorphic.

\underline{{\bf Theorem I.3.2}} An arbitrary Bol algebra is a $W$-algebra of
an analytic Bol loop.

\subsection{IMBEDDING OF A LOCAL ANALYTIC BOL LOOP INTO A LOCAL LIE GROUP}

Let $ (G, \Delta, e) $- be a local Lie group and $ H $- be one of its subgroups,
and let's  denoted the  corresponding Lie algebra and subalgebra by
$ \mathfrak{G}$ and $ \mathfrak{h}$. Consider a vector subspace
$\mathfrak{B}$ such that

$ \mathfrak{G}=\mathfrak{h}\dotplus \mathfrak{B}$.
 Let $ \Pi:G\longrightarrow G\setminus H $ be the canonical projection and let
$ \Psi$ be the restriction of mappings composition $ \Pi \circ \exp$,
to $ \mathfrak{B}$. Then there exists such a neighborhood $ \mathcal{U}$ of
the
point $ O $ in $\mathfrak{B}$ such that, $\Psi$ maps it diffeomorphically into
the
neighborhood $\Psi (u)$ of the coset $\Pi (e)$ in $ G \setminus H $ [46].\\

$$
\begin{CD}
\mathfrak{G} @<{i}<< \mathfrak{B}\\
@V{\exp}VV  @VV{\Psi}V \\
G @>{\pi}>> G\setminus H
\end{CD}
$$

Introduce a local composition law:

$$
 a \star b= \Pi_{B}(a\Delta b)
$$

On points of local cross-section $ B=\exp \mathcal{U} $ of left cosets of $G$
 mod $H$, where $ \prod_{B}=\exp \circ \Psi^{-1} \circ\Pi: G \longrightarrow B $
is the local projection on $ B $ parallel to $ H $ which, puts into correspondence
 every element $ a \in B $ so that $ g=a \Delta p $, where $  p \in H $.

\underline{{\bf Proposition I.3.1}}[46] Let's assume that for any
$ a, b \in B=\exp \mathcal{U} $ sufficiently close to the point $e$, and
$ a\Delta b \Delta a \in B $; then the local analytic loop
$ (B, \times, e) $ satisfies the left Bol condition.

\begin{center}

\subsection{IMBEDDING OF BOL ALGEBRA INTO LIE ALGEBRAS: ENVELOPING LIE ALGEBRA OF BOL ALGEBRAS}

\end{center}

Let the local cross-section $ B=\exp \mathcal{U} $ of left cosets $ G
\setminus H $ satisfy the condition of Proposition I.3.1 above. It is
interesting to calculate the, operations of the Bol algebra $
\mathfrak{B} $ tangent to the local analytic Bol loop $ ( B, \star, e)$
in terms of a Lie algebra $ \mathfrak{G} $, its subalgebra
$\mathfrak{h} $ and vector subspace $ \mathfrak{B} $. Introduce in $ G $
normal coordinates then:

$$
a \star b=a+b+\frac{1}{2}[a,b]_{\mathfrak{B}}+0(2)
$$
$\forall a $ ,$ b \in B $
where $ [a,b]_{\mathfrak{B}} $ is the projection of $ [a,b] $ on $ \mathfrak{B} $ parallel to $ \mathfrak{h} $ thus
$$
\xi \cdot \eta=[\xi, \eta]_{\mathfrak{B}}
$$
$$
(\xi, \eta, \chi)=[[\xi, \eta], \chi]
$$
$$
<\xi, \eta, \chi>=-\frac{1}{2}[[\xi, \eta], \chi]+\frac{1}{2}[[\xi, \eta]_{\mathfrak{B}}, \chi]_{\mathfrak{B}}.
$$
One can find the correctness of the following propositions [30,44]:\\

\underline{{\bf Proposition I.3.2.}} \cite{sabmik1} Any local analytic Bol
loops
$ B $ and
$ B' $ are locally isomorphic, if and only if their corresponding Bol algebras
$ \mathfrak{B} $ and $ \mathfrak{B'} $ are isomorphic.

\underline{{\bf Proposition I.3.3.}} \cite{sabmik1} An arbitrary Bol algebra is
 a tangent
$ W $-algebra of some local analytic loop with the left Bol loop.

\underline{{\bf Proposition I.3.4.}} \cite{sabmik1}  The property
$ a\Delta b \Delta a \in B $
is held for
 any $ a$, $b \in B $ sufficiently close to $e $, if and only if
 $ [[\xi, \eta], \zeta ] \in \mathfrak{B}
 \forall \xi, \eta, \zeta \in \mathfrak{B} $.

Let $\mathfrak{B}$ be a finite Bol algebra over $ \mathbb{R}$, the basic
 operations of which are $\xi \cdot \eta $ and $ (\xi, \eta, \zeta) $,
$ \mathfrak{G} $-finite Lie algebra over $\mathbb{R} $, $\mathfrak{h}$ -
 subalgebra of $ \mathfrak{G}$ and
$i$:$ \mathfrak{B} \longrightarrow \mathfrak{G} $ a linear mapping such that
 $ i(\mathfrak{B}) \in \mathfrak{B}$.

$ \mathfrak{G} \dotplus \mathfrak{h} $ (direct sum of vector spaces)\\

$ [[\mathfrak{B}],\mathfrak{B}], \mathfrak{B}] \in \mathfrak{B} $ and
$\forall \xi, \eta, \zeta \in \mathfrak{B}$

$ \xi \cdot \eta=[\xi, \eta]_{\mathfrak{B}} $
$ (\xi, \eta, \zeta)=[[\xi, \eta], \zeta] $

where $[\xi, \eta]$ denotes the result of commutation of vectors in $ \mathfrak{G}$
 and $ [ \xi, \eta]_{\mathfrak{B}} $ denotes projection of the vector
 $ [\xi, \eta] $ on $\mathfrak{B}$ parallel to $\mathfrak{h}$.

 In that case we will talk about the enveloping pair $ (\mathfrak{G}, \mathfrak{h})$ Lie algebra of Bol algebra$ \mathfrak{B}$ or, in other words enveloping
 Lie algebra $\mathfrak{G}$ of Bol algebra $\mathfrak{B}$.

Let   $ (\mathfrak{G}, \mathfrak{h}) $ be an enveloping pair of Lie algebra of
 Bol algebra$ \mathfrak{B}$. Let us identify $\mathfrak{B}$ with a vector
 subspace
 $ i(\mathfrak{B})$ into $\mathfrak{G}$, and let us consider the subalgebra
 $ \mathfrak{G'}=\mathfrak{B}\dotplus [ \mathfrak{B}, \mathfrak{B}] $ into
 $\mathfrak{G}$, and the subalgebra
 $ \mathfrak{h'}=\mathfrak{h}\cap \mathfrak{G'} $ into$ \mathfrak{G'}$, then
the pair $ ( \mathfrak{G'}, \mathfrak{h'} ) $ is also for enveloping for a
 Bol algebra $\mathfrak{B}$.

By the construction of the Lie algebra which is a canonical enveloping for
$ \mathfrak{B} $, it is better to use the construction made in [45]. For such
 Bol algebra $ \mathfrak{B}$, there exists a Lie algebra
 $  \mathfrak{\widetilde{G}} $, and an envoluting automorphism
 $ \tau \in Aut \mathfrak{\widetilde{G}} $ such that $ ( \tau^{2}=Id) $, linear
 injection map $ i: \mathfrak{B} \longrightarrow \mathfrak{\widetilde{G}} $
 and a subalgebra $ \mathfrak{\widetilde{h}} $ in $ \mathfrak{\widetilde{G}} $,
 such that ( we are identifying $ i(\mathfrak{B}) $ with $\mathfrak{B}$).

$ \mathfrak{\widetilde{G}} =\mathfrak{\widetilde{G_{-}}}+\mathfrak{\widetilde{G_{+}}} $, where $ \mathfrak{\widetilde{G_{-}}}=\mathfrak{B} $;
 $ \mathfrak{\widetilde{G}}=\mathfrak{B} \dotplus \mathfrak{\widetilde{h}} $,
 $ < \mathfrak{B}>=\mathfrak{\widetilde{G}} $

$$
\xi \cdot \eta= \prod_{\mathfrak{B}}[\xi, \eta]= [\xi, \eta]_{\mathfrak{B}}
$$
$$
(\xi, \eta, \zeta)=[[\xi, \eta], \zeta]
$$

Lie algebra $ \mathfrak{\widetilde{G}}$ is an enveloping Lie algebra for Bol
algebra $ \mathfrak{B}$, but in general not canonical enveloping because
$ \mathfrak{\widetilde{h}}$ may contain an ideal $ I $ of Lie algebra
$ \mathfrak{\widetilde{G}}$. Hence the canonical enveloping Lie algebra
$ \mathfrak{G}$ for a Bol algebra $ \mathfrak{B}$ is obtained by factorizing
$ \mathfrak{\widetilde{G}}$  with the ideal $ I $. Therefore after the
factorization of $ \mathfrak{\widetilde{G}}_{+} \setminus I $ and $\mathfrak{B}$
 ( we identify $\mathfrak{B}$ and $ \mathfrak{B} \setminus I $) in general
interest, let us note that the construction see [45] follows that $ dim \mathfrak{\widetilde{G}} \leq dim \mathfrak{B} \wedge \mathfrak{B}+dim \mathfrak{B} $.

In our case $ dim \mathfrak{B}=3 $, that is why under the examination of the
 corresponding canonical enveloping Lie algebra $\mathfrak{G}$ we must consider
 the case:
$$
dim[ \mathfrak{B}, \mathfrak{B}]=0, 1, 2, 3.
$$


\subsection{THREE-WEBS COORDINATES LOOP OF THREE-WEBS}


\underline{{\bf Definition I.4.1}} Let $ W $ be a $ C^{\infty} $
 smooth manifold of
 $ 2N $ dimension, equipped with such foliation $ \lambda_{i} $, $ i=1,2,3 $
codimension $ N $, for any point from $ W $, then there exists a neighborhood
 $\mathcal{U}$ containing a point, such that the fiber of any two different
 foliations have in $ \mathcal{U} $ not more than one common point. In that
 case $ ( W, \lambda_{1}, \lambda_{2}, \lambda_{3}) $ is called local three-web
.
 The three-web $ ( W, \lambda_{1}, \lambda_{2}, \lambda_{3}) $ is called
 global, if for any two fiber from the different families $ \lambda_{i} $ and
 $ \lambda_{j} $ $ (i,j=1,2,3, i \neq j )$ have exactly one common point.

About the three-web theory see $[2,3,4,5,6,9]$.

Let $ P $ be a fixed point of the three-Web $ W $ and let
 $ \mathcal{F}_{a} \in \lambda_{1} $ and  $ \mathcal{F}_{b} $ two fibers of web
, passing through $ P $ and let $ q(a,b)=c  \in  x_3 $. Let us examine the
 map:
$$
u=q(x,b), v=q(a,y), z=z. \; \;\; \; \; \; (1)
$$
The foliations $ \lambda_1$, $ \lambda_2 $, abide us to define, the operation of
 multiplication in $ \lambda_{3}$ such that:
$$
u \cdot v=q(x,y)=q(q^{-1}(u,b),q^{-1}(a,v))
$$
where $ c $ denotes the unity.\\
If the neighborhood $\mathcal{U} \ni p $ is so small, the map (1) is a
 bijection which defines some isotopic transformation of quasigroup $ q $. The
so obtained isotope, is called principal isotope. It's having a neutral element
 that is why it's a loop $ l(a,b) $.

\begin{equation}
\raisebox{-0.5\height}{\epsfig{file=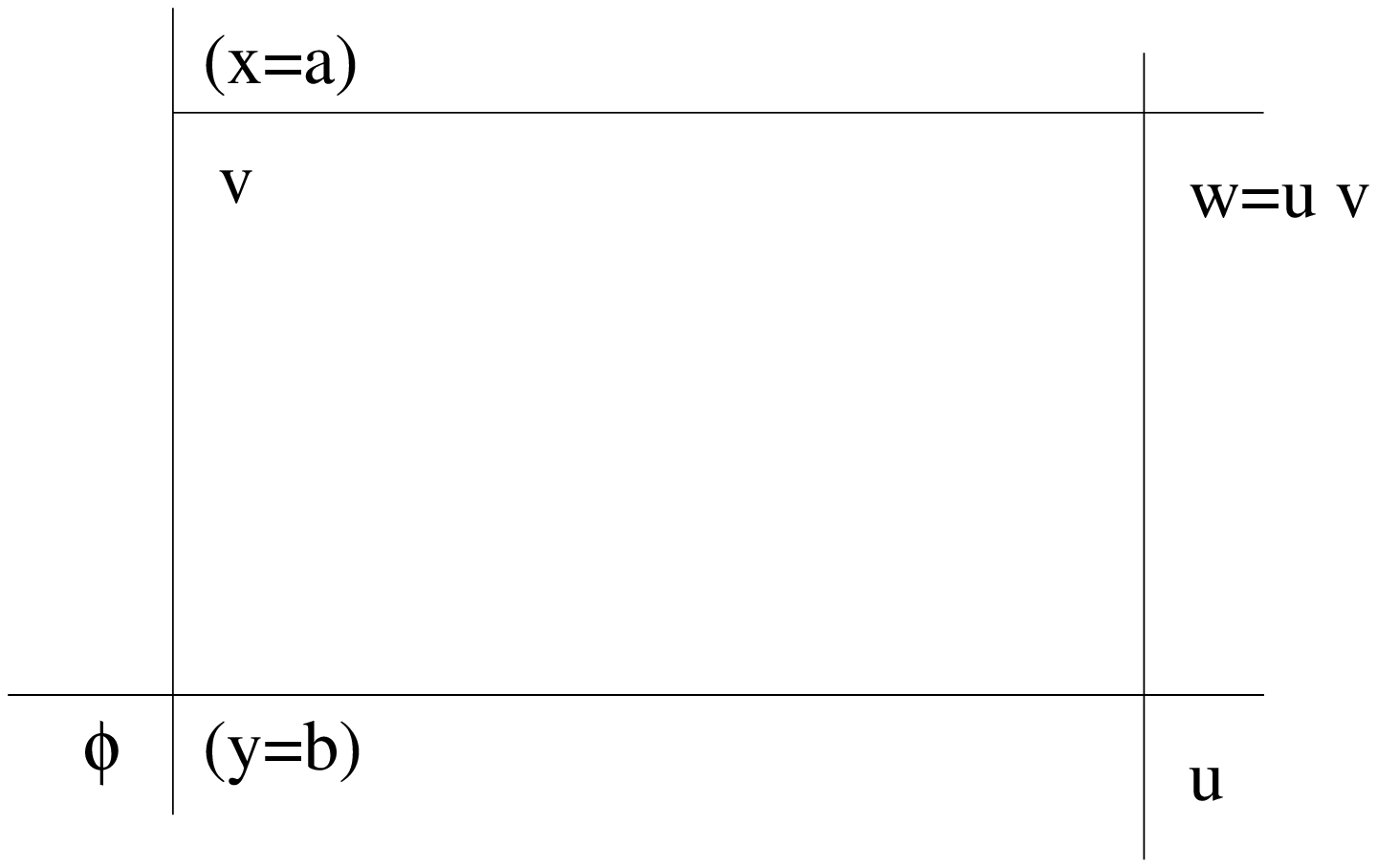,width=6cm}}
\end{equation}

\begin{center}

CONDITIONS FOR CLOSURES

\end{center}

In the three-web theory, the condition for closure of some figures formed by the wed and point in their intersection plays an important role.\\
Let us examine a three-web $ W $ and let $ z=q(x,y) $  it's coordinate loop. It is said that on the web $ W $ the figure of closure verifies one of the following types if:

\begin{equation}
\raisebox{-0.5\height}{\epsfig{file=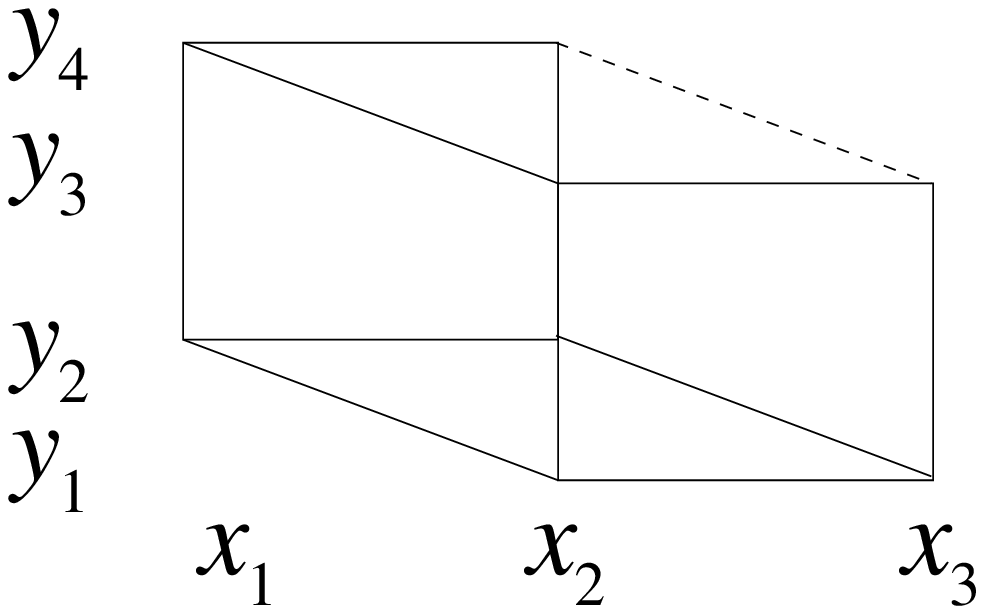,width=3cm}}
\left.
\begin{array}{c} q(x_1,y_2)=q(x_2,y_1)\\
                 q(x_2,y_2)=q(x_3,y_1)\\
                 q(x_1,y_4)=q(x_2,y_3)\\
\end{array}
\right\}\longrightarrow q(x_2,y_4)=q(x_3,y_3)
\end{equation}
\begin{equation}
\raisebox{-0.5\height}{\epsfig{file=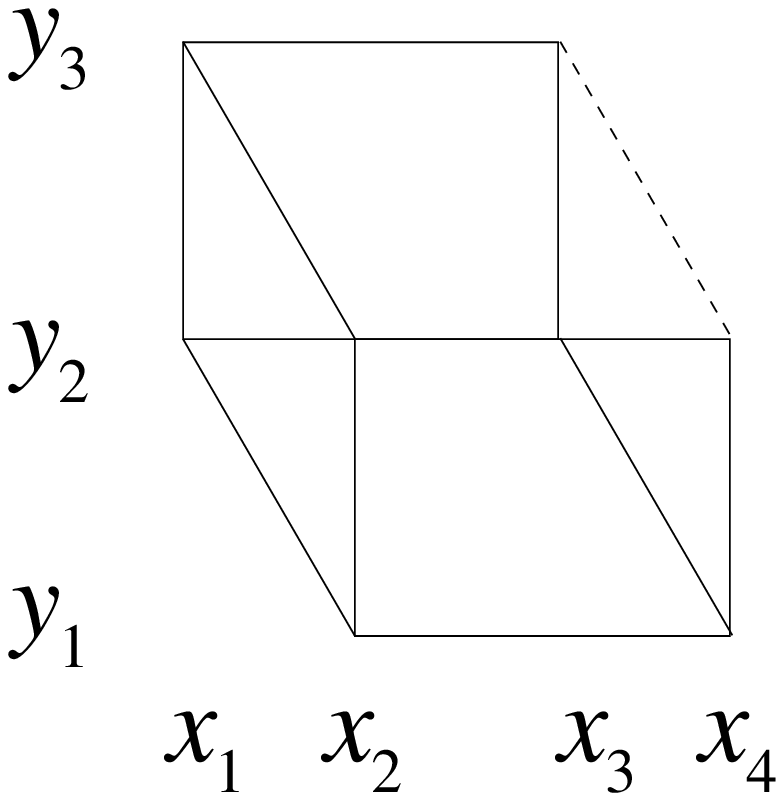,width=3cm}}
\left.
\begin{array}{c} q(x_1,y_2)=q(x_2,y_1)\\
                 q(x_2,y_2)=q(x_1,y_3)\\
                 q(x_4,y_1)=q(x_3,y_2)\\
\end{array}
\right\}\longrightarrow q(x_3,y_3)=q(x_4,y_2)
\end{equation}
\begin{equation}
\raisebox{-0.5\height}{\epsfig{file=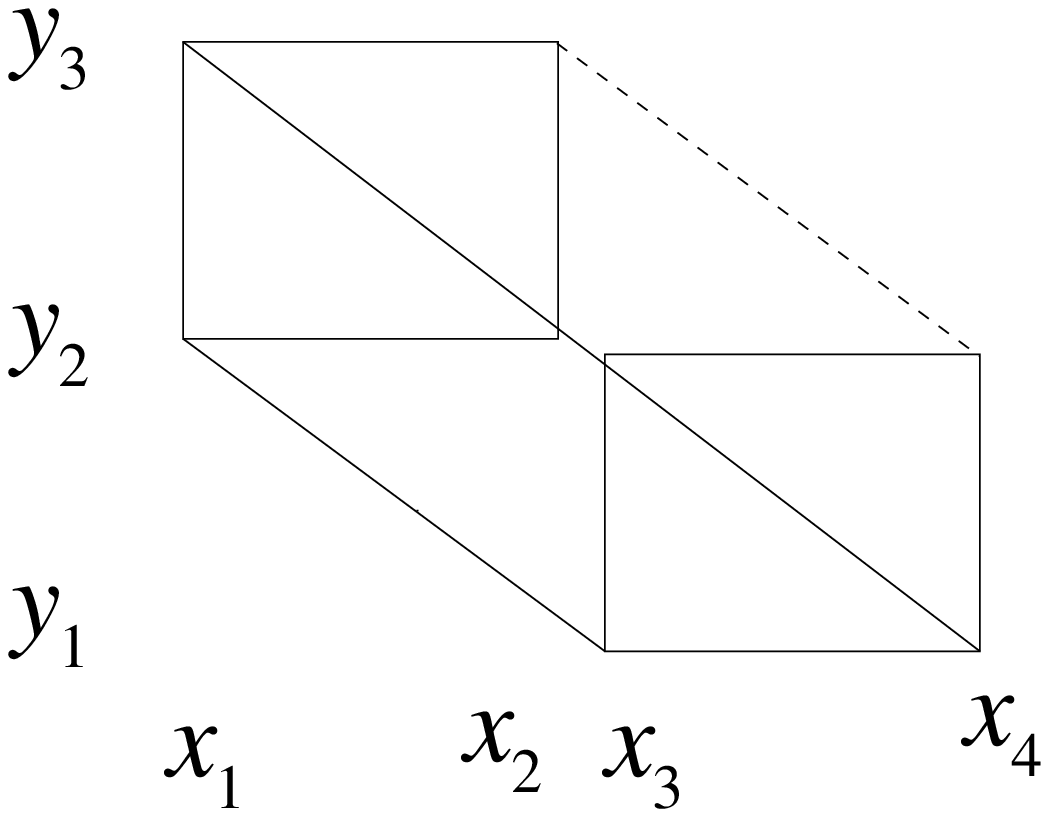,width=3cm}}
\left.
\begin{array}{c} q(x_1,y_2)=q(x_2,y_1)\\
                 q(x_1,y_4)=q(x_2,y_3)\\
                 q(x_3,y_2)=q(x_4,y_1)\\
\end{array}
\right\}\longrightarrow q(x_3,y_4)=q(x_4,y_3)
\end{equation}

These figures are respectively called left, right and mean Bol figure and it is denoted by $ B _{l} $, $ B_{r} $, $ B_{m} $.

\underline{{\bf Proposition I.4.1.}} \cite{ak4} Three-web $ W $ is a Bol web if
 any of its
coordinate loops is a Bol loop
Any local analytical loop can be realized as some coordinate loop of a certain three-web.\\
\underline{{\bf Proposition I.4.2}} \cite{ak4}
\begin{enumerate}
\item The coordinate loop of a global three-web is an isotope.
\item The three-web, corresponding (analytical) to isotopic (global) loops are isomorphic.
\end{enumerate}

\subsection{ISOTOPIC BOL ALGEBRAS}

Below we give a correct generalization of the notion of isotopy of (global)
loops to the case of local analytic Bol loops. Our approach is base on a
construction of embedding of Bol loops into a group and on an interpretation
of isotopic loops in terms of their enveloping groups [48,86].\\
\underline{{\bf Definition I.5.1.}} Let $ \mathfrak{B}$ and $\mathfrak{\widetilde{B}} $ be
two Bol algebras and let $ \mathfrak{G}=\mathfrak{B} \dotplus \mathfrak{h} $
and $ \mathfrak{\widetilde{G'}}=\mathfrak{\widetilde{B}}\dotplus \mathfrak{\widetilde{h}}$
be their canonical enveloping Lie algebras. The algebras $\mathfrak{B}$ and
$ \mathfrak{\widetilde{B}}$ will be called isotopic if there exit such a Lie
algebras isomorphism $ \Phi: \mathfrak{G} \longrightarrow \mathfrak{G} $,
such that $ \Phi(\mathfrak{G})=\mathfrak{G} $ and $ \Phi (\mathfrak{h}) $
coincide with the image of the subalgebra $\mathfrak{\widetilde{h}}$ in
$ \mathfrak{\widetilde{G}}$ under the action of an inner automorphism
 $ Ad \xi $, $ \xi \in \mathfrak{\widetilde{G}}$ i.e.
$$
 \Phi (\mathfrak{h})=( Ad_{\xi}) \mathfrak{\widetilde{h}}
$$

It is clear that the notion of isotopy is not an equivalence relation to Bol algebra manifolds.

\underline{{\bf Theorem I.5.1}} Let $ B(\times) $ and $ \widetilde{B}(\circ) $
be global analytic Bol loops, and let their tangent Bol algebras be isotopic, then $ \widetilde{B}(\circ) $ is locally isomorphic to an analytic Bol loop analytically isotopic to $ B(\times) $.\\
Proof.\\
Let  $ \mathfrak{G} $ be the canonical Lie algebra enveloping the Bol algebra  $ \mathfrak{B} $ tangent to the analytic Bol loop $ B(\times) $.  $ \mathfrak{G}=\mathfrak{B} \dotplus \mathfrak{h} $ ( direct sum of vector spaces). There exists a Lie group $ G $ with $ \mathfrak{G} $ as its tangent Lie algebra, a closed subgroup $ G_{O} $ corresponding to the subalgebra $ \mathfrak{h} $ and an analytical embedding $ i: B \longrightarrow G $ such that the composition law takes the following form [48,86]
$$
a \times b=\prod(a,b),
$$
where $a \times b$ denotes the composition of elements $a$, $ b $ in $ B $ and $ \prod: G \longrightarrow B $ is the projection on $ B $ parallel to  $ G_{O} $.\\ Let us suppose that the Bol algebra $ \widetilde{\mathfrak{B}} $ tangent to the loop $ \widetilde{B}(\circ) $ is isotopic to the Bol algebra $ \mathfrak{B} $, i.e. there exist Lie algebra isomorphism
 $ \Phi : \mathfrak{G} \longrightarrow  \mathfrak{\widetilde{G}}$, such that  $ \Phi ( \mathfrak{G})=  \mathfrak{\widetilde{G}} $ and  $ \Phi ( \mathfrak{h}) $ coincide with the image of the subalgebra $  \mathfrak{\widetilde{h}}$ in $  \mathfrak{\widetilde{G}} $ under the action of an inner automorphism $ Ad \xi $, $ \xi \in  \mathfrak{\widetilde{B}} $, $ \Phi ( \mathfrak{h})=(Ad \xi)\mathfrak{\widetilde{h}} $.

Let us introduce the element $ y=exp(\xi) $ and the subgroup $ \widetilde{G_{O}}=yG_{O}y^{-1} $. The analytic loop $ B(\star) $
$$
a \star b=\widetilde{\prod}(ab), a,b \in B
$$

where $ \widetilde{\prod}: G \longrightarrow B $ is the projection on $ B $ parallel to  $ \widetilde{G_{O}} $, is a Bol loop whose tangent Bol algebra is isomorphic to  $ \widetilde{\mathfrak{B}} $. In particular, $ B( \star) $ and   $ \widetilde{B}(\circ) $ are locally isomorphic. The operations on $ B $ are isotopic. Indeed let us introduce an analytic diffeomorphism $ \Omega: B \longrightarrow B $, $ a \longmapsto (y \times a) \setminus y $, then:\\
$$
\Omega^{-1}_{y}(\Omega_{y}(a) \star \Omega_{y}(b))=(L_{y})^{-1}(\Omega_{y}(a) \times L_{y}b)=Y^{-1}\times [((y\times a)\setminus y) \times (y \times b)]
$$
$$
=[y^{-1} \times (((y \times a) \setminus y) \times y^{-1})]\times (y^{2} \times b)
$$
$$
=[y^{-1} \times ((y^{-1} \times ((y^{2}\times a)\times y^{-1}))]\times(y^{2} \times b)
$$
$$
=[y^{-1} \times [y^{-1} \times ((y^{2}\times a)\times y^{-2})]\times(y^{2} \times b)
$$
$$
=[y^{-2} \times ((y^{2} \times a) \times y^{-2})]\times (y^{2} \times b)
$$
$$
=(a \setminus y^{2}) \times (y^{2} \times b)
$$
$$
=(a \setminus y) \times (y \times b)
$$
therefore $ \Omega_{y}((a \setminus y) \times (y \times b))= \Omega_{y}(a) \star \Omega_{y}(b) $.
hence $B(\times)$ is isotopic to  $ \widetilde{B}(\circ) $ and by the
diffeomorphism $\Omega$ they are isomorphic. Hence the theorem is proved.
\begin{center}
\subsection{ABOUT THE CLASSIFICATION OF BOL ALGEBRAS}
\end{center}

Let   $ \mathfrak{G}=\mathfrak{B} \dotplus \mathfrak{h} $ be Lie algebra in the
 involutive decomposition $  \mathfrak{h'} $ a subalgebra in  $ \mathfrak{G} $
such that  $ \mathfrak{G}=\mathfrak{B} \dotplus \mathfrak{h'} $
$$
a \cdot b=\prod [a,b]
$$
$$
(a,b,c)=[[a,b],c]
$$
where $ \prod:   \mathfrak{G} \longrightarrow \mathfrak{B}  $ is a projection
on
  $ \mathfrak{B}  $ parallel to  $ \mathfrak{h'} $ and [,] commutator in Lie
 algebra  $ \mathfrak{G} $.\\
That's why the classification of Bol algebra lead to the classification of subalgebra  $ \mathfrak{h'} $, in the enveloping Lie algebra  $ \mathfrak{G} $ (not necessarily canonical) for a Lie triple system  $ \mathfrak{M} $.\\
Below, we will examine the classification of Bol algebras with isomorphism accuracy and isotopic accuracy [66]. The classification with isotopic accuracy is more crude, than the classification with isomorphism accuracy. However the notion of isotopy of Bol algebra is opening a new connection between non isomorphic Bol algebras.
\subsection{ISOCLINE BOL ALGEBRAS}

One can prove that any vectorial space  $ \mathfrak{B} $, equipped with operation
$$
\forall \xi, \eta, \zeta \in \mathfrak{B}\; \; \xi \cdot \eta=\alpha (\xi)\eta-\alpha(\eta)\xi,
$$
$$
<\xi, \eta, \zeta>=\beta(\xi,\zeta)\eta-\beta(\eta, \zeta)\xi \; \; (I)
$$
where $ \alpha: \mathfrak{B} \longrightarrow \mathbb{R}$-linear form,  $ \beta: \mathfrak{B}\times \mathfrak{B} \longrightarrow \mathbb{R}$-bilinear symmetric
form. Is a Bol algebra.

\underline{{\bf Definition I.7.}} Bol algebra of view (I) is called isocline.

One can prove see [46], that any Bol algebra is called isocline if and only if
 it verifies the plane axiom that means any two-dimensional vectorial subspace
(plane) is a subalgebra.

In particular any two-dimensional Bol algebra is isocline.

If $ dim \mathfrak{B}=3 $ and $ \alpha=0 $, then depending on the rank, and the
signature of the form $ \beta $, we obtained 5 non trivial and non isomorphic
types of Lie triple systems.
\newpage
\section{CHAPTER II}
\begin{center}
CLASSIFICATION OF SOLVABLE 3-DIMENSIONAL LIE TRIPLE SYSTEMS EXAMPLE OF 3-DIMENSIONAL BOL ALGEBRAS
\end{center}
\subsection{ABOUT LIE TRIPLE SYSTEMS}

\underline{{\bf Definition II.1.1}} The vector space $ \mathfrak{M} $
(finite over the field of real numbers $ \mathbb{R}$) with trilinear operation (x,y,z) is called Lie
triples system if:
$$
(x,x,y)=0
$$
$$
(x,y,z)+(y,z,x)+(z,x,y)=0
$$
$$
(x,y,(u,v,w))=((x,y,u),v,w)+(u,(x,y,v),w)+(u,v,(x,y,w))
$$
Let $ \mathfrak{M} $ be a Lie triple system, the subspace
 $ \mathfrak{D}  \subset  \mathfrak{M} $ is called subsystem if
 $ (\mathfrak{D}, \mathfrak{D}, \mathfrak{D})\subset \mathfrak{D} $, and is ideal, if $ ( \mathfrak{D}, \mathfrak{M}, \mathfrak{M}) \subset \mathfrak{D} $.
 The ideals are the Kernel of the homomorphism of the Lie triple system [26,52].

\underline{\bf{Example}} For a typical way of construction of a Lie triple system see
 in [26,52].

Let $ \mathfrak{G}$ be a Lie algebra (finite over the field of real numbers
$ \mathbb{R}) $  and $ \sigma $-an involutive automorphism, then
$$
\mathfrak{G}=\mathfrak{G}^{+} \dotplus \mathfrak{G}^{-}
$$
where $ \sigma | \mathfrak{G^{+}}=Id $ and  $ \sigma | \mathfrak{G^{-}}=-Id $,
as for any element $x$ from $ \mathfrak{G}$ can be written in the form:
$$
x=\frac{1}{2}(x+\sigma x)+\frac{1}{2}(x-\sigma x),
$$
where  $x+\sigma x  \in \mathfrak{G}^{+} $,  $ x-\sigma x \in \mathfrak{G}^{-}$ and $ \mathfrak{G}^{+} \cap \mathfrak{G}^{-}=0 $.

The following inclusions are held:

$$
[\mathfrak{G^{+}}, \mathfrak{G^{+}}] \subset \mathfrak{G^{+}}, [\mathfrak{G^{+}}, \mathfrak{G^{-}}] \subset \mathfrak{G^{-}}, [\mathfrak{G^{-}}, \mathfrak{G^{-}}] \subset \mathfrak{G^{+}}.
$$
Such that the subspace $ \mathfrak{G^{-}} $ turns into a Lie triple system
relatively under the operation $ (x,y,z)=[[x,y],z] $.

The inverse construction [26].

Let $\mathfrak{M}$ be a Lie triple system and define by
$$h(X,Y): z \longrightarrow (X,Y,Z)$$
a linear transformation of the space $\mathfrak{M}$ into it selves where
$X,Y, \in \mathfrak{M}$.

Let $H$ be a subspace in the space of linear transformations, of Lie triple
systems $\mathfrak{M}$ of the transformations of the form $h(X,Y)$. The
vector space $\mathfrak{G}=\mathfrak{M} \dotplus H$, become a Lie algebra
relatively to the commutator $[A,B]=AB-BA$, $[A,X]=-[X,A]=AX$;
$[X,Y]=h(X,Y)$ where $A,B \in H$, $X,Y \in \mathfrak{M}$.

Let define the mapping $\sigma$ with the condition
$\sigma(A)=A, A\in H$ and $\sigma(X)=-X, X \in \mathfrak{M}$, then
$\sigma$- is an involutif automorphism of Lie algebra $\mathfrak{G}=\mathfrak{M} \dotplus H$.

The algebra $\mathfrak{G}$ constructed above from the Lie triple system, is
called universal enveloping Lie algebra of the Lie triple system $\mathfrak{M}$.

\underline{\bf Definition II.1.2} The derivation of the Lie triple system
$\mathfrak{M}$, is called
 the linear transformation
$\mathfrak{d}: \mathfrak{M} \longrightarrow \mathfrak{M}$ such that
$$ (X,Y,Z)\mathfrak{d}=(X\mathfrak{d},Y,Z)+(X,Y\mathfrak{d},Z)+(X,Y,Z\mathfrak{d}).$$

One can verify that, the set $\mathfrak{d}(\mathfrak{M})$ of all the
derivation of the Lie triple systems $\mathfrak{M}$ is a Lie algebra of the
linear transformations acting on $\mathfrak{M}.$

\underline{\bf Definition II.1.3} The imbedding of a Lie triple system
$\mathfrak{M}$ into a Lie
algebra $\mathfrak{G}$ is called the linear injection
$R: \mathfrak{M} \longrightarrow \mathfrak{G}$ such that
$(X,Y,Z)=[[X^R,Y^R],Z^R]$.

The imbedding $R$ of the Lie triple system $\mathfrak{M}$ into the Lie algebra
$\mathfrak{G}$ is called canonical, if the envelope of the image of the set
$\mathfrak{M}^R$ in the Lie algebra $\mathfrak{G}$ coincide with
$\mathfrak{G}$ and $h$ does not contain trivial ideals of Lie algebra
$\mathfrak{G}.$ Let us note that if the Lie triple system $\mathfrak{M}$ is
a subset of the Lie algebra $\mathfrak{G}$, then $(X,Y,Z)=[[X,Y],Z]$
and $[\mathfrak{M},\mathfrak{M}]$ is a subalgebra of the Lie algebra
$\mathfrak{G}$ hence
$\mathfrak{M}+[\mathfrak{M},\mathfrak{M}]$- is a Lie subalgebra of
$\mathfrak{G}$ and the initial imbedding $R$ can be consider as canonical in
$\mathfrak{M}^R+[\mathfrak{M}^R,\mathfrak{M}^R]$; this lead us to formulate
the following proposition:

\underline{{\bf Proposition II.1.1}} For any finite Lie triple system
$\mathfrak{M}$
over $\mathbb{R}$, there exist one and only up to automorphism accuracy, one
canonical imbedding to the Lie algebra.

\begin{center}
SOLVABLE AND SEMISIMPLE LIE TRIPLE SYSTEM
\end{center}

Following \cite{lis}:, let $\Omega$- be an ideal of the Lie triple system
$\mathfrak{M}$, we assume $\Omega^{(1)}=(\mathfrak{M},\Omega,\Omega)$ and,
  $\Omega^{(k)}=(\mathfrak{M},\Omega^{(k-1)},\Omega^{(k-1)})$

\underline{{\bf Proposition II.1.2}} \cite{lis} For all natural number $k$, the subspace
$\Omega^{(k)}$ is an ideal of $\mathfrak{M}$ and we have the following
inclusions:
$$ \Omega \supseteq \Omega^{(1)} \supseteq ........\supseteq \Omega^{(k)}$$

Proof

$$
(\Omega^{(1)},\mathfrak{M},\mathfrak{M})=((\mathfrak{M},\Omega,\Omega),\mathfrak{M},\mathfrak{M})\subseteq ((\mathfrak{M},\Omega,\mathfrak{M}),\Omega,\mathfrak{M})+[[[\mathfrak{M},\Omega],[\mathfrak,\Omega]],\mathfrak{M}]
$$
according to the definition of a Lie triple system
$$
(\Omega^{(1)},\mathfrak{M},\mathfrak{M})\subseteq (\Omega,\Omega,\mathfrak{M})+[(\mathfrak{M},\Omega,\mathfrak{M}),[\mathfrak{M},\Omega]]\subseteq (\mathfrak{M},\Omega,\Omega)+(\mathfrak{M},\Omega,\Omega)=\Omega^{(1)}
$$
that means  $ \Omega^{(1)}$ is an ideal of $\mathfrak{M}$
further more $\Omega^{(k)}=(\Omega^{(k-1)})^{(1)}$ hence each $\Omega^{(i)}$ is
an
ideal in $\mathfrak{M}$.

\underline{{\bf Definition II.1.4}} The ideal $\Omega$ of a Lie triple system
$\mathfrak{M}$ is
called solvable, if there exist a natural number $k$  such that $\Omega^{(k)}=0.$

\underline{{\bf Proposition II.1.3}} \cite{lis} If $\Omega$ and $\Theta$ are two solvable ideals
 of a Lie triple system $\mathfrak{M}$ then  $\Omega+\Theta$ is also a solvable
ideal in $\mathfrak{M}$.

Proof

using the definition of a Lie triple system, the following inclusion hold:
$(\Theta+\Omega)^{(1)}\subseteq (\mathfrak{M},\Theta,\Theta)+(\mathfrak{M},\Omega,\Omega)+(\mathfrak{M},\Theta,\Omega)+(\mathfrak{M},\Omega,\Theta)\subseteq \Theta^{(1)}+\Omega^{(1)}+\Theta \cap \Omega$.

Assume for every natural number $k$ the following inclusion holds:

$(\Theta+\Omega)^{(k)}\subseteq \Theta^{(k)}+\Omega^{(k)}+\Theta \cap \Omega$

by induction let's prove that its holds for $(k+1)$

$(\Theta+\Omega)^{(k+1)}=(\mathfrak{M},(\Theta+\Omega)^{(k)},(\Theta+\Omega)^{(k)})\subseteq(\mathfrak{M},\Theta^{(k)}+\Omega^{(k)}+\Theta \cap \Omega,(\Theta \cap \Omega))\subseteq \Theta^{(k+1)}+\Omega^{(k+1)}+\Theta \cap \Omega$

hence the result

\underline{{\bf Definition II.1.5}} The radical of a Lie triple system denoted
by
$\mathfrak{R}(\mathfrak{M})$, is called the maximal solvable ideal of the Lie
triple system $\mathfrak{M}.$

A Lie triple system $\mathfrak{M}$ is called semi-simple if $\mathfrak{R}(\mathfrak{M})=0$.

\underline{{\bf Theorem II.1.1}} \cite{lis} If $\mathfrak{R}$ is a radical in$\mathfrak{M}$
then $(\mathfrak{M}\setminus \mathfrak{R})$ is semisimple. And if
$\Omega $ is an ideal in $\mathfrak{M}$ such that $(\mathfrak{M}\setminus \mathfrak{R})$ is semisimple then $\Omega \supset \mathfrak{R}$.

\underline{{\bf Proposition II.1.4}} \cite{lis} The enveloping Lie algebra, of
a solvable Lie
triple system is solvable. And if a Lie triple system has some solvable
enveloping Lie algebra, it is solvable.

\underline{{\bf Theorem II.1.2}} If $\mathfrak{M}$ is a semisimple Lie
triple system, then
the universal enveloping Lie algebra $\mathfrak{G}$ is semisimple.

\underline{{\bf Theorem II.1.3}} Let $\mathfrak{M}$ be a Lie triple system and
$\mathfrak{G}=\mathfrak{M} \dotplus \mathfrak{h}$ his canonical enveloping Lie
algebra and $\mathfrak{r}$- the radical of the Lie algebra $\mathfrak{G}$.
In $\mathfrak{G}$ there exist a subalgebra $\mathfrak{P}$ semisimple
supplementary to with $\mathfrak{r}$ such that:
$$
\mathfrak{M}=\mathfrak{M}' \dotplus \mathfrak{M}''\; (direct\; sum\; of\; vectors \;spaces)
$$
where
$$
\mathfrak{M}'=\mathfrak{M} \cap \mathfrak{r} -radical\; of\; the\; Lie\; triple\; system\; \mathfrak{M}
$$
$$
\mathfrak{M}''=\mathfrak{M} \cap \mathfrak{P} -semisimple\; subalgebra\; of\; Lie\; triple \; system \; \mathfrak{M}
$$
$$
\mathfrak{h}=\mathfrak{h}' \dotplus \mathfrak{h}''\; (direct\; sum\; of\; vectors \; spaces)
$$
$$
\mathfrak{h}'=\mathfrak{h} \cap \mathfrak{r}
$$
and
$$
\mathfrak{h}''=\mathfrak{h} \cap \mathfrak{P}\; are\; subalgebra \;in \; \mathfrak{h}
$$
$$
\mathfrak{r}=\mathfrak{M}' \dotplus \mathfrak{h}'
$$
$$
\mathfrak{P}=\mathfrak{M}'' \dotplus \mathfrak{h}''.
$$

\subsection{PROBLEM SETTING}

Let $\mathfrak{M}$ be a Lie triple system and $dim \mathfrak{M}=3$. To be
consistent with the above Theorem the following cases are possible:
\begin{enumerate}
\item  semisimple case\\
$\mathfrak{M}$- semisimple Lie triple System (in fact simple).
 About the classification of such Lie triple systems see \cite{lis,ber, fed1}
\item Splitting case\\
$$\mathfrak{M}=\mathfrak{M}_1 \dotplus \mathfrak{M}_2 $$
where $\mathfrak{M}\equiv \mathbb{R}$- solvable ideal of dimension 1 in $\mathbb{R}$
and $\mathfrak{M}_2$ -semisimple Lie triple system of dimension 2 This type of Lie triple system is not considered in this survey.
\item Solvable case\\
$\mathfrak{M}$ is a solvable Lie triple system. The classification of such
system is given after the next paragraph.
\end{enumerate}
\begin{center}
CLASSIFICATION OF LIE TRIPLE SYSTEM OF DIMENSION 2
\end{center}

For a better survey of such Lie triple system, we will write their trilimear
operation in a special form.

Let $\mathfrak{M}$ be a 2-dimensional Lie triple system we write the trilinear
operation $(X,Y,Z)=\beta (X,Y)Y-\beta (Y,Z)X$
where $\beta:V \times V\longrightarrow \mathbb{R}$ is a symmetric form. The
choice of the basis $V=<e_1,e_2>$ one can reduce the symmetric form to the
view:
\begin{displaymath}
\left(\begin{array}{cc}
\alpha & 0 \\
0 & \nu \\
\end{array}\right),
\end{displaymath}

where $ \alpha, \nu= \pm 1;0$.

By introducing the notation of the derivation
$$
\mathfrak{D}_{x,y}:\mathfrak{M} \longrightarrow \mathfrak{M}
$$
$$
z\longmapsto (x,y,z)
$$
$$
\mathfrak{h}=\left\{\mathfrak{D}_{x,y}\right\}_{x,y \in \mathfrak{M}}.
$$
And

$\mathfrak{G}=\mathfrak{M} \dotplus \mathfrak{h}$- canonical enveloping
Lie algebra of the Lie triple system $\mathfrak{M}$.

Let $\mathfrak{M}=<e_1,e_2>$ then,

$ \mathfrak{h}=\left\{tD_{x,y}\right\}_{t \in \mathbb{R}}$,
$$
e_{1}D=(e_1,e_2,e_1)=\beta (e_1,e_1)e_2
$$
$$
e_{2}D=(e_1,e_2,e_2)=-\beta (e_2,e_2)e_1
$$

$\mathfrak{G}=<e_1,e_2,e_3>$

where
$[e_1,e_2]=e_3, \; [e_1,e_3]=-e_{1}D, \; [e_2,e_3]=-e_{2}D$

Therefore we can have the up to isomorphism accuracy the following five
cases:

\begin{enumerate}
\item (Spherical Geometry)\\

\begin{displaymath}
\left(\begin{array}{cc}
1 & 0 \\
0 & 1 \\
\end{array}\right),
\end{displaymath}

$\mathfrak{G}/\mathfrak{h}\cong so(3)/so(2)$
\item (Lobatchevski Geometry)\\

\begin{displaymath}
\left(\begin{array}{cc}
-1 & 0 \\
0 & -1 \\
\end{array}\right),
\end{displaymath}

$\mathfrak{G}/\mathfrak{h}\cong sl(2,\mathbb{R})/so(2)$
\item Lie triple system with non compact subalgebra $\mathfrak{h}$\\
\begin{displaymath}
\left(\begin{array}{cc}
1 & 0 \\
0 & -1 \\
\end{array}\right),
\end{displaymath}

$\mathfrak{G}/\mathfrak{h}\cong sl(2,\mathbb{R})/\mathbb{R}$
\item Solvable case\\
\begin{itemize}
\item a)

\begin{displaymath}
\beta=\left(\begin{array}{cc}
1 & 0 \\
0 & 0 \\
\end{array}\right),
\end{displaymath}

$e_1 \cdot e_2=e_3, \; e_1 \cdot e_3=e_2$

(This is a Lie algebra $\mathfrak{G}$ of type $g_{3,5}(p=0)$ in \cite{mub2})
\item b)
\begin{displaymath}
\beta=\left(\begin{array}{cc}
-1 & 0 \\
0 & 0 \\
\end{array}\right),
\end{displaymath}

$e_1 \cdot e_2=e_3, \; e_1 \cdot e_3=-e_2$

(This is a Lie algebra $\mathfrak{G}$ of type $g_{3,4}(h=-1)$ in \cite{mub2})
\end{itemize}
\item Abelian case\\
$\beta=0$ $\mathfrak{G}/\mathfrak{h}\cong (\mathbb{R})^2/\left\{0\right\}$
\end{enumerate}


\subsection{CLASSIFICATION OF SOLVABLE LIE TRIPLE SYSTEMS OF DIMENSION 3}


Let $\mathfrak{M}$- be a solvable Lie triple system of dimension 3, and
$\mathfrak{G} \dotplus \mathfrak{h}$ its canonical enveloping Lie algebra
then $\mathfrak{G}$ is solvable in particular $\mathfrak{G}$ posses
a characteristic ideal
$\mathfrak{G}'=[\mathfrak{G},\mathfrak{G}]\vartriangleright \mathfrak{G}$,

$\sigma \mathfrak{G}'=\mathfrak{G}'$,
$\mathfrak{G}' \cap \mathfrak{M}=\mathfrak{M}'=(\mathfrak{M},\mathfrak{M},\mathfrak{M})$ further more $\mathfrak{h} \subset \mathfrak{G}$ since
$\mathfrak{h}=[\mathfrak{M},\mathfrak{M}]$ then

$\mathfrak{G}'=[\mathfrak{G},\mathfrak{G}]=\mathfrak{M}'+ \mathfrak{h}$ where
$\mathfrak{M}'\subsetneq \mathfrak{M}$

Possible situations:

\begin{enumerate}
\item $ dim \mathfrak{M'}=0 $.
 Then $ [\mathfrak{h}, \mathfrak{M}]=\mathfrak{M'}=\{O \} $, that
 means $ \mathfrak{h} \vartriangleright \mathfrak{G} $- ideal, that is why
 $ \mathfrak{h}=\{O \} $ (since $ \mathfrak{G} $- is an enveloping
 Lie algebra) and $ \mathfrak{M}=\mathbb{R} \oplus \mathbb{R} \oplus \mathbb{R} $.
In this case, the Lie triple system is Abelian and we denote it (type I).
\item $   dim \mathfrak{M'}=1 $. Choosing the base $ e_1 $, $ e_2 $, $ e_3 $ in $ \mathfrak{M} $ such that, $ \mathfrak{M'}=<e_{1}> $ and $ \mathfrak{M}=\mathfrak{M'}+<e_{2}, e_{3}> $.

We will introduce in consideration the linear transformation
$ A, B, C: \mathfrak{M}\longrightarrow \mathfrak{M} $, define as:
\begin{displaymath}
A=(e_{1},e_{2},-)=\left(\begin{array}{ccc}
a & b & c\\
0 & 0 & 0\\
0 & 0 & 0\\
\end{array}\right),
B=(e_{2},e_{3},-)=\left(\begin{array}{ccc}
\alpha & \beta & \gamma\\
0 & 0 & 0\\
0 & 0 & 0\\
\end{array}\right),
\end{displaymath}
\begin{displaymath}
C=(e_{3},e_{1},-)=\left(\begin{array}{ccc}
x & -\alpha -c & y\\
0 & 0 & 0\\
0 & 0 & 0\\
\end{array}\right).
\end{displaymath}

And if a skew symmetric form defined as
$ \Phi(-,-): \mathfrak{M} \times \mathfrak{M} \longrightarrow \mathbb{R} $,
such that $ (x,y,e_{1})= \Phi(x,y)e_{1}$.
The dimension of $ \mathfrak{M} $ is 3, that is why there exists
 $ z \in \mathfrak{M} $, $ z \neq 0 $, such that $ \Phi(-,z)=0 $. The following
cases are possible:
\begin{itemize}
\item b.I. The skew-symmetric form $ \Phi $ is non zero and $ z $ is parallel
to $ e_1 $ $ (z \parallel e_{1}) $, then in the base $e_1 $, $ e_2 $, $ e_3 $
the skew-symmetric form $ \Phi $ has the corresponding matrix:

\begin{displaymath}
\left(\begin{array}{ccc}
0 & 0 & 0\\
0 & 0 & \delta\\
0 & -\delta  & 0\\
\end{array}\right),
\end{displaymath}

where $ \delta \neq 0 $. Adjusting $ e_3 $ to $ 1 \setminus \delta e_{3} $, then
$ \Phi( e_{2},e_{3})=1 $, $ \Phi( e_{3},e_{2})=-1 $, so that $ \alpha=1 $,
$ a=x=0$ and

\begin{displaymath}
A=(e_{1},e_{2},-)=\left(\begin{array}{ccc}
0 & b & c\\
0 & 0 & 0\\
0 & 0 & 0\\
\end{array}\right),
B=(e_{2},e_{3},-)=\left(\begin{array}{ccc}
1 & \beta & \gamma\\
0 & 0 & 0\\
0 & 0 & 0\\
\end{array}\right),
\end{displaymath}
\begin{displaymath}
C=(e_{3},e_{1},-)=\left(\begin{array}{ccc}
0 & -1 -c & y\\
0 & 0 & 0\\
0 & 0 & 0\\
\end{array}\right).
\end{displaymath}

The verification of the defined relation of Lie triple system shows that,
with accuracy to the choice of the vector basis $ e_2 $ and $ e_3 $, it is
possible to afford the following realization of the operators $ A $, $ B $, $ C$
 as:

\begin{displaymath}
A=0,
B=(e_{2},e_{3},-)=\left(\begin{array}{ccc}
1 & 0 & 0\\
0 & 0 & 0\\
0 & 0 & 0\\
\end{array}\right),
\end{displaymath}
\begin{displaymath}
C=(e_{3},e_{1},-)=\left(\begin{array}{ccc}
0 &  -1 & 0\\
0 & 0 & 0\\
0 & 0 & 0\\
\end{array}\right).
\end{displaymath}  \; \; \; \; \; \; \; \; (type VII)
\item b.II. The skew-symmetric form $ \Phi $ is non zero and $ z $ is not
parallel to $ e_1 $, let $ z=e_{2 } $, then

\begin{displaymath}
A=(e_{1},e_{2},-)=\left(\begin{array}{ccc}
0 & b & c\\
0 & 0 & 0\\
0 & 0 & 0\\
\end{array}\right),
B=(e_{2},e_{3},-)=\left(\begin{array}{ccc}
0 & \beta & \gamma\\
0 & 0 & 0\\
0 & 0 & 0\\
\end{array}\right),
\end{displaymath}
\begin{displaymath}
C=(e_{3},e_{1},-)=\left(\begin{array}{ccc}
-1 &  -c & y\\
0 & 0 & 0\\
0 & 0 & 0\\
\end{array}\right).
\end{displaymath}

The verification of the defined relations of Lie triple system, show that the
 indicated case has no realization.
\item b.III. The skew-symmetric form $ \Phi $ is trivial. By completing the
 vector $ e_1 $ with the arbitrary choose vector $ e_2 $ and $ e_3 $ up to
 the base, it is possible to realize the operator $ A $, $ B $, and $ C $:

\begin{displaymath}
A=(e_{1},e_{2},-)=\left(\begin{array}{ccc}
0 & b & c\\
0 & 0 & 0\\
0 & 0 & 0\\
\end{array}\right),
B=(e_{2},e_{3},-)=\left(\begin{array}{ccc}
0 & \beta & \gamma\\
0 & 0 & 0\\
0 & 0 & 0\\
\end{array}\right),
\end{displaymath}
\begin{displaymath}
C=(e_{3},e_{1},-)=\left(\begin{array}{ccc}
0 &  -c & y\\
0 & 0 & 0\\
0 & 0 & 0\\
\end{array}\right).
\end{displaymath}

The verification of the defined relations of Lie triple system, show that by
a suitable choice of basis vectors $e_2, e_3$ the following realization
of operators $A,B,C$ is possible:
\begin{itemize}
\item Abelian Type (Type above)
\item $A=C=0$,
\begin{displaymath}
B=(e_{2},e_{3},-)=\left(\begin{array}{ccc}
0 &  0 & 1\\
0 & 0 & 0\\
0 & 0 & 0\\
\end{array}\right).
\end{displaymath}\; \; \; \; \; \; \; \; \;  (Type II)

This Lie triple system, is obtained by a direct multiplication of a Lie triple
system of dimension two $<e_1,e_2>$, by an Abelian one dimensional $<e_3>$.
\item - \begin{displaymath}
A=(e_{1},e_{2},-)=\left(\begin{array}{ccc}
0 & \pm 1 & 0\\
0 & 0 & 0\\
0 & 0 & 0\\
\end{array}\right)
\end{displaymath} B=C=0. \; \; \; \; \; \; \;  (Type III)
\item \begin{displaymath}
A=(e_{1},e_{2},-)=\left(\begin{array}{ccc}
0 &  \pm 1 & 1\\
0 & 0 & 0\\
0 & 0 & 0\\
\end{array}\right)
\end{displaymath}, B=0,
\begin{displaymath}
C=(e_{3},e_{1},-)=\left(\begin{array}{ccc}
0 &  -1 & \mp 1\\
0 & 0 & 0\\
0 & 0 & 0\\
\end{array}\right).
\end{displaymath} \; \; \; \; \; \; \; \; \; \; (Type IV)
\end{itemize}

\end{itemize}
\item $dim \mathfrak{M}'=2$ in particular, $ \mathfrak{M}'$ is a
subsystem of dimension two in $\mathfrak{M}$. one can consider
(refer to \S 2.) $\forall a,b,c \in \mathfrak{M}'$
$$
(a,b,c)=\beta (a,c)b-\beta (b,c)a
$$
where
\begin{displaymath}
\beta=\left(\begin{array}{cc}
\pm 1 & 0\\
0 & 0 \\
\end{array}\right)
\end{displaymath}
and $\mathfrak{M}'$ is a two-dimensional Abelian ideal in $\mathfrak{M}$.
 In the first case the choice of the base $ \mathfrak{M}=<e_{1}, e_{2},e_{3}>$
such that $ \mathfrak{M'}=<e_{1},e_{2}> $, the operations of the Lie triple
 system are reduced to:

\begin{displaymath}
A=(e_{1},e_{2},-)=\left(\begin{array}{ccc}
0 & \pm 1 & x\\
0 & 0 & y\\
0 & 0 & 0\\
\end{array}\right),
B=(e_{2},e_{3},-)=\left(\begin{array}{ccc}
\alpha & \gamma & \mu\\
\beta & \delta & \nu\\
0 & 0 & 0\\
\end{array}\right),
\end{displaymath}
\begin{displaymath}
C=(e_{3},e_{1},-)=\left(\begin{array}{ccc}
\kappa &-x -\alpha & \xi\\
\chi & -y- \beta & \beta\\
0 & 0 & 0\\
\end{array}\right).
\end{displaymath}

The verification of the defined relation of Lie triple system, leads to the
contradiction of the condition that $ dim \mathfrak{M'}=2 $.

Let $ \mathfrak{M'}=<e_{1},e_{2}> $-be a two-dimensional Abelian ideal and
$ e_3 $- the vector completing $ e_1 $, $e_2 $ up to the basis. Then:

\begin{displaymath}
A=(e_{1},e_{2},-)=\left(\begin{array}{ccc}
0 & 0 & a\\
0 & 0 & b\\
0 & 0 & 0\\
\end{array}\right),
B=(e_{2},e_{3},-)=\left(\begin{array}{ccc}
\alpha & \gamma & \mu\\
\beta & \delta & \nu\\
0 & 0 & 0\\
\end{array}\right),
\end{displaymath}
\begin{displaymath}
C=(e_{3},e_{1},-)=\left(\begin{array}{ccc}
\kappa &-a -\alpha & \xi\\
\chi & -b- \beta & \beta\\
0 & 0 & 0\\
\end{array}\right).
\end{displaymath}

 Deforming the vector $ e_1 $ in the limit of the subspace $ <e_{1},e_{2}>$,
the matrix $ A $ can be reduced to the form $ a=b=0 $ or $ a=1 $, $ b=0 $.

The verification of the defined relation of the Lie triple system, in the second
 case leads to the following realization of the operators $ A $, $ B $, $ C $:

\begin{displaymath}
A=0,
B=(e_{2},e_{3},-)=\left(\begin{array}{ccc}
0 & 0 & 1\\
0 & 0 & 0\\
0 & 0 & 0\\
\end{array}\right),
\end{displaymath}
\begin{displaymath}
C=(e_{3},e_{1},-)=\left(\begin{array}{ccc}
0 & 0 & 0\\
0 & 0 & 1\\
0 & 0 & 0\\
\end{array}\right).
\end{displaymath} \; \; \; \; \; \; \; \; \; \; \; \; \; ( type V)

\begin{displaymath}
A=0,
B=(e_{2},e_{3},-)=\left(\begin{array}{ccc}
0 & 1 & 0\\
0 & 0 & \pm 1\\
0 & 0 & 0\\
\end{array}\right),\;
C=0.
\end{displaymath}\; \; \; \; \; \; \; \; \; \; \; \;  (type VI)
\end{enumerate}
In conclusion to the conducted examination we have the following theorem:

\underline{{\bf Theorem II.3.1.}} Let $ \mathfrak{M}=<e_{1},e_{2},e_{3}> $- be
a solvable
Lie triple system of dimension 3, $ \mathfrak{G}$- its canonical enveloping Lie
algebra(solvable), and let $ A, B, C:\mathfrak{M} \longrightarrow \mathfrak{M}$
the linear transformations of the view: $ A=(e_{1},e_{2},-) $,  $ A=(e_{1},e_{2},-) $,  $ B=(e_{2},e_{3},-) $,  $ C=(e_{3},e_{1},-) $: with isomorphism
accuracy, one can find the possibility of the following types:
\begin{itemize}
\item Type I. $ \mathfrak{M}$- Abelian Lie triple system.

\item Type II. \begin{displaymath}
A=0,
C=0,
B=(e_{2},e_{3},-)=\left(\begin{array}{ccc}
0 & 0 & 1\\
0 & 0 & 0\\
0 & 0 & 0\\
\end{array}\right)
\end{displaymath}
$ \mathfrak{G}=<e_{1}, e_{2}, e_{3}, e_{4}> $- four-dimensional
non-decomposable nilpotent Lie algebra with defined relations
$$
[e_{2},e_{3}]=e_{4}, [e_{3},e_{4}]=-e_{1}
$$
(this is $ g_{4,1} $ algebra in Mubaraczyanov classification\cite{mub2}).
\item Type III. $ \mathfrak{M} $ is a direct product of a two-dimensional
solvable Lie triple system $ <e_{1},e_{2}> $, and a one-dimensional Abelian
$ <e_{3}> $ :

\begin{displaymath}
A=(e_{1},e_{2},-)=\left(\begin{array}{ccc}
0 & \pm 1 & 1\\
0 & 0 & b\\
0 & 0 & 0\\
\end{array}\right),
B=0,
C=0
\end{displaymath}

$ \mathfrak{G}=<e_{1}, e_{2}, e_{3}, e_{4}> $  four-dimensional solvable and
decomposable Lie algebra, with defined relations:
$$
[e_{1},e_{2}]=e_{4}, [e_{2}, e_{4}]= \pm e_{1}
$$
moreover $ \mathfrak{G}=<e_{1},e_{2}, e_{4}> \oplus <e_{3}> $, where
$ <e_{1}, e_{2}, e_{4}> $- three-dimensional solvable Lie (algebra
 $ g_{3,4 \setminus 5}$ in Mubaraczyanov classification \cite{mub2}).
\item Type IV. \begin{displaymath}
A=(e_{1},e_{2},-)=\left(\begin{array}{ccc}
0 & \pm 1 & 1\\
0 & 0 & 0\\
0 & 0 & 0\\
\end{array}\right),
B=0,
C=(e_{3},e_{1},-)=\left(\begin{array}{ccc}
0 & -1 & \pm 1\\
0 & 0 & 0\\
0 & 0 & 0\\
\end{array}\right)
\end{displaymath},

$ \mathfrak{G}=<e_{1}, e_{2}, e_{3}, e_{4}> $- four-dimensional solvable and
non-decomposable Lie algebra, with defined relations:
$$
[e_{1},e_{2}]=e_{4}, [e_{2}, e_{4}]= \pm e_{1}
$$
$$
[e_{1},e_{3}]=\pm e_{4}, [e_{3}, e_{4}]= - e_{1}
$$
(algebra $ g_{4,5 \setminus 6}$ in Mubaraczyanov classification \cite{mub2}).

\item Type V
\begin{displaymath}
B=\left(\begin{array}{ccc}
0 &  1 & 0\\
0 & 0 & \pm 1\\
0 & 0 & 0\\
\end{array}\right),
A=C=0
\end{displaymath}

$ \mathfrak{G}=<e_{1}, e_{2}, e_{3}, e_{4}> $- four-dimensional solvable
non-decomposable Lie algebra with defined relations:
$$
[e_{2},e_{3}]=e_{4}, [e_{2}, e_{4}]= - e_{1}
$$
$$
[e_{3},e_{4}]=\mp e_{2}
$$
(algebra $ g_{8 \setminus 9}$ in Mubaraczyanov classification \cite{mub2}).

\item Type VI
\begin{displaymath}
A=0,
B=(e_{2},e_{3},-)=\left(\begin{array}{ccc}
0 & 0 & 1\\
0 & 0 & 0\\
0 & 0 & 0\\
\end{array}\right),
C=(e_{3},e_{1}, -)=\left(\begin{array}{ccc}
0 & 0 & 0\\
0 & 0 & 1\\
0 & 0 & 0\\
\end{array}\right)
\end{displaymath}

$ \mathfrak{G}=<e_{1}, e_{2}, e_{3}, e_{4},e_{5}> $- five-dimensional solvable
non-decomposable Lie algebra, with defined relations:
$$
[e_{1}, e_{2}]=e_{4}, [e_{1}, e_{3}]= - e_{5}
$$
$$
[e_{3}, e_{4}]=- e_{1}, [e_{3}, e_{5}]=-e_{2}
$$
(as a result we obtain an extension of four-dimensional Abelian ideal
 $ \mathfrak{G}=< e_{1}, e_{2}, e_{4}, e_{5}>$  by means of $ <e_{3}>$,
algebra $ g_{4, 13} $ in Mubaraczyanov classification \cite{mub2}).
\item Type VII.\begin{displaymath}
A=0,
B=(e_{2},e_{3},-)=\left(\begin{array}{ccc}
1 & 0 & 0\\
0 & 0 & 0\\
0 & 0 & 0\\
\end{array}\right),
C=(e_{3},e_{1}, -)=\left(\begin{array}{ccc}
0 & -1 & 0\\
0 & 0 & 0\\
0 & 0 & 0\\
\end{array}\right)
\end{displaymath}

\end{itemize}

Lie algebra $ \mathfrak{G}=<e_{1}, e_{2}, e_{3}, e_{4},e_{5}> $- five-
dimensional solvable non-decomposable Lie algebra, with defined relations:
$$
[e_{2}, e_{3}]=e_{4}, [e_{1}, e_{3}]=  e_{5}
$$
$$
[e_{1}, e_{4}]=- e_{1}, [e_{2}, e_{5}]=-e_{1}, [e_{4},e_{5}]=e_{5}
$$
(algebra $ g_{4, 11} $ in Mubaraczyanov classification \cite{mub2}).
\newpage
\subsection{EXAMPLES OF 3-DIMENSIONAL BOL ALGEBRAS WITH SOLVABLE TRILINEAR OPERATIONS}
\begin{itemize}
\item \underline{Example I.} Solvable Lie algebras

Let $ \mathfrak{L} $- be a 3-dimensional Lie algebra with basis $ e_1 $, $e_2 $
$e_3 $ and $ C^{i}_{jk}(C^{i}_{jk}=-C^{i}_{kj} $ as its structural constants.
The structural constant  $ C^{i}_{jk} $ can be represented as:
$$
  C^{i}_{jk}=\varepsilon_{ijk}b^{li}+\delta^{i}_{k}a_{k}-\delta^{i}_{j}a_{k}
$$

where $ \varepsilon_{ijk} $- discriminant tensor, $ b^{li} $ symmetrical tensor
such that $ (b^{li}=b^{il} $, $ \delta ^{i}_{k} $ - Kronecker symbol and $ a_j $
covector, defined such that $ b^{ij}a_{j}=0 $. Reducing $ b^{li} $ to the
diagonal view and, by examination of the following 3-dimensional Lie algebras
(called Bianchi classification \cite{dou}). Can be obtained the following non
isomorphic Lie algebras):
\begin{itemize}
\item I. $ [e_{i}, e_{j}=0$, $i$, $j$=1,2,3
\item II.$[e_{1},e_{2}]=0$, $[e_{1},e_{3}]=0$,$[e_{2},e_{3}]=e_{1}$
\item III.$[e_{1},e_{2}]=0$, $[e_{1},e_{3}]=e_{1}$,$[e_{2},e_{3}]=0$
\item IV.$[e_{1},e_{2}]=0$, $[e_{1},e_{3}]=e_{1}$,$[e_{2},e_{3}]=e_{1}+e_{2}$
\item V.$[e_{1},e_{2}]=0$, $[e_{1},e_{3}]=e_{1}$,$[e_{2},e_{3}]=e_{2}$
\item VI.$[e_{1},e_{2}]=0$, $[e_{1},e_{3}]=e_{1}$,$[e_{2},e_{3}]= \lambda e_{2}$
$\lambda \neq 0,1$
\item VII.$[e_{1},e_{2}]=0$, $[e_{1},e_{3}]=e_{2}$,$[e_{2},e_{3}]=-e_{1}\mu+e_{2}$  $(\mu^{2}<4)$
\item VIII.$[e_{1},e_{2}]=e_{1}$, $[e_{1},e_{3}]=2e_{2}$,$[e_{2},e_{3}]= e_{3}$
\item XI.$[e_{1},e_{2}]=e_{3}$, $[e_{1},e_{3}]=-e_{2}$,$[e_{2},e_{3}]= e_{1}$
\end{itemize}
Let's introduce in this examination the trilinear operation corresponding to
the Lie triple systems:
$$
(e_{i},e_{j},e_{k})=[[e_{i},e_{j}],e_{k}]=C^{t}_{ij} \cdot C^{s}_{tk}e_{s}
$$
and matrix
$$
A=(e_{1},e_{2},-), B=(e_{2},e_{3},-), C=(e_{1},e_{3},-).
$$
Applying it to the distinguished Lie algebras above (I-IX) we obtain
the following
(with the omition of the zero relation):
\begin{itemize}
\item I-II. $ A=B=C=0 $ (Abelian case).
\item III. $ (e_{1},e_{3},e_{3})=e_{1}$ (type II in Ch.II \S 3).
\item IV.$ (e_{1},e_{3},e_{3})=e_{1}$,$ (e_{2},e_{3},e_{3})=e_{2}$ (type VI in
Ch. II \S 3).
\item V. $(e_{1},e_{3},e_{3})=e_{1}$, $(e_{2},e_{3},e_{3})=e_{2}$ (subcase of
 type VI Ch. II \S 3).
\item VI. $ (e_{1},e_{3},e_{3})=e_{1}$, $(e_{2},e_{3},e_{3})=\lambda^{2}e_{2}$
$\lambda \neq 0,1 $.
\item VII. $(e_{1},e_{3},e_{3})=-e_{1}+ \mu e_{2}$,
$(e_{2},e_{3},e_{3})=- \mu e_{1}+ (1-\mu^{2}) e_{2}$ $(\mu^{2}<4)$ (also
allows imbedding in Type VI of classification in Ch. II \S 3).
\item $(e_{1},e_{3},e_{3})=-e_{1}$, $(e_{1},e_{3},e_{3})=2e_{3}$,
$(e_{1},e_{2},e_{3})=2e_{2}$, $(e_{1},e_{3},e_{1})=-2e_{1}$,
$(e_{2},e_{3},e_{1})=-2e_{2}$, $(e_{2},e_{3},e_{2})=-e_{3}$ (simple Lie triple
system).
\item IX. $(e_{1},e_{2},e_{2})=-e_{1}$,$(e_{1},e_{3},e_{1})=e_{3}$,
$(e_{1},e_{2},e_{3})=e_{2}$, $(e_{1},e_{3},e_{3})=-e_{1}$,
$(e_{2},e_{3},e_{2})=e_{3}$, $(e_{},e_{3},e_{3})=-e_{2}$ (simple Lie triple
system).
\end{itemize}

Choosing the solvable Lie triple system we obtain the following table:

\begin{tabular}{|r|c|l|}
\hline
  & Defining relation of Lie algebra &
  \parbox{40mm}{Corresponding type\\ of LTS}\\
\hline
I & Abelian & Abelian\\
\hline
II & $ [e_{2},e_{3}]=e_{1}$ & Abelian \\
\hline
III & $ [e_{1},e_{3}]=e_{1}$ & $<e_{1}> \oplus )$
\footnotemark
Type III \\
\hline
IV &  $[e_{2},e_{3}]=e_{1}+e_{2},\; [e_{1},e_{3}]=e_{1}$ &
   $\begin{array}{c} <e_{3}> \oplus <e_{1},e_{2}>\\ \mbox {Type VI}\end{array}$ \\
\hline
V & $ [e_{2},e_{3}]=e_{2}$,$[e_{1},e_{3}]=e_{1}$ & $ \mathbb{R} \oplus )$
  \footnotemark
  Type VI \\
\hline
VI & $ [e_{2},e_{3}]= \lambda e_{2}$ $\lambda \neq 0,1$,$[e_{1},e_{3}]=e_{1}$ & Type Vi \\
\hline
VII & $\begin{array}{c} [e_{2},e_{3}]=-e_{1}+ \mu e_{2}\\
              \mu^{2}<4$,$[e_{1},e_{3}]=e_{2} \end{array}$ &
 Type VI \\
\hline
\end{tabular}

\addtocounter{footnote}{-1}\footnotetext{2-dimensional solvable Lie triple system}
\addtocounter{footnote}{1}\footnotetext{2-dimensional solvable Lie triple system}

\item \underline{Example 2.} Classification of the right-alternative
3-dimensional algebras

A linear algebra $ V$ with the composition law $ \{x,y\} $ is called
right-alternative if $\forall x$, $y \in V$
$$
\{y,\{x,x\}\}=\{\{y,x,\}x\}
$$
In the work \cite{mike1} it was shown that, over any field of characteristic
 zero with isomorphism accuracy there exist only five types of 3-dimensional
right-alternative algebras. The indicated algebras are noted $A$, $B$, $C$, $E$, $H$.

Algebras  $A$, $B$, $C$, $E$ and $H$ are defined by the following relations:

\begin{itemize}
\item $A$:$\{e_{2},e_{3}\}=e_{2}$,$\{e_{3},e_{2}\}=e_{1}$, $\{e_{3},e_{3}\}=e_{3}$
\item  $B$:$\{e_{1},e_{3}\}=e_{1}$,$\{e_{3},e_{1}\}=e_{1}$, $\{e_{3},e_{2}\}=e_{1}+e_{2}$, $\{e_{3},e_{3}\}=e_{3} $
\item  $C$:$\{e_{1},e_{3}\}=e_{1}$,$\{e_{2},e_{3}\}=e_{2}$, $\{e_{3},e_{1}\}=e_{1}+e_{2}$, $\{e_{3},e_{3}\}=e_{3} $
\item  $E$:$\{e_{2},e_{2}\}=e_{1}$,$\{e_{2},e_{3}\}=e_{2}$, $\{e_{3},e_{2}\}= \beta e_{1}$, $\{e_{3},e_{3}\}=e_{3} $
\item  $H$:$\{e_{1},e_{3}\}=e_{1}$,$\{e_{2},e_{2}\}=e_{1}$, $\{e_{2},e_{3}\}= \gamma e_{1}$,   $\{e_{3},e_{3}\}=e_{3} $,$\{e_{3},e_{1}\}=e_{1}$, $\{e_{3},e_{2}\}=e_{2}$
\end{itemize}
for any arbitrary $ \beta $, $ \gamma $.

According to \cite{mike4}, any right-alternative algebra is a Bol algebra
under the operations:
$$
x \cdot y=\{x,y\}-\{y,x\}
$$
$$
<x,y,z>=\{\{x,y\},z\}-\{x,\{y,z\}\}
$$
$$
(x,y,z)=z \cdot (x \cdot y)+2<z,x,y>.
$$
In the case of the algebra $ A $ we obtain

$$ e_{2} \cdot e_{3}=e_{2}-e_{1}, \; (e_{2},e_{3},e_{3})=-e_{1}-e_{2}.$$

In the case of algebras $ B $ we obtain

$ e_{2} \cdot e_{3}=e_{2}+e_{1}$, $ (e_{2},e_{3},e_{3})=e_{1}-e_{2}$.

Here we see that the linear mapping with matrix

\begin{displaymath}
\left(\begin{array}{ccc}
-1 & 0 & 0\\
0 & 1 & 0\\
0 & 0 & 1\\
\end{array}\right)
\end{displaymath}

realizes the isomorphism of Bol algebra $A$ and $B$. The so obtained Bol
algebra has Type $ III^{+}$ in chapter III \S 3 . classifications.

The composition law, corresponding to the local analytical Bol loop, can be follows:

$$ \begin{pmatrix} x_1\\x_2\\x_3\end{pmatrix}\star
   \begin{pmatrix} y_1\\y_2\\y_3\end{pmatrix}
   =
   \begin{pmatrix}
     x_{1}+y_{1}+x_{2}y_{3}\\x_{2}+y_{2}+x_{3}y_{2}\\x_{3}+y_{3}+x_{3}y_{3}\end{pmatrix}
$$

In the case of algebra $C$

$$ e_{1} \cdot e_{3}=e_{2},\;\;  (e_{2},e_{3},e_{3})=-e_{2}$$
$$ e_{2} \cdot e_{3}=-e_{2},\; \;  (e_{1},e_{3},e_{3})=e_{2}.$$

In the base $ \widetilde{e_{1}}=e_{1}+e_{2}$, $ \widetilde{e_{2}}=e_{2}$,
 $ \widetilde{e_{3}}=e_{3}$ the defined relations turn to

$$ (\widetilde{e_{1}},\widetilde{e_{3}},\widetilde{e_{3}})=0, \; \;
\widetilde{e_{1}} \cdot \widetilde{e_{3}}=0 $$
$$ (\widetilde{e_{1}},\widetilde{e_{3}},\widetilde{e_{3}})=-\widetilde{e_{2}},\; \; \widetilde{e_{1}} \cdot \widetilde{e_{3}}=-\widetilde{e_{2}}. $$

Transforming the base again $ \bar{e_{1}}=\widetilde{e_{2}}$,
$ \bar{e_{2}}=\widetilde{e_{3}}$, $ \bar{e_{3}}=\widetilde{e_{2}}$, we obtain
an isomorphism of the given Bol algebra with algebra type $ III^{+}$ $(x=1)$,
see Chapter III \S 3.

The composition law corresponding to the local analytic Bol loop, can be
described as follows:

$$ \begin{pmatrix} x_1\\x_2\\x_3\end{pmatrix}\star
   \begin{pmatrix} y_1\\y_2\\y_3\end{pmatrix}
   =
   \begin{pmatrix}
     x_{1}+y_{1}+x_{1}y_{3}+x_{3}y_{1}\\x_{2}+y_{2}+x_{1}y_{3}+x_{3}y_{2}\\x_{3}+y_{3}+x_{3}y_{3}\end{pmatrix}
$$

In the case of algebra $ E $

$$ e_{2} \cdot e_{3}=e_{2}- \beta e_{1},\;  (e_{2},e_{3},e_{2})=-e_{1}$$
 $$ (e_{2},e_{3},e_{3})=-e_{2}- \beta e_{1}.$$

After transforming the base

$ \bar{e_{1}}=e_{1}$; $ \bar{e_{2}}=-e_{2}- \beta e_{1} $; $ \bar{e_{3}}=e_{3} $

we obtain

$(\bar{e_{2}},\bar{e_{3}}, \bar{e_{3}})=- \bar{e_{2}}$;
 $(\bar{e_{2}},\bar{e_{3}}, \bar{e_{3}})=- \bar{e_{1}}$;
$ \bar{e_{2}} \cdot \bar{e_{3}}=\bar{e_{2}}+2 \beta \bar{e_{1}}$.

In the case of algebra $ H $

$ e_{2} \cdot e_{3}=-e_{2}+ \gamma e_{1}$; $ (e_{2},e_{3},e_{2})=2e_{1}$;
  $ (e_{2},e_{3},e_{3})=-e_{2}+ \gamma e_{1}$.

After the first transforming of the base

$ \bar{e_{1}}=2e_{1}$; $ \bar{e_{2}}=e_{2} $; $ \bar{e_{3}}=e_{3}, $

we obtain

 $(\bar{e_{2}},\bar{e_{3}}, \bar{e_{3}})=\alpha \bar{e_{1}}- \bar{e_{2}}$;

 $(\bar{e_{2}},\bar{e_{3}}, \bar{e_{2}})= \bar{e_{1}}$;

$ \bar{e_{2}} \cdot \bar{e_{3}}=-\bar{e_{2}}+ \alpha \bar{e_{1}}$.

After the second transformation of the base

$ \widehat{e_{1}}= \bar{e_{1}}$,

 $ \widehat{e_{2}}=\alpha \bar{e_{1}}- \bar{e_{2}}$,

 $ \widehat{e_{3}}=-\bar{e_{3}}$, we obtain:

$ ( \widehat{e_{2}}, \widehat{e_{3}}, \widehat{e_{3}})=- \widehat{e_{2}}$,

$ ( \widehat{e_{2}}, \widehat{e_{3}}, \widehat{e_{2}})=- \widehat{e_{1}}$,

$ \widehat{e_{2}} \cdot \widehat{e_{3}}= \widehat{e_{2}}$.

The linear map with matrix

\begin{displaymath}
\left(\begin{array}{ccc}
-1 & 2 \beta & 0\\
0 & 1 & 0\\
0 & 0 & 1\\
\end{array}\right)
\end{displaymath}
for any $ \beta $ realize the isomorphism of $ H $ and $ E $. Then the
trilinear composition law has Type V in Chapter III \S 3. classification.
The composition law of local analytical Bol loop, can be described as follows:

$$ \begin{pmatrix} x_1\\x_2\\x_3\end{pmatrix}\star
   \begin{pmatrix} y_1\\y_2\\y_3\end{pmatrix}
   =
   \begin{pmatrix}
     x_{1}+y_{1}+x_{2}y_{2}+x_{3}y_{2}+x_{2}y_{3}\\x_{2}+y_{2}+x_{3}y_{2}\\x_{3}+y_{3}+x_{3}y_{3}\end{pmatrix}
$$

\item REMARK ABOUT THE CLASSIFICATION OF BOL 3-WEBS

In the work \cite{fedo} the classification of six-dimensional Bol 3-Webs was
examined. Let $ M^6$ be a smooth manifold given in the six-dimensional Bol
3-Webs $ W_6 $. Then three families of surfaces will define three perfectly
integrable Pfaffian system of equations: $ w^{i}_{1}=0 $,$ w^{i}_{2}=0 $
 $ w^{i}_{1}+ w^{i}_{2}=0 $, where $i$=1,2,3. The equations of structure of
3-Webs can be written as follows:
$$
dw^{i}_{1}=w^{j}_{1} \wedge \theta^{i}_{j}+\varepsilon_{jkl}a^{il}w^{j}_{1}\wedge w^{k}_{1}
$$
$$
dw^{i}_{2}=w^{j}_{2} \wedge \theta^{i}_{j}-\varepsilon_{jkl}a^{il}w^{j}_{2}\wedge w^{k}_{2}
$$
$$
d \theta^{i}_{j}=\theta^{k}_{j} \wedge \theta^{i}_{k}+\varepsilon_{klm}b^{im}w^{k}_{1}\wedge w^{l}_{2}
$$
$$
a^{ij}=a^{ij} \theta^{p}_{p}+ \frac{1}{2} (b^{ij}_{m}-b^{ip}_{p}\delta^{i}_{m})(w^{m}_{1}-w^{m}_{2})
$$

Moreover the tensor $ a^{ij} $, $( i,j=1,2,3)$ can be represented as
$ a^{ij}=b^{ij}+c^{ij} $, where $ b^{ij}=a^{(ij)} $- symmetrical tensor,
$ c^{ij}=a^{[ij]} $- skew-symmetrical tensor. As the tensor $ c^{ij} $ is
skew-symmetrical, then its rank can only be equal to zero or two.

Let the rang of $ c^{ij} $ be equal to two, then after a suitable
transformation and normalizing its frame, the tensor $ a^{ij} $ can be reduced to one of the following types:

\begin{displaymath}
a^{ij}=\left(\begin{array}{ccc}
\epsilon_1 & 0 & 0\\
0 & \epsilon & 1\\
-1 & -1 & 0\\
\end{array}\right)
\end{displaymath}
where $ \epsilon =\pm 1$; \; $\epsilon_1=\pm 1;0$
\begin{displaymath}
a^{ij}=\left(\begin{array}{ccc}
\epsilon & 0 & 0\\
0 & a & 1\\
0 & -1 & \pm a\\
\end{array}\right)
\end{displaymath}
where $ \epsilon =\pm 1,0$ and for any $a$,
\begin{displaymath}
a^{ij}=\left(\begin{array}{ccc}
\epsilon & 0 & 0\\
0 & 0 & 1\\
0 & -1 & \epsilon_{1}\\
\end{array}\right)
\end{displaymath}
where $ \epsilon =\pm 1,0 $ $ \epsilon_{1}=\pm 1,0 $.

In this connection, only the last one verifies the condition of compatibility
of the equations (1) \cite{fedo}.
$$
b^{ik}_{k}=2\varepsilon_{jkl}a^{il}a^{jk}
$$
$$
b^{ij}_{p}a^{pk}-b^{jk}_{p}a^{ip}-b^{ik}_{p}a^{pj}=b^{ip}_{p}a^{jk}-b^{pk}_{p}a^{ij}
$$
which represents equations (1) ordinarily as
$$
a^{i}_{jk}=\varepsilon_{jkl}a^{il}, \; \;b^{i}_{jkl}=\varepsilon_{klp}b^{ip}_{j}
$$
where $ \varepsilon_{jkl} $, $ \varepsilon_{klp} $ are discriminant tensors.
As shown in \cite{fedo} the skew-symmetric tensor $ a^{ij}$

\begin{displaymath}
a^{ij}=\left(\begin{array}{ccc}
0 & 0 & 0\\
0 & 0 & 1\\
0 & -1 &0\\
\end{array}\right)
\end{displaymath}

is the necessary and sufficient condition for six-dimensional isocline Bol Web.

In other words, if the tensor $ a^{ij} $ (six-dimensional Bol three-Web) is
symmetric and the the rank of $ a^{ij} $ equals zero, the Web in this case is parallelizable.

The trilinear operations of Bol algebras corresponding to the given three-Webs,
with tensor $ a^{ij} $ of type above, define in\cite{fedo}, prove to be
simple and splittable(in the sense of the terminology of Chapter I \S 5.).
\end{itemize}
\newpage
\section{CHAPTER III}
\begin{center}
CLASSIFICATION OF 3-DIMENSIONAL BOL ALGEBRAS AND CALCULUS OF CORRESPONDING BOL 3-WEBS
\end{center}
\subsection{BOL ALGEBRAS WITH TRIVIAL TRILINEAR OPERATIONS OF TYPE I}

Let $ \mathfrak{M} $ be a Bol algebra [46] with a trivial trilinear operation.
The structure of $ \mathfrak{M} $ will be reduced, to the representation of an
anti-commutative bilinear multiplication
$( \cdot ) $$ \mathfrak{M} \times \mathfrak{M} \longrightarrow \mathfrak{M} $,
such that:
$$
(X \cdot Y) \cdot (Z \cdot U)=O
$$
$ \forall  X, Y,Z,U \in  \mathfrak{M}$

Possible case:
\begin{enumerate}
\item $ \xi \cdot \eta =0 $ (abelian case)
\item $ \xi \cdot \eta \zeta =0 $ the algebra $ \mathfrak{M}$ - not abelian
(but anticommutative 2-nilpotent algebra),
\item $ ( \xi \eta ) \cdot ( \zeta \kappa)=0 $, but the algebra $ \mathfrak{M}
( \cdot )$ not 2-nilpotent algebra).
\end{enumerate}

 (In particular, $ \mathfrak{M}$ $( \cdot )$ can be anticommutative 3-nilpotent
algebra) hence:
$$
( \xi \eta \cdot \zeta) \cdot \kappa=0
$$
$ \forall \xi, \eta, \zeta, \kappa $ $\in  \mathfrak{M} $

  In this investigation we are considering the question {\bf how long will ( with
accuracy to the isomorphism) exist the operation viewed in (1)-(3) in the case
of 3-dimensional algebra $ \mathfrak{M}$}.

  We will consider (with accuracy to the isomorphism) the uniquness of the abelian algebra, denoted Type III.1 and hence investigate case 2.

Case 2. Let's denote the subspace $ \mathfrak{M} \cdot \mathfrak{M} $ in
$ \mathfrak{M} $ by $ V$ then $ V \cdot \mathfrak{M}=0 $ and $ V  \neq \mathfrak{M} $, $ V \neq 0$. The following variants can be obtained:
\begin{itemize}
\item 2.a $ dim V=1$, then  consider that,
$ \mathfrak{M}$=$<e_{1}, e_{2}, e_{3}>$, $ V=<e_{1}>$,
$$
e_{1} \cdot e_{2}=e_{1} \cdot e_{3}=0
$$
$$
e_{2} \cdot e_{3}=\alpha e_{1}, \alpha \neq 0
$$
by adjusting the base $e_{1}$,$ e_{2}$, $e_{3}$, we  make
$$
e_{2} \cdot e_{3}= e_{1},
$$
the so obtained algebra will be denoted by algebra Type III.2.
\item 2.b. $ dim V=2 $, then $ \mathfrak{M}=<e_{1}, e_{2}, e_{3}>$,
$ V=<e_{1}, e_{2}>$, and $ e_{1} \cdot e_{2}=e_{1} \cdot e_{3}=0$,
$ e_{2} \cdot e_{3}=0$, \; \; \;(1)

The relation in (1) are in contradiction with the condition
$ \mathfrak{M} \cdot  \mathfrak{M}=V$. Hence there exist with isomorphism
accuracy only one Bol algebra of case 2.
\end{itemize}

  Now we past to the examination of Bol algebra case 3.

  Let us suppose that $ \mathfrak{M} \cdot \mathfrak{M}=V \subseteq \mathfrak{M}$, thus:
\begin{center}
$ V \cdot V = \{O\}$ but $ V  \cdot  \mathfrak{M} \neq  \{O\}$
 $ V \neq  \mathfrak{M}$, $V \neq \{O\} $.
\end{center}

\begin{itemize}
\item 3.a Let $ dim V =1$ and $ V =<e_{1}>$, where $e_{1},e_{2}, e_{3}$- are
basis vectors in $ \mathfrak{M}$, then:

$$
e_{1}\cdot  e_{2}=\alpha e_{1}, e_{1} \cdot e_{3}=\beta e_{1},
$$
$$
e_{2} \cdot e_{3}= \gamma e_{1}, \alpha^{2}+ \beta^{2} \neq 0.
$$
Without limiting our selves, one can consider $ \beta \neq 0$, by changing
the base  $e_{1},e_{2}, e_{3}$ the defined relation of anti-commutative algebra
$ \mathfrak{M}$ can be reduced to:
$$
e_{1}\cdot  e_{2}=0, e_{1} \cdot e_{3}= e_{1},
$$
$$
e_{2} \cdot e_{3}= \varepsilon e_{1}, where \; \;\varepsilon=1, 0.
$$

Let us denote the defined algebra through $\mathfrak{M}_{0}$ and
$\mathfrak{M}_{1}$ correspondently, and examine its isomorphism. We will note:
$$
 \mathfrak{M}_{0}=\{O\} \;  \;  \mathfrak{M}_{1}=\{e_{1}\},
$$
that means their centers are different, hence algebras  $\mathfrak{M}_{0}$ and
$\mathfrak{M}_{1}$  are not isomorphic. We will denote them by, Type III.3 and
Type III.4 correspondently.
\item 3.b. Let $ dim V=2 $, and let's consider  $e_{1},e_{2}, e_{3}$ the basis
in $ \mathfrak{M}$ such that $ V=< e_{1},e_{2}>$, thus:
$$
e_{1}\cdot  e_{2}=0
$$
and
$$
e_{1} \cdot e_{3}, e_{2} \cdot e_{3} \in V
$$
Deforming $e_{3}$ to $ \widetilde{e_{3}}=te_{3}+ve_{2}+ue_{1}$, where $t\neq 0$
$ u, v \in \mathbb{R} $ we obtain
$$
e_{1} \cdot \widetilde{e_{3}}=e_{1} \cdot (te_{3}+ve_{2}+ue_{1})=te_{1} \cdot e_{3},
$$
$$
e_{2} \cdot \widetilde{e_{3}}=e_{2} \cdot (te_{3}+ve_{2}+ue_{1})=te_{2} \cdot e_{3},
$$

we can consider that, the transformation

$ \Phi=ade_{3}|: V \longrightarrow V$
\begin{displaymath}
\Phi=\left(\begin{array}{cc}
\alpha & \beta \\
\gamma & \delta \\
\end{array}\right)
\end{displaymath}
is defined with accuracy to the multiplicative scalar.

The choice of the basis $ e_{1}, e_{2} $ in $V$ and vector $e_{3}$ we can
reduce $ \phi$ to one of the following:
\begin{displaymath}
\left(\begin{array}{cc}
1 & 0 \\
0 & \mu \\
\end{array}\right)
\end{displaymath}
or
\begin{displaymath}
\left(\begin{array}{cc}
1 & 0 \\
\alpha & 1 \\
\end{array}\right)
\end{displaymath}
where $\alpha \neq 0$ and $ \forall \mu \in \mathbb{R}$

Respectively
$$
e_{1} \cdot e_{3}=e_{1},
$$
$$
e_{2} \cdot e_{3}= \mu e_{2}
$$
or
$$
e_{1} \cdot e_{3}=e_{1}+ \alpha e_{2}
$$
$$
e_{2} \cdot e_{3}=e_{2}
$$

Let us examine the first case: In that case $ \mu \neq 0 $ (otherwise we will
obtain the contradiction with $ dim V=2$), by adjusting $ e_{2}$ to
$ \widetilde{e_{2}}= \mu e_{2} $ we obtain:

$$
e_{1} \cdot e_{3}=e_{1},
$$
$$
e_{2} \cdot e_{3}= e_{2}
$$
(algebra Type III.5)

In the second case the change of $ (1 \setminus a) e_{1} \longrightarrow e_{1}$
reduces the defined relations of algebra $ \mathfrak{M}$ to:

$$
e_{1} \cdot e_{3}=e_{1}+  e_{2}
$$
$$
e_{2} \cdot e_{3}=e_{2}
$$
(algebra Type III.6).
The algebra $ \mathfrak{M} $ of Type III.5 and Type III.6 are obtained by the
extension 1-dimensional Abelian algebra by means of 2-dimensional Abelian
algebra $ V= \mathfrak{M} \cdot \mathfrak{M}=<e_{1},e_{2}>$, hence for
$ a \in \mathfrak{M} \setminus V $ the structure of the operator $ ad | V $ is:

\begin{displaymath}
\left(\begin{array}{cc}
\alpha & 0 \\
0 & \alpha \\
\end{array}\right)
\end{displaymath}
 $\alpha \neq 0$ (algebras Type III.5)
\begin{displaymath}
\left(\begin{array}{cc}
\alpha & 0 \\
\alpha & \alpha \\
\end{array}\right)
\end{displaymath}
 $\alpha \neq 0 $(algebras Type III.6)

that's why algebra of Type III.5 and Type III.6 are not isomorphic.

\end{itemize}

In result we obtain the following theorem:

\underline{{\bf Theorem III.1}} There exits with isomorphism accuracy 6 Bol
algebras, with
a trivial trilinear operation defined as follows:
\begin{itemize}
\item III.1 trivial bilinear operation
\item III.2 $ e_{2} \cdot e_{3}=e_{1} $
\item III.3 $ e_{1} \cdot e_{3}=e_{1} $
\item III.4 $ e_{2} \cdot e_{3}=e_{1} $, $ e_{1} \cdot e_{3}=e_{1} $
\item III.5 $ e_{2} \cdot e_{3}=e_{2} $, $ e_{1} \cdot e_{3}=e_{1} $
\item III.6 $ e_{2} \cdot e_{3}=e_{2} $, $ e_{1} \cdot e_{3}=e_{1}+e_{2} $
\end{itemize}

\; \;We will note the enveloping Lie algebras of Bol algebras III.3 and Bol algebras
III.4 are 4-dimensional, but these Bol algebras are not isotopics by their
structure. The enveloping Lie algebras of Bol algebras III.5 and Bol algebras
III.6 are identical, but these Bol algebras are not isotopic by definition.

\; \; Below we will give the description of 3-webs, corresponding to the
 defined Bol algebras.

\begin{itemize}
\item Type III.1. Bol algebra $\mathfrak{B}$ with trivial trilinear and bilinear operation. Here we obtain grouped 3-Webs (Abelian group$<\mathbb{R}^{3},+,0>$
).
\item Type III.2. Bol algebras $\mathfrak{B}$ with bilinear anti-commutative
operation of view:
$$
e_{2} \cdot e_{3}=e_{1}.
$$

We also obtain here a grouped 3-Webs (global 3-Web) corresponding to the Lie
group, which is isomorphic to the upper triangular unipotent matrix, which
means matrix of form:

\begin{displaymath}
\left(\begin{array}{ccc}
1 & x & y\\
0 & 1 & z\\
0 & 0 & 1\\
\end{array}\right)
\end{displaymath}

\item Type III.3 Bol algebras $ \mathfrak{B} $ with bilinear operation of View:

$$
e_{1} \cdot e_{3}=e_{1}.
$$

and trivial trilinear operation, has a 4-dimensional canonical enveloping Lie
algebra $ \mathfrak{G}=<e_{1}, e_{2},e_{3}, e_{4}>$,
$$
\mathfrak{G}=\mathfrak{B} \dotplus \mathfrak{h}, \mathfrak{B}=<e_{1}, e_{2},e_{3}>,  \mathfrak{h}=<e_{4}-e_{4}>
$$

with a composition law having the following constant of structure
$$
[e_{1}, e_{3}]= e_{4}
$$

the composition law $ (\triangle)$ corresponding to the Lie group $ G$ is
defined as follows:

$$ \begin{bmatrix} x_1\\x_2\\x_3\\x_4\end{bmatrix} \triangle
   \begin{bmatrix} y_1\\y_2\\y_3\\y_4\end{bmatrix}
   =
   \begin{bmatrix}
     x_{1}+y_{1}\\x_{2}+y_{2}\\x_{3}+y_{3}\\x_{4}+y_{4}+\frac{x_{1}y_{3}-y_{1}x_{3}}{2}\end{bmatrix}
$$

Moreover the subgroup $ H= \exp \mathfrak{h}$, realized as a collection of
elements:
$$
\{\exp t(e_{4}-e_{1})\}_{t \in \mathbb{R}}=\{(\exp te_{4})(exp te_{1})^{-1}\}_{t \in \mathbb{R}}=\{-t,0,0,t\}_{t \in \mathbb{R}}.
$$

The collection of elements
$$
 B=\exp \mathfrak{B}=\{\exp (te_{1}+ue_{2}) \cdot \exp ve_{3}\}_{t,u,v \in \mathbb{R}}=\{t,u,v,t\}_{t,u,v \in \mathbb{R}}
$$
form a local section of space left coset $ G \bmod H $.

The subgroup $ H $ is defined as:
$$
H=\exp \mathfrak{h}=\{\exp \alpha (e_{4}-e_{3})\}_{\alpha \in \mathbb{R}}=\{0,0,-\alpha,\alpha\}_{\alpha \in \mathbb{R}},
$$

$ B=\{t, u, v, 0 \}. $

Any element $ ( x_{1}, x_{2}, x_{3}, x_{4})$ from $ G $ such that \mbox{$ ( X_{3}<-2)$}
can be uniquely represented in the form:

$$ \begin{pmatrix} x_1\\x_2\\x_3\\x_4\end{pmatrix} =
   \begin{bmatrix}\frac{2x_1+2x_4-x_1x_3}{2-x_3}\\x_2\\x_3\\0\end{bmatrix}
   \triangle
   \begin{bmatrix}
     -\frac{2x_4}{2-x_3}\\0\\0\\\frac{2x_4}{2-x_3}\end{bmatrix}
$$
$$
=\prod _{B}\begin{bmatrix} x_1\\x_2\\x_3\\x_4\end{bmatrix}
   \triangle
   \begin{bmatrix}
     -\frac{2x_4}{2-x_3}\\0\\0\\\frac{2x_4}{2-x_3}\end{bmatrix}.
$$

\; \; The composition law$( \star )$ of local analytic Bol loop $B( \star )$ is
defined as:

$$ \begin{pmatrix} t\\u\\v\end{pmatrix} \star
   \begin{pmatrix}t'\\u'\\v'\end{pmatrix}
   = \prod_{B}  \left(
    \begin{bmatrix}t\\u\\v\\0\end{bmatrix}\triangle \begin{bmatrix}
     t'\\u'\\v'\\0\end{bmatrix} \right)
$$
$$
= \begin{bmatrix}t+t'+\frac{tv'-vt'}{2-(v+v')}\\u+u'\\v+v'\end{bmatrix}.
$$

The corresponding local Bol 3-Webs can be realized in the neighborhood of the
point $ (O,O)$ in $ \mathbb{R}^{6}=\{(\bar{X},\bar{Y}), \bar{X},\bar{Y} \in\mathbb{R}^3\}$ as a space form by:
$$
\bar{X}=const
$$
$$
\bar{Y}=const
$$
$$
\bar{X} \star \bar{Y}=const.
$$
\item Type III.4 Bol algebra $ \mathfrak{B}$ with bilinear operation of view:
$$
e_{1} \cdot e_{3}=e_{1}, e_{2} \cdot e_{3}=e_{1}
$$

an trivial trilinear operation, has a 4-dimensional canonical enveloping Lie
algebra $\mathfrak{G}=<e_1, e_2, e_3, e_4>,$
$\mathfrak{G}=\mathfrak{B} \dotplus \mathfrak{h}$, $ \mathfrak{B}= <e_1, e_2, e_3>$,$ \mathfrak{h}=<e_{4}-e_{1}>$

With composition law having the following structural equations
$$
[e_1, e_3]=e_4, \; [e_2, e_3]=e_4.
$$

\; \; The composition law $( \triangle)$, corresponding to the Lie group $ G $ is defined as follows:

$$ \begin{bmatrix} x_1\\x_2\\x_3\\x_4\end{bmatrix} \triangle
   \begin{bmatrix} y_1\\y_2\\y_3\\y_4\end{bmatrix}
   =
   \begin{bmatrix}
     x_{1}+y_{1}\\x_{2}+y_{2}\\x_{3}+y_{3}\\x_{4}+y_{4}+\frac{(x_{1}+x_{2})y_{3}-(y_{1}+y_{2})x_{3}}{2}\end{bmatrix}
$$

Moreover the subgroup $ H=\exp \mathfrak{h} $, is realized as a collection of elements:

$$
\{\exp t(e_{4}-e_{1})\}_{t \in \mathbb{R}}=\{(\exp te_{4})\cdot (\exp te_{1})^{-1}\}_{t \in \mathbb{R}}=\{-t,0,0,t\}_{t \in \mathbb{R}}.
$$

\; \; The collection

$$
 B=\exp \mathfrak{B}=\{\exp (te_{1}+ue_{2}) \cdot \exp ve_{3}\}_{t,u,v \in \mathbb{R}}=\{t,u,v,0\}_{t,u,v \in \mathbb{R}}
$$

form a local section of space left coset $ G \bmod H $.

\; \; Any element $(x_1, x_2, x_3, x_4 )$ from $G$ such that $(x_3>-2)$, can be
uniquely represented in the form:

$$ \begin{pmatrix} x_1\\x_2\\x_3\\x_4\end{pmatrix} =
   \begin{bmatrix}\frac{2x_1+2x_4+x_1x_3}{2+x_3}\\x_2\\x_3\\0\end{bmatrix}
   \triangle
   \begin{bmatrix}
     -\frac{2x_4}{2+x_3}\\0\\0\\\frac{2x_4}{2+x_3}\end{bmatrix}
$$
$$
=\prod _{B}\begin{bmatrix} x_1\\x_2\\x_3\\x_4\end{bmatrix}
   \triangle
   \begin{bmatrix}
     -\frac{2x_4}{2+x_3}\\0\\0\\\frac{2x_4}{2+x_3}\end{bmatrix}.
$$

\; \; The composition law( $\star$ ), of local analytic Bol loop $B( \star )$ is
defined as:

$$ \begin{pmatrix} t\\u\\v\end{pmatrix} \star
   \begin{pmatrix}t'\\u'\\v'\end{pmatrix}
   = \prod_{B}\left(
    \begin{bmatrix}t'\\u'\\v'\\0\end{bmatrix}\triangle \begin{bmatrix}
     t'\\u'\\v'\\0\end{bmatrix} \right)
$$
$$
= \begin{bmatrix}t+t'+\frac{(t+u)v'-v(t'+u')}{2+(v+v')}\\u+u'\\v+v'\end{bmatrix}.
$$

$$
= \begin{bmatrix}\frac{2t+2t'+tv+tv'+t'v'+t'v+uv'-vu'}{2+(v+v')}\\u+u'\\v+v'\end{bmatrix}.
$$

\; \; The corresponding local Bol 3-Web, can be realized in the neighborhood of
 the point $(O,O)$ in  $ \mathbb{R}^{6}=\{(\bar{X},\bar{Y}), \bar{X},\bar{Y} \in\mathbb{R}^3\}$ as a space  of second order.

\item Type III.5 \; Bol algebra $ \mathfrak{B}$, with bilinear operation of view:
$$
e_{1} \cdot e_{3}=e_{1},\; \; e_{2} \cdot e_{3}=e_{2}
$$

an trivial trilinear operation, has a 5-dimensional canonical enveloping Lie
algebra $\mathfrak{G}=<e_1, e_2, e_3, e_4, e_5>,$
$\mathfrak{G}=\mathfrak{B} \dotplus \mathfrak{h}$, $ \mathfrak{B}=<e_1, e_2, e_3>$,$ \mathfrak{h}=<e_{4}-e_{1},e_{5}-e_{2}>$

With composition law having the following structural equations
$$
[e_1, e_3]=e_4, \; [e_2, e_3]=e_5.
$$

\; \; The composition law $( \triangle)$ corresponding to the Lie group $ G $ is defined as follows:

$$ \begin{bmatrix} x_1\\x_2\\x_3\\x_4\\x_5\end{bmatrix} \triangle
   \begin{bmatrix} y_1\\y_2\\y_3\\y_4\\y_5\end{bmatrix}
   =
   \begin{bmatrix}
     x_{1}+y_{1}\\x_{2}+y_{2}\\x_{3}+y_{3}\\x_{4}+y_{4}+\frac{x_{1}y_{3}-y_{1}x_{3}}{2}\\x_{5}+y_{5}+\frac{x_{2}y_{3}-x_{3}y_{2}}{2}\end{bmatrix}
$$

Moreover the subgroup $ H=\exp \mathfrak{h} $, can be realized as a collection of elements:

$$
\{\exp t(e_{4}-e_{1}), \exp p(e_{5}-e_{2}\}_{t,p \in \mathbb{R}}=\{-t,-p,0,t,p\}_{t,p \in \mathbb{R}}.
$$

\; \; The collection

$$
 B=\exp \mathfrak{B}=\{\exp (te_{1}+ue_{2}) \cdot \exp ve_{3}\}_{t,u,v \in \mathbb{R}}=\{t,u,v,0,0\}_{t,u,v \in \mathbb{R}}
$$

form a local section of space left coset $ G \bmod H $.

\; \; Any element $(x_1, x_2, x_3, x_4,x_5 )$ from $G$ such that $(x_3>-2)$, can be
uniquely represented in the form:

$$ \begin{pmatrix} x_1\\x_2\\x_3\\x_4\\x_5\end{pmatrix} =
   \begin{bmatrix}\frac{2x_1+2x_4+x_1x_3}{2+x_3}\\\frac{2x_2+2x_5+x_2x_3}{2+x_3}\\x_3\\0\\0\end{bmatrix} \triangle  \begin{bmatrix} -\frac{2x_4}{2+x_3}\\-\frac{2x_5}{2+x_3}\\0\\\frac{2x_4}{2+x_3}\\\frac{2x_5}{2+x_3}\end{bmatrix}
$$
$$
=\prod _{B}\begin{bmatrix} x_1\\x_2\\x_3\\x_4\\x_5\end{bmatrix}
   \triangle
   \begin{bmatrix}
     -\frac{2x_4}{2+x_3}\\-\frac{2x_5}{2+x_3}\\0\\\frac{2x_4}{2+x_3}\\\frac{2x_5}{2+x_3}\end{bmatrix}.
$$

\; \; The composition law($ \star$ ), of local analytic Bol loop $B( \star )$ is
defined as:

$$ \begin{pmatrix} t\\u\\v\end{pmatrix} \star
   \begin{pmatrix}t'\\u'\\v'\end{pmatrix}
   = \prod_{B}\left(
    \begin{bmatrix}t\\u\\v\\0\end{bmatrix}\triangle \begin{bmatrix}
     t'\\u'\\v'\\0\end{bmatrix}\right)
$$
$$
= \begin{bmatrix}t+t'+\frac{tv'-vt'}{2+(v+v')}\\u+u'+\frac{uv'-vu'}{2+(v+v')}\\v+v'\end{bmatrix}.
$$

$$
= \begin{bmatrix}\frac{2t+2t'+tv+2tv'+t'v'}{2+(v+v')}\\\frac{2u+2u'+uv+2uv'+u'v'}{2+(v+v')}\\v+v'\end{bmatrix}.
$$

\; \; The corresponding local Bol 3-Web, can be realized in the neighborhood of
 the point $(O,O)$ in  $ \mathbb{R}^{6}=\{(\bar{X},\bar{Y}), \bar{X},\bar{Y} \in\mathbb{R}^3\}$ as a space of second order.

\item Type III.6 \; Bol algebra $ \mathfrak{B}$ with bilinear operation of view:
$$
e_{1} \cdot e_{3}=e_{1}+e_{2},\; e_{2} \cdot e_{3}=e_{2}
$$

an trivial trilinear operation, has a 5-dimensional canonical enveloping Lie
algebra $\mathfrak{G}=<e_1, e_2, e_3, e_4, e_5>,$
$\mathfrak{G}=\mathfrak{B} \dotplus \mathfrak{h}$, $ \mathfrak{B}=<e_1, e_2, e_3>$,$ \mathfrak{h}=<e_{4}-e_{1}-e_{2},e_{5}-e_{2}>$.

With composition law having the following structural equations
$$
[e_1, e_3]=e_4, \; [e_2, e_3]=e_5.
$$

\; \; The composition law $( \triangle)$, corresponding to the Lie group $ G $ is defined as follows:

$$ \begin{bmatrix} x_1\\x_2\\x_3\\x_4\\x_5\end{bmatrix} \triangle
   \begin{bmatrix} y_1\\y_2\\y_3\\y_4\\y_5\end{bmatrix}
   =
   \begin{bmatrix}
     x_{1}+y_{1}\\x_{2}+y_{2}\\x_{3}+y_{3}\\x_{4}+y_{4}+\frac{x_{1}y_{3}-y_{1}x_{3}}{2}\\x_{5}+y_{5}+\frac{x_{2}y_{3}-x_{3}y_{2}}{2}\end{bmatrix}.
$$

Moreover the subgroup $ H=\exp \mathfrak{h} $, can be realized as a collection of elements:

$$
\{\exp t(e_{4}-e_{1}-e{2}), \exp p(e_{5}-e_{2}\}_{t,p \in \mathbb{R}}=\{-t,-t-p,0,t,p\}_{t,p \in \mathbb{R}}.
$$

\; \; The collection

$$
 B=\exp \mathfrak{B}=\{\exp (te_{1}+ue_{2}) \cdot \exp ve_{3}\}_{t,u,v \in \mathbb{R}}=\{t,u,v,0,0\}_{t,u,v \in \mathbb{R}}
$$

form a local section of space left coset $ G \bmod H $.

\; \; Any element $(x_1, x _2, x_3, x_4,x_5 )$ from $G$ such that $(x_3>-2)$, can be
uniquely represented in the form:

$$ \begin{pmatrix} x_1\\x_2\\x_3\\x_4\\x_5\end{pmatrix} =
   \begin{bmatrix}\frac{2x_1+2x_4+x_1x_3}{2+x_3}\\\frac{x_2 (2+x_3)^2+2x_4(2+x_3)+4x_5+2x_3 x_5-2x_3 x_4}{(2+x_3)^2}\\x_3\\0\\0\end{bmatrix}
   \triangle
   \begin{bmatrix}
     -\frac{2x_4}{2+x_3}\\-\frac{4x_4+4x_5+2x_3 x_5}{(2+x_3)^2}\\0\\\frac{2x_4}{2+x_3}\\\frac{4x_5+2x_3 x_5-2x_3 x_4}{(2+x_3)^2}\end{bmatrix}
$$
$$
=\prod _{B}\begin{bmatrix} x_1\\x_2\\x_3\\x_4\\x_5\end{bmatrix}
   \triangle
   \begin{bmatrix}

     -\frac{2x_4}{2+x_3}\\-\frac{4x_4+4x_5+2x_3 x_5}{(2+x_3)^2}\\0\\\frac{2x_4}{2+x_3}\\\frac{4x_5+2x_3 x_5-2x_3 x_4}{(2+x_3)^2}\end{bmatrix}.
$$

\; \; The composition law $( \star )$, of local analytic Bol loop $B( \star )$ is
defined as:

$$ \begin{pmatrix} t\\u\\v\end{pmatrix} \star
   \begin{pmatrix}t'\\u'\\v'\end{pmatrix}
   = \prod_{B}\left(
    \begin{bmatrix}t\\u\\v\\0\end{bmatrix}\triangle \begin{bmatrix}
     t'\\u'\\v'\\0\end{bmatrix}\right)
$$
$$
= \begin{bmatrix}t+t'+\frac{tv'-vt'}{2+(v+v')}\\
u+u'+\frac{tv'-vt'}{2+(v+v')}+\frac{uv'-vu'}{2+(v+v')}-\frac{(v+v')(tv'-vt')}{(2+(v+v'))^2}\\v+v'\end{bmatrix}.
$$

\; \; The corresponding local Bol 3-Web can be realized in the neighborhood of
 the point $(O,O)$ in  $ \mathbb{R}^{6}=\{(\bar{X},\bar{Y}), \bar{X},\bar{Y} \in\mathbb{R}^3\}$ as a space of third order.
\end{itemize}
\newpage

\subsection{BOL ALGEBRAS WITH TRILINEAR OPERATION OF TYPE II}

\;We will base our investigation of 3-dimensional Bol algebras on the
examination
of their canonical enveloping Lie algebras. In what follows, we consider Bol
algebras of dimension 3, from their construction see \cite{mike6,sabmik1};
it follows that
the dimension of their canonical enveloping Lie algebras can not be more than 6
. Below we limit ourselves to the classification of Bol algebras
(and their corresponding 3-Webs) with canonical enveloping Lie
algebras of dimension$ \leq 5$.

\; \; \; Let $\mathfrak{B}$ be a 3-dimensional Bol algebras with a trilinear
operation Type II in \cite{boue2}, and
$ \mathfrak{G}=\mathfrak{B} \dotplus \mathfrak{h} $ - its canonical enveloping
Lie  algebra according to \cite{boue2}, we note that the situation
$ dim \mathfrak{G}=3$ is impossible that means the case of Web-group
is excluded.

\; \; Let us examine the case $ dim \mathfrak{G}=4 $, the structural constants
of Lie algebra $ \mathfrak{G}=<e_1,e_2, e_3, e_4>$,
$ \mathfrak{B}=<e_1,e_2, e_3>$ are defined as follows:
$$
[e_2,e_3]=e_4, \; \; [e_3,e_4]=-e_1
$$
here $\mathfrak{G}=\mathfrak{B} \dotplus [ \mathfrak{B}, \mathfrak{B}]$,
$[ \mathfrak{B}, \mathfrak{B}]=<e_4>$.

Introducing the subalgebra
$$
\mathfrak{h}_{ \alpha, \beta, \gamma}=<e_4+ \alpha e_1+ \beta e_2+ \gamma e_3>
$$
were $ \alpha, \beta, \gamma \in \mathbb{R} $.

We are getting a collection of Bol algebras with the structural constants
defined as follows:
$$
e_2 \cdot e_3 =- \alpha e_1 - \beta e_2 -\gamma e_3, \; \; (e_2,e_3,e_3)=e_1.
$$

\; \; Below we give an isomorphical and an isotopical classification of this
collection.

\; \; The group of automorphism $ F$ of Lie triple system $ \mathfrak{B} $ is
defined as follows:

****

\begin{displaymath}
ad ( \xi )= \left(\begin{array}{cccc}
0 & 0 & 0 & -z \\
0 & 0 & 0 & 0  \\
0 & 0 & 0 & 0 \\
0 & -z & y & 0 \\
\end{array}\right)
\end{displaymath}

\begin{displaymath}
Ad ( \xi )= \left(\begin{array}{cccc}
1 & \frac{z^2}{2} &- \frac{zy}{2} & -z \\
0 & 1 & 0 & 0  \\
0 & 0 & 1 & 0 \\
0 & -z & y & 1 \\
\end{array}\right)
\end{displaymath}

and transformations $ Ad \xi$, $ \xi \in \mathfrak{B} $ allow us to combine
only case 4 from (3)

\underline{{\bf Theorem}} Any Bol algebra of dimension 3, with the trilinear operation
of Type II, and the canonical enveloping Lie algebra of dimension 4 is isotopic
to one of the algebras 2-4.

Below we will give the description of three-Webs corresponding to the
distinguished Bol algebras.

\; \; The composition law of local Lie group $ G=< \mathbb{R}^4, \star, O> $,
tangent to the Lie algebra $ \mathfrak{G}$ is defined as follows:

$$ \begin{bmatrix} x_1\\x_2\\x_3\\x_4\end{bmatrix} \triangle
   \begin{bmatrix} y_1\\y_2\\y_3\\y_4\end{bmatrix}
   =
   \begin{bmatrix}
     x_{1}+y_{1}+ \frac{x_{4}y_{3}-x_{3}y_{4}}{2}+\frac{x^{2}_{3}y_{2}-x_{3}x_{2}y_{3}}{12}+ \frac{x_{2}y^{2}_{3}-x_{3}y_{2}y_{3}}{12}\\x_{2}+y_{2}\\x_{3}+y_{3}\\x_{4}+y_{4}+\frac{x_{2}y_{3}-y_{2}x_{3}}{2}\end{bmatrix}.
$$

Case 2. Bol algebra with bilinear and trilinear operation is defined as:

$$
e_2 \cdot e_3 =-e_3, \; (e_2, e_3, e_3)=e_1
$$
The subgroup $H$ is defined as

$$
H=\exp \mathfrak{h}=\{\exp \alpha (e_4 -e_3)\}_{ \alpha \in \mathbb{R}}.
$$
The subgroup $H$ is defined as

$$
H=\exp \mathfrak{h}=\{\exp \alpha (e_4 -e_3)\}_{ \alpha \in \mathbb{R}}=\{0,0,\alpha , \alpha \},
$$
$$
B=\{t,u,v,0\}.
$$

Any element $ (x_1,x_2, x_3, x_4) $ from $ G $ such that $ ( x_2 <2)$, can be
uniquely represented in the form:

$$ \begin{pmatrix} x_1\\x_2\\x_3\\x_4\end{pmatrix} =
   \begin{bmatrix}T\\x_2\\ \frac{2x_3+2x_4-x_2 x_3 }{2-x_2}\\0\end{bmatrix}
   \triangle
   \begin{bmatrix} 0\\0\\- \frac{2x_4}{2-x_2}\\\frac{2x_4}{2+x_2}\end{bmatrix}
$$
$$
=\prod_{B}\begin{bmatrix} x_1\\x_2\\x_3\\x_4\end{bmatrix}
  \triangle
  \begin{bmatrix}0\\0\\- \frac{2x_4}{2-x_2}\\\frac{2x_4}{2+x_2} \end{bmatrix}.
$$

Where
$$ T= \frac{12x^2_4 -4x_2 (x_4)^2 +(x_2)^2 x_4 x_3 -8x_2 x_3 x_4 +12x_4 x_3 +6x_1 (2-x_2)^2}{ (6(2-x_2)^2)}.
$$
\; \; The composition law $( \star )$, of local analytic Bol loop $B( \star )$
 is defined as:

$$ \begin{pmatrix} t\\u\\v\end{pmatrix} \star
   \begin{pmatrix}t'\\u'\\v'\end{pmatrix}
   = \prod_{B}\left(
    \begin{bmatrix}t\\u\\v\\0\end{bmatrix}\triangle \begin{bmatrix}
     t'\\u'\\v'\\0\end{bmatrix}\right)
$$
$$
= \begin{bmatrix}\frac{L}{6(2-(u+u'))^2}\\u+u'\\v+v'+ \frac{uv'-vu'}{2-u-u'}\end{bmatrix}.
$$

Where
\begin{gather*}
L=24t+24t'-24tu-24t'u-24tu'-24t'u'+12tuu'+12t'uu'\\
 +6tu^2 +6t'u^2+6t(u')^2 +6t'(u')^2 -2uvu'-2uu'v'\\
 +6uvv'-6vu'v'-6uu'v'v+2u^2 vu'+2u^2 u'v'-4u^2 vv'\\
 +2u(v')^2-4v^2 u'+2uv^2 u'+2uv(u')^2 -6uu'(v')^2\\
 +6u(v')^2+2u(u')^2 v'+4v(u')^2 v'-3u^2 (v')^2 \\
 +5v^2 (u')^2 -u^2 v(u')^2 +u^2 u'(v')^2\\
 -u^2 (u')^2 v'+u(u')^2 (v')^2 -v^2 (u')^3 +\frac{5u^2 vu'v'}{2}-uv^2 (u')^2\\
 +\frac{3uvv'({u'}^2)}{2}-\frac{u^3 vu'}{2}-\frac{u'v'u^3}{2}-\frac{uv(u')^3}{2}\\
-\frac{u(u')^3}{2}-\frac{vv'(u')^3}{2}+\frac{vv'(u)^3}{2}.
\end{gather*}

The corresponding local 3-Webs can be realized as a 6-order space.

 Bol algebra  in case 3 with bilinear and trilinear operation is defined as:

$$
e_2 \cdot e_3 =-e_2, \; (e_2, e_3, e_3)=e_1
$$
has a 4-dimensional canonical enveloping Lie algebra
$ \mathfrak{G}=<e_1,e_2,e_3,e_4> $
$$
\mathfrak{G}=\mathfrak{B} \dotplus \mathfrak{h}, \;\mathfrak{B}=<e_1,e_2,e_3>, \mathfrak{h}=<e_4+e_2>
$$
with the composition law indicated above the subgroup

$$
H=\exp \mathfrak{h}=\{\exp \alpha (e_4 +e_2)\}_{ \alpha \in \mathbb{R}}=\{0,\alpha ,0, \alpha \}
$$
$$
B=\{t,u,v,0\}.
$$

Any element $ (x_1,x_2, x_3, x_4) $ from $ G $ such that $ ( x_2 <2)$, can be
uniquely represented in the form:

$$ \begin{pmatrix} x_1\\x_2\\x_3\\x_4\end{pmatrix} =
   \begin{bmatrix} \frac{-6x_1 x_3 +6x_4 x_3-x_4 x^2_3 +12x_1}{6(2-x_3)}\\ \frac{2x_2 -2x_4 -x_2 x_3}{2-x_3}\\ x_3\\0\end{bmatrix}
   \triangle
   \begin{bmatrix} 0\\\frac{2x_4}{2-x_3}\\0\\\frac{2x_4}{2-x_3}\end{bmatrix}
$$
$$
=\prod_{B}\begin{bmatrix} x_1\\x_2\\x_3\\x_4\end{bmatrix}
  \triangle\begin{bmatrix} 0\\\frac{2x_4}{2-x_3}\\0\\\frac{2x_4}{2-x_3}\end{bmatrix}.
$$
\; \; The composition law$ ( \star )$, of local analytic Bol loop $B( \star )$ is
defined as:

$$ \begin{pmatrix} t\\u\\v\end{pmatrix} \star
   \begin{pmatrix}t'\\u'\\v'\end{pmatrix}
   = \prod_{B}\left(
    \begin{bmatrix}t\\u\\v\\0\end{bmatrix}\triangle \begin{bmatrix}
     t'\\u'\\v'\\0\end{bmatrix}\right)
$$
$$
= \begin{bmatrix}t+t'+\frac{B''}{6(2-v-v')}\\u+u'+\frac{u'v-v'u}{2-v-v'}\\v+v'\end{bmatrix}.
$$

Where
$$
B''=v^2 u'-uvv'+u(v)^2 -vu'v'-u'v^3 +uv'v^2 +3uvv'-3u'v^2 -3vu'v'+uvv'-u'v'v^2 .
$$

The corresponding local 3-Webs can be realized as a 5-order space.

 Bol algebra  in case 4 with bilinear and trilinear operation is defined as:

$$
e_2 \cdot e_3 =-e_1, \; (e_2, e_3, e_3)=e_1,
$$
has a 4-dimensional canonical enveloping Lie algebra
$ \mathfrak{G}=<e_1,e_2,e_3,e_4> $,
$$
\mathfrak{G}=\mathfrak{B} \dotplus \mathfrak{h}, \;\mathfrak{B}=<e_1,e_2,e_3>, \mathfrak{h}=<e_4-e_1>
$$
with the composition law indicated above. The subgroup

$$
H=\exp \mathfrak{h}=\{\exp \alpha (e_4 -e_1)\}_{ \alpha \in \mathbb{R}}=\{-\alpha ,0,0, \alpha \}
$$
$$
B=\{t,u,v,0\}.
$$

Any element $ (x_1,x_2, x_3, x_4) $ from $ G $,  can be uniquely represented in
the form:

$$ \begin{pmatrix} x_1\\x_2\\x_3\\x_4\end{pmatrix} =\begin{bmatrix} x_1+ x_4\\ x_2\\ x_3\\0\end{bmatrix}\triangle \begin{bmatrix} -x_4\\0\\0\\x_4\end{bmatrix}
$$
$$
=\prod _{B}\begin{bmatrix} x_1\\x_2\\x_3\\x_4\end{bmatrix}
  \triangle\begin{bmatrix} 0\\\frac{2x_4}{2-x_3}\\0\\\frac{2x_4}{2-x_3}\end{bmatrix}.
$$
\; \; The composition law $( \star )$, of local analytic Bol loop $B( \star )$ is
defined as:

$$ \begin{pmatrix} t\\u\\v\end{pmatrix} \star
   \begin{pmatrix}t'\\u'\\v'\end{pmatrix}
   = \prod_{B}\left(
    \begin{bmatrix}t\\u\\v\\0\end{bmatrix}\triangle \begin{bmatrix}
     t'\\u'\\v'\\0\end{bmatrix}\right)
$$
$$
= \begin{bmatrix}t+t'+\frac{uv'-vu'}{2}+\frac{v^2 u'-uvv'}{12}+\frac{u(v')^2 -vu'v'}{12}\\u+u'\\v+v'\end{bmatrix}.
$$

The corresponding 3-web is global and can be realized as a 3-order space.

{\bf Remark:} case 1 has not been investigated, because it is combined with case 4 see proposition 2.

We pass to the classification of Bol algebras of dimension 3, with trilinear
operation of type II with 5-dimensional canonical enveloping Lie algebra.

The structural constants of Lie algebra
$ \mathfrak{G}=<e_1,e_2,e_3,e_4,e_5>$
$ \mathfrak{G}=\mathfrak{B} \dotplus V $, must satisfy the relations:
$\mathfrak{B}=<e_1,e_2,e_3>$, $[\mathfrak{B},\mathfrak{B}] \subset V=<e_4,e_5>$,
$[V,\mathfrak{B}] \subset \mathfrak{B} $, $[e_5,\mathfrak{B}]=0$
$$
[e_2,e_3]=e_4 +re_5,
$$
$$
[e_3,e_1]=se_5,
$$
$$
[e_4,e_3]=e_1,
$$
(trivial relations have been omitted) in case where $s=0$, we have
$ <\mathfrak{B} \neq \mathfrak{G}$, which means Lie algebras $ \mathfrak{G}$
is not a canonical enveloping. That is why by adjusting, if necessary vector
$e_5$ we have $s=1$ (for any $r$).

Here we will have two families of subalgebras $ \mathfrak{h} $ in
$ \mathfrak{G}$, such that $ \mathfrak{G}=\mathfrak{B} \dotplus \mathfrak{h}$.
$$
\mathfrak{h}=<e_1+ \alpha e_2 + \beta e_3 + \gamma e_4 ,e_5>, \; \forall  \alpha , \beta , \gamma \in \mathbb{R}
$$
$$
\mathfrak{h}_{ \alpha, \beta,\bar{\alpha} ,\bar{\beta} } =<e_4+ \alpha e_1 + \beta e_2  ,e_5 +\bar{\alpha}e_1  +\bar{\beta}e_2>, \; \forall  \alpha , \beta , ,\bar{\alpha} ,\bar{\beta} \in \mathbb{R}.
$$

The first subalgebras contain ideal $<e_5>$, by factorization with the ideal we
get the case considered above (Theorem III.2.1, and 2). That is why below we
will examine only the subalgebras of the second kind which do not contain any
ideals of the Lie algebra $ \mathfrak{G}$, and are not ideals themselves
(That is why the description of correspondent algebras is not reduced to the
 case $ dim \mathfrak{G}=4)$.

\; \; Let us note that the automorphism $A$ of Lie algebra $ \mathfrak{G} $,
which is the extension of an automorphism of the Lie triple system
$\mathfrak{G} $ (with the trilinear operation Type II [11]), is defined as
follows:

\begin{displaymath}
ad ( \xi )= \left(\begin{array}{ccccc}
bf^2 & a & d & 0 & 0 \\
0 & b & l & 0 & 0 \\
0 & 0 & f & 0 & 0 \\
0 & 0 & 0 & bf & 0 \\
0 & 0 & 0 & -af & bf^3 \\
\end{array}\right),
\end{displaymath}

$ b,f \neq 0 $, \; $a,d,l \in \mathbb{R}$

so that

$$
 A(e_4 +\alpha e_1 +\beta e_2)=e_4 +\alpha_n e_1 +\beta_n e_2,
$$
$$
A(e_5 +\bar{\alpha} e_1 +\bar{\beta} e_2)=e_4 +\bar{\alpha_n} e_1 +\bar{\beta_n} e_2,
$$
where,
$$
\alpha_n=\alpha f+\frac{\beta}{f}\frac{a}{b}+\frac{\bar{\beta}}{f^2}\frac{a^2}{b^2},
$$
$$
\beta_n=\frac{\beta}{f}+\frac{\bar{\beta}}{f^3}\frac{a}{b},
$$
$$
\bar{\alpha_n}=\frac{\bar{\alpha}}{f}+\frac{\bar{\beta}}{f^3}\frac{a}{b},
$$
$$
\bar{\beta_n}=\frac{\bar{\beta}}{f^3}.
$$

\begin{itemize}
\item I. if $ \bar{\beta}=0$ then $ \bar{\beta_n}=0$ since $ \bar{\alpha} \neq 0$
otherwise we will get an ideal, by factorization through it, we will get the
above case. Therefore we can make $ \bar{\alpha_n}=1 $ then we will have two
possibilities:
\begin{itemize}
\item $ \beta_n =-1$ and $\alpha_n=\alpha f+\frac{\beta}{\bar{\alpha}}\frac{a}{b}$ is any number
\item $\beta_n=\frac{\beta}{\bar{\alpha}}$ is any number not equal to $-1$ and we make $ \alpha_n =0$.
\end{itemize}
\item II. If $ \bar{\beta} \neq 0 $, then $ \bar{\beta_n}=1$
\begin{itemize}
\item If $ \frac{a}{b}=-\frac{\beta}{f}$, then $ \beta_n =0$,
$ \bar{\alpha_n}=\frac{\bar{\alpha}}{f}+\frac{a}{b}$- is any number,
 and $ \alpha_n=\alpha f+\frac{a}{b}\frac{\beta}{f}+\frac{a}{b}\frac{\bar{\alpha}}{f}+(\frac{a}{b})^2 f $ is any number;
\item if $ \beta =0$, then $ \beta_n =\frac{a}{b}$ we can take $ a=0$, $ \beta_n =0$, then $ \bar{\alpha_n}=0$, $ \beta_n = \alpha f$- any number.
\end{itemize}
\end{itemize}

\; \; \; Here we note that the situation $b$ above is the particular case of
the situation $a$ hence within the isomorphism there are
$$
(\nu, -1, 1, 0), (0,\beta_n=\frac{\beta}{\bar{\alpha}} \neq -1,0)
$$
$$
(\mu, 0, \theta, 1) \; \; \; (4)
$$
where
$$
\nu=\alpha f+\frac{\beta}{\alpha}\frac{a}{b},\mu=\alpha f+\frac{\beta}{\alpha}\frac{a}{b}+\frac{\bar{\alpha}}{f}\frac{a}{b}+({a}{b})^2 f , \theta=\frac{\bar{\alpha}}{f}+\frac{a}{b}.
$$

\; \; The elements (4) form a section of orbits which is obtained under the
action of the automorphism $A$ in the 4-dimensional space $( \alpha, \beta, \bar{\alpha}, \bar{\beta})$ respectively.

The following proposition holds

\underline{{\bf Theorem III.2.3}} Any Bol algebra of dimension 3, with the
trilinear
operation of Type II, and the canonical enveloping Lie algebra of dimension 5, is isomorphic to one of the Bol algebras of the form (4)

1.$$
e_2 \cdot e_3 =- \nu e_1 +e_2, \; - \nu \geq 0,
$$
$$
e_1 \cdot e_3 =e_1 \; \; \; \; \; \; \; (e_2, e_3, e_3)=e_1,
$$
2.$$
e_2 \cdot e_3 =- \frac{\beta}{\bar{\alpha}}e_2, \; - \frac{\beta}{\bar{\alpha}} \geq 0, \beta \neq - \bar{\alpha},
$$
$$
e_1 \cdot e_3 =e_1 \; \; \; \; \; \; \; (e_2, e_3, e_3)=e_1,
$$
3.$$
e_2 \cdot e_3 =- \mu e_1 ,\; \; \; - \mu \geq 0,
$$
$$
e_1 \cdot e_3 =- \theta e_1 -e_2, \; \; \theta \geq 0 \; \;  \; (e_2, e_3, e_3)=e_1.
$$

{\bf Remark.} The classification of Bol algebra up to the isotopy is not
examined here.

\; \; \; Below, we give an example of 3-Web corresponding to the case
$(0, \beta , 1, 0)$ where $\beta \neq -1$.

Bol algebra $ \mathfrak{B}$ with the trilinear and bilinear operations
defined as follows:

$$
e_2 \cdot e_3 =\beta e_2,
$$
$$
e_1 \cdot e_3 =e_1 \; \; \; \; \; \; \; (e_2, e_3, e_3)=e_1
$$

has a 4-dimensional canonical enveloping Lie algebra of dimension 5,
$\mathfrak{G}=<e_1,e_2,e_3,e_4,e_5>$, $\mathfrak{G}=\mathfrak{B} \dotplus \mathfrak{h}$, $\mathfrak{h}=<e_4 -\beta e_2,e_5 -e_1>$.

The composition law $(\triangle)$, corresponding to the Lie group $G$, is defined as follows:

$$ \begin{bmatrix} x_1\\x_2\\x_3\\x_4\\x_5\end{bmatrix} \triangle
   \begin{bmatrix} y_1\\y_2\\y_3\\y_4\\y_5\end{bmatrix}
   =
   \begin{bmatrix}
     x_{1}+y_{1}+ \frac{x_{4}y_{3}-x_{3}y_{4}}{2}-\frac{x_{2}x_{3}y_{3}-x^{2}_{3}y_{2}}{12}+ \frac{x_{3}y_{2}y_{3}-x_{2}y^{2}_{3}}{12}\\x_{2}+y_{2}\\x_{3}+y_{3}\\x_{4}+y_{4}+\frac{x_{2}y_{3}-y_{2}x_{3}}{2}\\ x_{5}+y_{5}+ \frac{x_{1}y_{3}-x_{3}y_{1}}{2}-\frac{x_{3}x_{4}y_{3}-x^{2}_{3}y_{4}}{12}+ \frac{x_{3}y_{3}y_{4}-x_{4}y^{2}_{3}}{12}\end{bmatrix}.
$$

The subgroup H is defined as the collection of elements

$$
H=\exp \mathfrak{h}=\{\exp t(e_4 - \beta e_2), \exp q(e_5 -e_1)\}_{ t,q \in \mathbb{R}}=\{-q, -t \beta, 0 , t, q \}_{ t,q \in \mathbb{R}}.
$$

 The collection of elements
$$
 B=\{t,u,v,0,0\}_{ t,u,v  \in \mathbb{R}}
$$
form a local section of left coset space $G \bmod H$.

Any element $ (x_1,x_2, x_3, x_4,x_5) $ from $ G $,  can be
uniquely represented as follows:

$$ \begin{pmatrix} x_1\\x_2\\x_3\\x_4\\x_5\end{pmatrix} =
   \begin{bmatrix}\frac{T}{6(1+x_3)(2+x_3)}\\\frac{2 \beta x_4}{2+ \beta x_3}+x_2\\  x_3 \\0\\0\end{bmatrix}
   \triangle
   \begin{bmatrix} \frac{x_4 \cdot x^2_3}{6(2+ \beta x_3)}-x_5\\-\frac{2 \beta x_4}{2+ \beta x_3}\\0\\\frac{2x_4}{2+ \beta x_3}\\-\frac{x_4 \cdot x^2_3}{6(2+ \beta x_3)}+x_5\end{bmatrix}
$$
$$
=\prod _{B}\begin{bmatrix} x_1\\x_2\\x_3\\x_4\\x_5\end{bmatrix}
  \triangle    \begin{bmatrix} \frac{x_4 \cdot x^2_3}{6(2+ \beta x_3)}-x_5\\-\frac{-2 \beta x_4}{2+ \beta x_3}\\0\\\frac{2x_4}{2+ \beta x_3}\\-\frac{x_4 \cdot x^2_3}{6(2+ \beta x_3)}+x_5\end{bmatrix}
$$

where

\begin{multline}
 T= 12x_1 +6 \beta x_1 x_3 +12x_1 x_3+ 6 \beta x_1 (x_3)^2 +12x_5+ 6 \beta x_3 x_5 -6(x_3)^2 x_4+\\ +6x_3 x_4
+ \beta x_4 (x_3)^2 + \beta (x_3)^2 x_4.
\end{multline}

\; \; The composition law $( \star )$ of local analytic Bol loop $B( \star )$ is
defined as:

$$ \begin{pmatrix} t\\u\\v\end{pmatrix} \star
   \begin{pmatrix}t'\\u'\\v'\end{pmatrix}
   = \prod_{B}\left(
    \begin{bmatrix}t+t'-\frac{uvv'-u'v^2}{12}-\frac{vu'v'-u(v')^2}{12}\\u+u'\\v+v'\\\frac{uv'-u'v}{2}\\\frac{tv'-vt'}{2}\end{bmatrix}\right)
$$
$$
= \begin{bmatrix}X\\u+u'+\frac{\beta (uv'-u'v)}{2+ \beta (v+v')}\\v+v'\end{bmatrix}.
$$
Where

\begin{multline}
 X=t+t'-\frac{uvv'-vu'}{12}-\frac{vu'v'-u(v')^2}{12}+\frac{tv'-vt'}{2}-
\frac{(uv'-vu')(v+v')}{2}+\\
+(\frac{v+v}{2}+\frac{\beta(v+v')^2}{12})\frac{(uv'-vu')}{2+ \beta (v+v')}.
\end{multline}

The corresponding local 3-Webs can be realized as a 6-order space.
\newpage
\subsection{BOL ALGEBRAS WITH TRILINEAR OPERATION OF TYPE III}

\; \; As in the previous chapter we will base our investigation of
3-dimensional Bol algebras, on the examination of their canonical enveloping
Lie algebras. In what follows, we consider Bol algebras of dimension 3, from
their construction see [34,36]; it follows that the dimension of their
canonical enveloping Lie algebras can not be more than 6. Below we limit
ourselves to the classification of Bol algebras (and their corresponding
3-webs) with canonical enveloping Lie algebras of dimension$ \leq 4$.

Let $ \mathfrak{B} $ be a 3-dimensional Bol algebra with a trilinear
operation Type III see chapter II \S 3, and
$\mathfrak{G}=\mathfrak{B}+\mathfrak{h}$-its canonical enveloping Lie algebra according to the table given in chapter II (case 4). We note that the situation
$dim \mathfrak{G}=3$ is possible that means we obtain a total grouped 3-Web the
corresponding Lie group $G$, is isomorphic to the matrix of the view:

\begin{displaymath}
 \left(\begin{array}{ccc}
1 & 0 & 0 \\
-(y-z)\frac{e^2x -1}{2x} & \frac{e^2x +1}{2} & -\frac{e^2x -1}{2} \\
(y-z)\frac{e^2x -1}{2x} & -\frac{e^2x +1}{2} & \frac{e^2x -1}{2} \\
\end{array}\right),
\end{displaymath}
with $ x,y,z \in \mathbb{R}$.

\; \; Let us examine the case $dim \mathfrak{G}=4$, the structural constants of Lie algebra $ \mathfrak{G}=<e_1,e_2,e_3,e_4>$, $\mathfrak{B}=<e_1,e_2,e_3>$
are defined as follows:
$$
[e_1,e_2]=e_4, \; \; [e_2,e_4]= \mp e_1
$$
in addition $ \mathfrak{G}=\mathfrak{B} \dotplus [\mathfrak{B},\mathfrak{B}]$,
$[\mathfrak{B}, \mathfrak{B}]=<e_4>$.

By introducing in consideration the 3-dimensional subspaces of subalgebras
$$
\mathfrak{h}_{x,y,z}=<e_4 +xe_1 +ye_2 +ze_3>,\; \; x,y,z, \in \mathbb{R}
$$
we obtain a collection of Bol algebras of view:

*****

if $ y \neq 0$ then choosing $  \alpha , a $ and $c$ one can make
$y'=1$, $x'=0$, $z'=0$.

If $y=0$, then $y'=0$,  choose $x'= \pm x \geq 0$, and $z'=0,1$.

Thus within the isomorphism are obtained two families of Bol algebras and one
exceptional Bol algebra.

\underline{{\bf Theorem$ III$ .3.1}} Any Bol algebra of dimension 3, with the
trilinear
operation of Type $III^-$, and the canonical enveloping Lie algebra of dimension
 4, is isomorphic to one of Bol algebras below:
\begin{itemize}
\item $ III^- .1\;\; e_1 \cdot e_2 =-e_2, \; \; (e_1,e_2,e_2)=e_1$,
\item $ III^- .1\;\; e_1 \cdot e_2 =-xe_1, \; \; (e_1,e_2,e_2)=e_1$ $x \geq 0$,
\item $ III^- .1\;\; e_1 \cdot e_2 =-xe_1-e_3, \; \;(e_1,e_2,e_2)=e_1$ $x \geq 0 $.
\end{itemize}

In this the distinguished Bol algebras are not isomorphic among themselves.

Similarly one can establish the correctness of the following Theorem.

\underline{{\bf Theorem III .3.2}} Any Bol algebra of dimension 3, with the
trilinear
operation of Type $III^+$, and the canonical enveloping Lie algebra of dimension
 4, is isomorphic to one of Bol algebras below:
\begin{itemize}
\item $ III^+ .1\;\; e_1 \cdot e_2 =-e_2, \; \; (e_1,e_2,e_2)-=e_1$,
\item $ III^+ .1\;\; e_1 \cdot e_2 =-xe_1, \; \; (e_1,e_2,e_2)-=e_1$ $x \geq 0$,
\item $ III^+ .1\;\; e_1 \cdot e_2 =-xe_1-e_3, \; \;(e_1,e_2,e_-2)=e_1$ $x \geq 0$.
\end{itemize}

In this the distinguished Bol algebras are not isomorphic among themselves.

We pass to the classification within isotopic of Bol algebras, enumerated in
Theorem III.31.

We will note that for every $ \xi=te_1 +ue_2 +ve_3 $ from $ \mathfrak{B}$
$ t,u,v, \in \mathbb{R} $

\begin{displaymath}
ad ( \xi )= \left(\begin{array}{cccc}
0 & 0 & 0 & -u \\
0 & 0 & 0 & 0  \\
0 & 0 & 0 & 0 \\
-u & t & 0 & 0 \\
\end{array}\right),
\end{displaymath}

\begin{displaymath}
Ad ( \xi )= \left(\begin{array}{cccc}
\cosh u & \frac{t}{u}(1-\cosh u) & 0 & -\sinh u \\
0 & 1 & 0 & 0  \\
0 & 0 & 1 & 0 \\
-\sinh u & \frac{t}{u}\sinh u & 0 & \cosh u \\
\end{array}\right).
\end{displaymath}

Let us calculate the image of $\Phi(\mathfrak{h})$ under the action of
$ \Phi=Ad \xi$ on the one-dimensional subalgebra $\mathfrak{h}$ with a
directing vector
$$
e_4 +xe_1 +ze_3, x \geq 0;1
$$
$$
\Phi (e_4 +xe_1 +ze_3)=(\cosh u -x\sinh u)e_4 +ze_3 + (x \cosh u -sinh u)e_1,
$$
in addition one can define the mapping
$$
x \longrightarrow \frac{x\cosh u-\sinh u}{-x\sinh u+\cosh}=x',
$$
$$
z \longrightarrow \frac{z}{-x\sin hu+\cosh u}=z'.
$$

By choosing $u$ such that $ \tanh u=x$, we obtain $x'=0$. We notice that
in addition $ \coth u \neq 0$, that is, the mapping is defined correctly.

By applying the automorphism (9) to the so obtained Bol algebras, one can consider $z'=0;1$. Moreover, the two isolated case are not isotopic (in the sense
of the definition see chapter I). By virtue of the formulas (5) and (9).

We note that the application of the isotopic transformation $ \Phi $ to the exceptional Bol algebra of Theorem III.3.1 is not changed.

Summarizing the conducted examination one can formulate the Theorem

\underline{{\bf Theorem III .3.3}} Any Bol algebra of dimension 3, with the
trilinear
operation of Type $III^-$, and the canonical enveloping Lie algebra$ \mathfrak{G}$ of dimension 4, is isotopic to one of the following Bol algebras:
\begin{itemize}
\item $  e_1 \cdot e_2 =-e_2, \; \; (e_1,e_2,e_2)=e_1$,
\item $  e_1 \cdot e_2 =-ze_1, \; \; (e_1,e_2,e_2)=e_1$ $z= 0;1$.
\end{itemize}

Analogically one can state the correctness of the Theorem.

\underline{{\bf Theorem III .3.4}} Any Bol algebra of dimension 3, with the
trilinear
operation of Type $III^+$, and the canonical enveloping Lie algebra
$ \mathfrak{G}$ of dimension 4, is isotopic to one of Bol algebras below:
\begin{itemize}
\item $  e_1 \cdot e_2 =-e_2, \; \; (e_1,e_2,e_2)-=e_1$,
\item $  e_1 \cdot e_2 =-ze_1, \; \; (e_1,e_2,e_2)-=e_1$ $x=0, 1$.
\end{itemize}

 Below we reduce to description of 3-Webs corresponding to the isolated
Bol algebras of Type$ III^-$ and Type$ III^+$.

 The composition law $(\triangle)$,
corresponding to the Lie group $G$ of Lie algebra enveloping Bol algebra is
defined as follows:

$$ \begin{bmatrix} x_1\\x_2\\x_3\\x_4\end{bmatrix} \triangle
   \begin{bmatrix} y_1\\y_2\\y_3\\y_4\end{bmatrix}
   =
   \begin{bmatrix}
     x_{1}+y_{1}+ \cosh x_{2} -y_{4}sinh x_{2}\\x_{2}+y_{2}\\x_{3}+y_{3}\\x_{4}-y_{1}+\sinh x_{2} +y_{4}\cosh x_{2}\end{bmatrix}.
$$

 In case $ III^-.1$ the subgroup $ H=\exp \mathfrak{h} $ can be realized as
the
collection of elements

$$
H=\exp \mathfrak{h}=\{\exp \alpha (e_4 +e_2)\}_{ \alpha \in \mathbb{R}}=\{0,\alpha ,0, \alpha \}_{ \alpha \in \mathbb{R}}.
$$
The collection of elements
$$
B=\exp \mathfrak{B}=\{\frac{t}{u}\sinh u,u,v,\frac{t}{u}(1-\cosh u)\}_{t,u.v \in \mathbb{R}}
$$
form a local section of left space coset $G \bmod H $.
$$
\exp^{-1}  \begin{bmatrix} x_1\\x_2\\x_3\\x_4\end{bmatrix}= \begin{bmatrix} \frac{x_{1}x_{2}}{\sinh x_{2}}\\x_2\\x_3\end{bmatrix},
$$
$  x_1,x_2,x_3,x_4 \in \mathbb{R}$.

Any element $ (x_1,x_2, x_3, x_4) $ from $ G $, in the neighborhood  $e$, can be
uniquely represented as follows:

$$ \begin{pmatrix} x_1\\x_2\\x_3\\x_4\end{pmatrix} =
   \begin{bmatrix}(x_1 \cosh p +x_4 \sinh p)\frac{p}{\sinh p}\\p\\ x_3\\(x_1 \cosh p +x_4 \sinh p)\frac{1-\cosh p}{\sinh p}\end{bmatrix}
   \triangle
   \begin{bmatrix} 0\\(x_4 -x_1 (\frac{1}{\sinh u}-\coth u))\\(x_4 -x_1 (\frac{1}{\sinh u}-\coth u))\end{bmatrix},
$$

where $p$ is defined from the relation
$$
p+x_4 -x_1 \frac{1- \cosh p}{\sinh p}=x_2.
$$

\; \; The composition law $( \star )$, corresponding to the local analytical Bol
loop $B( \star )$, is defined as follows:

$$ \begin{pmatrix} t\\u\\v\end{pmatrix} \star
   \begin{pmatrix}t'\\u'\\v'\end{pmatrix}
   = \exp^{-1} \left(\prod_{B}\left(
    \begin{bmatrix}t\\u\\v\\0\end{bmatrix}\triangle \begin{bmatrix}
     t'\\u'\\v'\\0\end{bmatrix}\right) \right)
$$
$$
= \exp^{-1}
  \left( \begin{bmatrix}
           \left[\cosh p +t'\cosh (u-p)\right]\frac{p}{\sinh p}\\
           p\\
           v+v'\\
           \left[t\cosh p +t'\cosh (u-p)\right] \frac{(1-\cosh p)}{\sinh p}
  \end{bmatrix} \right)
$$.

$$
= \begin{bmatrix}\left[t\cosh p +t'\cosh (u-p)\right]
     (\frac{p}{\sinh p})^2\\p\\v+v'\end{bmatrix},
$$

where $p$ is defined from the relation:
$$
p+\frac{-t'\sinh u \sinh p -(t+t'\cosh u)(1-\cosh p)}{\sinh p}=u+u'.
$$

 The corresponding local analytic 3-Web, can be realized as a hyperbolic space.

 In case $ III^-.2.$  the subgroup $ H=\exp \mathfrak{h} $ can be realized
as the
collection of elements

$$
H=\exp \mathfrak{h}=\{\exp \alpha (e_4 +xe_1)\}_{\alpha \in \mathbb{R}}=\{x\alpha ,0,0, \alpha \}_{ \alpha \in \mathbb{R}}.
$$
Bol loop can be realized from the elements
$$
B= \exp \mathfrak{B}=\{\frac{t}{u}\sinh u,u,v,\frac{t}{u}(1-\cosh u)\}_{t,u.v \in \mathbb{R}}
$$
which form a local section of left space coset $G \bmod H $.
$$
\exp^-1  \begin{bmatrix} x_1\\x_2\\x_3\\x_4\end{bmatrix}= \begin{bmatrix} \frac{x_{1}x_{2}}{\sinh x_{2}}\\x_2\\x_3\end{bmatrix}
  \quad x_1,x_2,x_3,x_4 \in \mathbb{R}$$.

Any element $ (x_1,x_2, x_3, x_4) $ from $ G $, in the neighborhood  $e$, can be
uniquely represented as follows:

\begin{multline}
 \begin{pmatrix} x_1\\x_2\\x_3\\x_4\end{pmatrix} =
   \begin{bmatrix}\frac{\left[x_1 (x\sinh x_2 -\cosh x_2)+x_4 (x\cosh x_2)- \sinh x_2)\right]\sinh x_2 }{x(-1+\cosh x_2)-\sinh x_2}\\x_2\\ x_3\\\frac{\left[x_1 (x\sinh x_2 -\cosh x_2)+x_4 (x\cosh x_2)- \sinh x_2)\right](1-\cosh x_2) }{x(-1+\cosh x_2)-\sinh x_2} \end{bmatrix}
   \triangle \\
   \begin{bmatrix}
 \frac{x(x_1 (1-\cosh x_2)-x_4 \sinh x_2)}{x(-1+\cosh x_2)-\sinh x_2}\\0\\0\\ \frac{x(x_1 (1-\cosh x_2)+x_4 \sinh x_2)}{x(-1+\cosh x_2)-\sinh x_2}\end{bmatrix}.
\end{multline}

\; \; The composition $law( \star )$ corresponding to the local analytical Bol
loop $B( \star )$ is defined as follows:

$$ \begin{pmatrix} t\\u\\v\end{pmatrix} \star
   \begin{pmatrix}t'\\u'\\v'\end{pmatrix}
   = \exp^{-1} \left(\prod_{B}
    \begin{bmatrix}t+t'\\u+u'\\v+v'\\-t'\sinh u\end{bmatrix}\right)
$$
$$
= \exp^{-1} \begin{bmatrix}\frac{\left[(t+t'\cosh u)(x\sinh (u+u')-\cosh (u+u'))-t'\sinh u (x\cosh (u+u') -\sinh (u+u'))\right]\sinh (u+u')}{x(\cosh (u+u') -1)-\sinh (u+u')}\\u+u'\\v+v'\\ \frac{\left[(t+t'\cosh u)(x\sinh (u+u')-\cosh (u+u'))-t'\sinh u (x\cosh (u+u') -\sinh (u+u'))\right](1-\cosh (u+u'))}{x(\cosh (u+u') -1)-\sinh (u+u')}\end{bmatrix}
$$

$$
= \begin{bmatrix}\frac{L_1 (u+u')}{L_2}\\u+u'\\v+v'\end{bmatrix},
$$

where

\begin{multline}
L_1=[(t+t'\cosh u)(x\sinh (u+u')-\cosh (u+u'))-\\
t'\sinh u (x\cosh (u+u') -\sinh (u+u'))],
\end{multline}

$$
L_2= x(\cosh (u+u') -1)-\sinh (u+u').
$$
The corresponding local analytic 3-Web, can be realized as a hyperbolic space
type.

 In case $ III^-.3.$  the subgroup $ H=\exp \mathfrak{h} $ can be realized as
the
collection of elements

$$
H=\exp \mathfrak{h}=\{exp \alpha (e_4 +xe_1)\}_{ \alpha \in \mathbb{R}}=\{x\alpha ,0,0, \alpha \}_{ \alpha \in \mathbb{R}}.
$$
Bol loop can be realized from the elements
$$
B=\exp \mathfrak{B}=\{\frac{t}{u}\sinh u,u,v,\frac{t}{u}(1-\cosh u)\}_{t,u.v \in \mathbb{R}}
$$
which form a local section of left space coset $G \bmod H $.

Any element $ (x_1,x_2, x_3, x_4) $ from $ G $, in the neighborhood  $e$, can be
uniquely represented as follows:

*****

 In case $ III^+.1$ the subgroup $ H=\exp \mathfrak{h} $, can be realized as the
collection of elements

$$
H=\exp \mathfrak{h}=\{\exp \alpha (e_4 +e_2)\}_{ \alpha \in \mathbb{R}}=\{0,\alpha ,0, \alpha \}_{ \alpha \in \mathbb{R}}.
$$
The collection of elements
$$
B=\exp \mathfrak{B}\{\frac{t}{u}\sinh u,u,v,\frac{t}{u}(1-\cosh u)\}_{t,u.v \in \mathbb{R}}
$$
form a local section of left space coset $G \bmod H $.
$$
\exp^{-1}  \begin{bmatrix} x_1\\x_2\\x_3\\x_4\end{bmatrix}= \begin{bmatrix} \frac{x_{1}x_{2}}{\sinh x_{2}}\\x_2\\x_3\end{bmatrix}
\quad  x_1,x_2,x_3,x_4 \in \mathbb{R}.$$

Any element $ (x_1,x_2, x_3, x_4) $ from $ G $, in the neighborhood  $e$, can be
uniquely represented as follows:

$$ \begin{pmatrix} x_1\\x_2\\x_3\\x_4\end{pmatrix} =
   \begin{bmatrix}(x_1 \cos p +x_4 \sin p)\frac{p}{\sin p}\\p\\ x_3\\(x_1 \cos p +x_4 \sin p)\frac{1-\cos p}{\sin p}\end{bmatrix}
   \triangle
   \begin{bmatrix} 0\\(x_4 -x_1 (\frac{1}{\sin u}-\cot u))\\(x_4 -x_1 (\frac{1}{\sin u}-\cot u))\end{bmatrix},
$$

where $u$ is defined from the relation
$$
u+x_4 -x_1( \cot u- \csc u)=x_2.
$$

\; \; The composition law $( \star )$, corresponding to the local analytical Bol
loop $B( \star )$, is defined as follows:

$$ \begin{pmatrix} t\\u\\v\end{pmatrix} \star
   \begin{pmatrix}t'\\u'\\v'\end{pmatrix}
   = \exp^{-1} \left(\prod_{B}
    \begin{bmatrix}t+t'\\u+u'\\v+v'\\-t'\sin u\end{bmatrix}\right)
$$

$$
= \exp^{-1} \left( \begin{bmatrix}\left[t\cos p +t'\cos (u+p)\right]\frac{p}{\sin p}\\p\\v+v'\\ \left[t\cos p +t'\cos (u+p)\right]\frac{(1-\cos p)}{\sin p}\end{bmatrix}\right)
$$.
$$
= \begin{bmatrix}\left[t\cos p +t'\cos (u+p)\right](\frac{p}{\sin p})^2\\p\\v+v'\end{bmatrix},
$$

where $p$ is defined from the relation:
$$
p+\frac{-t'\sin u \sin p -(t+t'\cos u)(1-\cos p)}{\sin p}=u+u'.
$$

 In case $ III^+.2.$  The subgroup $ H=\exp \mathfrak{h} $, can be realized as
the
collection of elements

$$
H=\exp \mathfrak{h}=\{\exp \alpha (e_4 +xe_1)\}_{ \alpha \in \mathbb{R}}=\{x\alpha ,0,0, \alpha \}_{ \alpha \in \mathbb{R}}.
$$
Bol loop can be realized from the elements
$$
B= \exp \mathfrak{B}=\{\frac{t}{u}\sin u,u,v,\frac{t}{u}(1-\cos u)\}_{t,u.v \in \mathbb{R}}
$$
which form a local section of left space coset $G \bmod H $.
$$
\exp^{-1}  \begin{bmatrix} x_1\\x_2\\x_3
\\x_4\end{bmatrix}= \begin{bmatrix} \frac{x_{1}x_{2}}{\sin x_{2}}\\x_2\\x_3\end{bmatrix}
$$
$  x_1,x_2,x_3,x_4 \in \mathbb{R}$.

Any element $ (x_1,x_2, x_3, x_4) $ from $ G $ in the neighborhood  $e$ can be
uniquely represented as follows:

$$ \begin{pmatrix} x_1\\x_2\\x_3\\x_4\end{pmatrix} =
   \begin{bmatrix}x_1 -\frac{(x\cos  x_2 -\sin  x_2)(x_1 (1-\cos x_2)+x_4 \sin  x_2) }{x(1-\cos x_2)-\sin x_2}\\x_2\\ x_3\\x_4 -\frac{(\cos x_2 -x\sin x_2)(x_1 (1-\cos x_2)+x_4 \sin x_2)}{x(1-\cos x_2)-\sin x_2} \end{bmatrix}
   \triangle
   \begin{bmatrix}
 \frac{x(x_1 (1-\cos x_2)+x_4 \sin x_2)}{x(1-\cos x_2)-\sin x_2}\\0\\0\\ \frac{x(x_1 (1-\cos x_2)+x_4 \sin x_2)}{x(1-\cos x_2)-\sin x_2}\end{bmatrix}.
$$

\; \; The composition law $( \star )$, corresponding to the local analytical Bol
loop $B( \star )$, is defined as follows:

$$ \begin{pmatrix} t\\u\\v\end{pmatrix} \star
   \begin{pmatrix}t'\\u'\\v'\end{pmatrix}
   = \exp^{-1} \left(\prod_{B}
    \begin{bmatrix}t+t'\cos u\\u+u'\\v+v'\\-t'\sin u\end{bmatrix}\right)
$$
$$
= \exp^{-1}  \begin{bmatrix}t+t' \cos u -\frac{N_1}{x(1-\cos(u+u'))-\sin (u+u')}\\u+u'\\v+v'\\ -t'\sin u  -\frac{N_2}{x(1-\cos (u+u'))-\sin (u+u')}\end{bmatrix}
$$

$$
= \begin{bmatrix}\left[t+t' \cos u -\frac{N_1}{x(1-\cos(u+u'))-\sin (u+u')}\right]\frac{u+u'}{\sin (u+u')}\\u+u'\\v+v'\end{bmatrix},
$$

where
$$
N_1=(x\cos (u+u')+\sin (u+u'))((t+t'\cos u)(\cos (u+u')-1)+t'\sin u \sin (u+u')),
$$
$$
N_2=(\cos (u+u')+x\sin (u+u'))((t+t'\cos u)(\cos (u+u') -1)+t'\sin u \sinh (u+u')).
$$

 In case $ III^+.3.$  The subgroup $ H=\exp \mathfrak{h} $, can be realized as the
collection of elements

$$
H=\exp \mathfrak{h}=\left\{\exp \alpha (e_4 +xe_1 +e_3) \right\}_{ \alpha \in \mathfrak{R}}=\left\{x\alpha ,0,\alpha, \alpha \right\}_{ \alpha \in \mathbb{R}}
$$
Bol loop can be realized from the elements
$$
B= \exp \mathfrak{B}=\left\{\frac{t}{u}\sin u,u,v,\frac{t}{u}(1-cos u)\right\}_{t,u.v \in \mathbb{R}}
$$
which form a local section of left space coset $G \bmod H $.
$$
\exp^{-1}  \begin{bmatrix} x_1\\x_2\\x_3\\x_4\end{bmatrix}= \begin{bmatrix} \frac{x_{1}x_{2}}{\sin x_{2}}\\x_2\\x_3\end{bmatrix}
$$
$  x_1,x_2,x_3,x_4 \in \mathbb{R}$.

Any element $ (x_1,x_2, x_3, x_4) $ from $ G $, in the neighborhood  $e$, can be
uniquely represented as follows:

$$ \begin{pmatrix} x_1\\x_2\\x_3\\x_4\end{pmatrix} =
   \begin{bmatrix}x_1 -\frac{(x\cos  x_2 +\sin  x_2)(x_1 (-1+\cos x_2)+x_4 \sin  x_2) }{x(1-\cos x_2)+\sin x_2}\\x_2\\ x_3-\frac{x(x_4 \sin x_2 -x_1 (\cos x_2-1))}{x(1-\cos x_2)+\sin x_2}\\x_4 -\frac{(\cos x_2 -x\sin x_2)(x_1 \cos x_2)+x_4 \sin x_2)}{x(1-\cos x_2)+\sin x_2} \end{bmatrix}
   \triangle
   \begin{bmatrix}
 \frac{x(-x_1 (-1+\cos x_2)+x_4 \sin x_2)}{x(1-\cos x_2)-\sin x_2}\\0\\ \frac{(-x_1 (-1+\cos x_2)+x_4 \sin x_2)}{x(1-\cos x_2)-\sin x_2}\\ \frac{(x_1 (-1+\cos x_2)+x_4 \sin x_2)}{x(1-\cos x_2)+\sin x_2}\end{bmatrix}.
$$

\; \; The composition law $( \star )$, corresponding to the local analytical Bol
loop $B( \star )$, is defined as follows:

$$ \begin{pmatrix} t\\u\\v\end{pmatrix} \star
   \begin{pmatrix}t'\\u'\\v'\end{pmatrix}
   = exp^-1 \left(\prod_{B}
    \begin{bmatrix}t+t'\cos u\\u+u'\\v+v'\\-t'\sin u\end{bmatrix}\right)
$$

$$
= \begin{bmatrix}\left[t+t' \cos u -\frac{P_1}{x(1-\cos(u+u'))+\sin (u+u')}\right]\frac{u+u'}{\sin (u+u')}\\u+u'\\v+v'+\frac{t'\sin u \sin (u+u')+(t+t'\cos u)(\cos(u+u')-1)}{x(1-\cos(u+u'))+\sin (u+u')}\end{bmatrix},
$$

where
$$
P_1=(x\cos (u+u')+\sin (u+u'))(t'\sinh u \sin (u+u')-(t+t'\cos u)( \cos (u+u')-1)).
$$
\newpage

\subsection{BOL ALGEBRAS WITH TRILINEAR OPERATION OF TYPE IV}

 As in the previous chapter we will base our investigation of
3-dimensional Bol algebras on the examination of their canonical enveloping
Lie algebras. In what follows, we consider Bol algebras of dimension 3, from
their construction see \cite{mike2,mike6}. Below we limit ourselves to the
classification
of Bol algebras (and their corresponding 3-Webs) with canonical enveloping Lie
algebras of dimension $ \leq 4$.

 Let $\mathfrak{B}$ be a 3-dimensional Bol algebra with a trilinear
operation Type IV see chapter II \S 3, and
$\mathfrak{G}=\mathfrak{B} \dotplus \mathfrak{h} $- its canonical enveloping
Lie algebra. According to the table given in chapter II(case 4). W note that
the situation $ dim \mathfrak{G}=3$ is impossible, that means the case of grouped 3-Webs is excluded.

 Let us examine the case   $ dim \mathfrak{G}=4$, the structural constants of Lie algebra  $ \mathfrak{G}=<e_1,e_2,e_3,e_4>$, $\mathfrak{B}=<e_1,e_2,e_3>$
are defined as follows:
$$
[e_1,e_2]=e_4, \; \; [e_2,e_4]= \mp e_1
$$
\; \; \; \; \; \; \; \; \; \; \; \; \; \; \; \; \; \; \; \; (1)
$$
[e_1,e_3]=\pm e_4, \; \; [e_3,e_4]=- e_1
$$
in addition $ \mathfrak{G}=\mathfrak{B} \dotplus [\mathfrak{B},\mathfrak{B}]$,
$[\mathfrak{B}, \mathfrak{B}]=<e_4>$.

By introducing in consideration the 3-dimensional subspaces of subalgebras

$$
\mathfrak{h}_{x,y,z}=<e_4 +xe_1 +ye_2 +ze_3>,\; \; x,y,z, \in \mathbb{R}
$$
we obtain a collection of Bol algebras of view:
$$
e_1 \cdot e_2=xe_1 +ye_2 +ze_3 ,
$$
$$
e_1 \cdot e_3= \mp xe_1 \mp ye_2  \mp ze_3 ,
$$
\; \; \; \; \; \; \; \; \; \; \; \; \; \; \; \; \; \; \; \; \; \; (2)
$$
(e_1 ,e_2 ,e_3)=\pm e_1, \; \; (e_1 ,e_3 ,e_2)=-e_1
$$
$$
(e_1 ,e_2 ,e_3)=e_1, \; \; (e_1 ,e_3 ,e_2)=\mp e_1.
$$
 Our main problem will be to give, an isomorphical and isotopical
classification of Bol algebras of view (2).

 For the full examination of this case, we will split it into cases
Type $ IV^-$ and Type $IV^+$, corresponding to the upper and the lower signs
of the formulas (1) and (2).

 The group of automorphisms $F$ of Lie triple system $\mathfrak{B}$
relatively to a fixed base $ e_1, e_2 e_3 $ from Type $IV^-$ is defined as
follows:

\begin{displaymath}
 F=\left\{P=\left(\begin{array}{ccc}
\alpha & 0 & 0 \\
0 & \pm1 & 0 \\
0 & 0 & \pm1 \\
\end{array}\right)\right\}. \; \; \; \; \; (3)
\end{displaymath}

 The extension of automorphism $A$ from $F$ to the automorphism of Lie
algebra $\mathfrak{G}$, transforming the subspace $\mathfrak{B}$ into itself
can be realized as follows:

$$
Ae_4=A[e_1,e_2]=[e_1,Ae_2]= \mp e_4,
$$
\; \; \; \; \; \;or \; \; \; \; \; \; \; \; \; \; \;\; \; \; \; \; \; \; \;  (4)
$$
\mp Ae_4=A[e_1,e_3]=[Ae_1,Ae_3]=- \alpha e_4.
$$
In addition
$$
A(e_4 +xe_1 +ye_2 +ze_3)=\pm \alpha e_4 +xe_1 \pm ye_2 \pm ze_3=\pm \alpha(e_4 \pm xe_1 +\frac{y}{\alpha}e_2 +\frac{z}{\alpha}e_3),
$$
that is
$$
A(\mathfrak{h}_{x,y,z})=\mathfrak{h}_{x',y',z'},
$$
where
$$
x'=\pm x,
$$
$$
y'=\frac{y}{\alpha}, \; \; \; \; \; \; \; \; \; \; \; \; (5)
$$
$$
z'=\frac{z}{\alpha}.
$$

If $ z \neq 0$ then at the expense of choice of  $\alpha $ one can make
$y'=p$, $x'=\pm x \geq 0$.

If $z=0$, then at the expense of choice of $\alpha$ one can make $y'=1$
and  $x'= \pm x \geq 0$ and $y'=1$, $z'=0$.

Thus within the isomorphism we obtain two families of Bol algebras.

\underline{{\bf Theorem $ IV^-$ .3.1}} Any Bol algebra of dimension 3, with the
trilinear
operation of Type $IV^-$ and the canonical enveloping Lie algebra of dimension 4, is isomorphic to one of Bol algebras below:
\begin{itemize}
\item $ IV^- .1. e_1 \cdot e_2 =xe_1 +pe_2 +e_3, \; (e_1,e_2,e_2)=e_1$, $(e_1,e_3,e_2)=-e_1 x \geq 0 $
$e_1 \cdot e_2 =-xe_1 -pe_2 -e_3, \; (e_1,e_2,e_3)=e_1$, $(e_1,e_3,e_3)=-e_1 x \geq 0 $ for any $p$,
\item $ IV^- .2. e_1 \cdot e_2 =xe_1 +pe_2 , \; (e_1,e_2,e_2)=e_1$, $(e_1,e_3,e_2)=-e_1  x \geq 0 $
$e_1 \cdot e_2 =-xe_1 -pe_2 , \; (e_1,e_2,e_3)=e_1$, $(e_1,e_3,e_3)=-e_1 x \geq 0 $ for any $p$.
\end{itemize}

the distinguished Bol algebras are not isomorphic among themselves.

Similarly one can establish the correctness of the following Theorem.

\underline{{\bf Theorem $ IV^+$ .3.1}} Any Bol algebra of dimension 3, with
the trilinear
operation of Type $IV^+$ and the canonical enveloping Lie algebra of dimension 4, is isomorphic to one of Bol algebras below:
\begin{itemize}
\item $ IV^+ .1. e_1 \cdot e_2 =xe_1 +pe_2 +e_3, \; (e_1,e_2,e_2)=-e_1$, $(e_1,e_3,e_2)=-e_1 x \geq 0 $
$e_1 \cdot e_2 =-xe_1 -pe_2 -e_3, \; (e_1,e_2,e_3)=e_1$, $(e_1,e_3,e_3)=e_1 x \geq 0 $ for any $p$,
\item $ IV^+ .2. e_1 \cdot e_2 =xe_1 +pe_2 , \; (e_1,e_2,e_2)=-e_1$, $(e_1,e_3,e_2)=-e_1 x \geq 0 $
$e_1 \cdot e_2 =-xe_1 -pe_2 , \; (e_1,e_2,e_3)=e_1$, $(e_1,e_3,e_3)=e_1 x \geq 0 $.
\end{itemize}

Also this distinguished Bol algebras are not isomorphic among themselves.

Let us  pass to the isotopic classification of Bol algebras given in
Theorems $IV^-$.3.1.

We note that for every $ \xi=ae_1 +be_2 +ce_3 $ from $ \mathfrak{B}$
$ a,b,c, \in \mathbb{R} $

\begin{displaymath}
ad ( \xi )= \left(\begin{array}{cccc}
0 & 0 & 0 & -b-c \\
0 & 0 & 0 & 0  \\
0 & 0 & 0 & 0 \\
-b+c & a & -a & 0 \\
\end{array}\right),
\end{displaymath}

\begin{displaymath}
Ad ( \xi )= \left(\begin{array}{cccc}
\cosh \sqrt{df} & \frac{a}{f}(-1+\cosh \sqrt{df}) & \frac{a}{f}(1-\cosh \sqrt{df} & \sqrt{\frac{d}{f}}\sinh \sqrt{df} \\
0 & 1 & 0 & 0  \\
0 & 0 & 1 & 0 \\
 \sqrt{\frac{d}{f}}\sinh \sqrt{df} &  \sqrt{\frac{a^2}{df}}\sinh \sqrt{df} &  \sqrt{\frac{a^2}{df}}\sinh \sqrt{df} & cosh \sqrt{df} \\
\end{array}\right)
\end{displaymath}
where $d=-b-c,f=-b+c$.

Let us fin the the image of $\Phi(\mathfrak{h})$ under the action of
$ \Phi=Ad \xi$ on the one-dimensional subalgebra $\mathfrak{h}$ with a
direction vector
$$
e_4 +xe_1 +pe_2 +ze_3, \; \;or\; \; e_4 +xe_1 +e_2.
$$
Moreover by the local character of the consideration we will limit ourselves
to the values $a,b$ and $c$ such that:

$$
\cosh \sqrt{df}-\frac{a}{\sqrt{fd}}+x\sqrt{\frac{f}{d}}\sinh \sqrt{fd}+\frac{pa}{\sqrt{df}}\sinh \sqrt{df} \neq 0,
$$
\begin{multline}
\Phi (e_4 +xe_1 +pe_2 +ze_3)=(\cosh \sqrt{df}-\frac{a}{\sqrt{fd}}+x\sqrt{\frac{f}{d}}\sinh \sqrt{fd}+\frac{pa}{\sqrt{df}}\sinh \sqrt{df})e_4+\\ +e_3 +pe_2 + (\sqrt{\frac{d}{f}}\sinh \sqrt{df}+x \cosh \sqrt{df}+p \frac{a}{f}(1+ \cosh \sqrt{fd})+\frac{a}{f}(1- \cosh \sqrt{df}))e_1
\end{multline}
\begin{multline}
\Phi (e_4 +xe_1 +e_2)=(\cosh \sqrt{df}+x\sqrt{\frac{f}{d}}\sinh \sqrt{fd}+\frac{a}{\sqrt{df}}\sinh \sqrt{df})e_4+\\+(\sqrt{\frac{d}{f}}\sinh \sqrt{df}+x \cosh \sqrt{df}+\frac{a}{f}(1+ \cosh \sqrt{fd}))e_1 +e_2
\end{multline}
Here we note that $e_2$ and $e_3$ are not changing only $e_4$ and $e_1$ are
changing.

We denote
$$
\widetilde{\alpha}=\cosh \sqrt{df}-\frac{a}{\sqrt{fd}}+x\sqrt{\frac{f}{d}}\sinh \sqrt{fd}+\frac{pa}{\sqrt{df}}\sinh \sqrt{df} \neq 0,
$$
$$
\widetilde{\beta}=\sqrt{\frac{d}{f}}\sinh \sqrt{df}+x \cosh \sqrt{df}+p \frac{a}{f}(1+ \cosh \sqrt{fd})+\frac{a}{f}(1- \cosh \sqrt{df}),
$$
$$
\Phi (\mathfrak{h})=\widetilde{\alpha}e_4 +\widetilde{\beta}e_1 +pe_2 +e_3=
$$
$$
=\widetilde{\alpha}\left(e_4 + \frac{\widetilde{\beta}}{\widetilde{\alpha}}\right)e_1 +\frac{p}{\widetilde{\alpha}}e_2 + \frac{1}{\widetilde{\alpha}}.
$$

After the application of the automorphism $A$ we obtain an expression of type:

$$
\pm \frac{\widetilde{\beta}}{\widetilde{\alpha}}e_1 +pe_2 +e_3
$$
The checking shows that $ \frac{\widetilde{\beta}}{\widetilde{\alpha}}$ are any
 values. By a nice choice of values $a,b,c$. One can note that,
the application of the isotopic transformation $\Phi$ to the Bol algebras of
Theorem $IV^-$ is not changing them. In what follows denote
$$
\delta =\frac{\widetilde{\beta}}{\widetilde{\alpha}}.
$$

 Summarizing the conducted examination one can formulate the Theorem:

\underline{{\bf Theorem $ IV^-$ .3.1}} Any Bol algebra of dimension 3, with the
 trilinear
operation of Type $IV^-$ and the canonical enveloping Lie algebra $\mathfrak{G}$ of dimension 4, is isotopic to one of the following Bol algebras:
\begin{itemize}
\item $  e_1 \cdot e_2 =\pm \delta e_1 +pe_2 +e_3, \; (e_1,e_2,e_2)=e_1$, \;$ (e_1,e_3,e_2)=-e_1, \delta \geq 0$
 $  e_1 \cdot e_3 =\mp \delta e_1 -pe_2 -e_3, \; (e_1,e_2,e_3)=e_1$, \;$ (e_1'e_3,e_2)=-e_1, \delta \geq 0$ where $p$ is any nomber,
\item $  e_1 \cdot e_2 =\pm \delta e_1 +e_2, \; (e_1,e_2,e_2)=e_1$, \;$ (e_1'e_3,e_2)=-e_1, \delta \geq 0,$
 $  e_1 \cdot e_3 =\mp \delta e_1 -e_2, \; (e_1,e_2,e_3)=e_1$, \;$ (e_1'e_3,e_2)=-e_1, \delta \geq 0$.
\end{itemize}

Analogically one can state the correctness of the Theorem:

\underline{{\bf Theorem $ IV^+$ .3.1}} Any Bol algebra of dimension 3, with
the trilinear
operation of Type $IV^+$ and the canonical enveloping Lie algebra $ \mathfrak{G}$ of dimension 4, is isotopic to one of the following Bol algebras:
\begin{itemize}
\item $  e_1 \cdot e_2 =\pm \delta e_1 +pe_2 +e_3, \; (e_1,e_2,e_2)=-e_1$, \;$ (e_1'e_3,e_2)=-e_1, \delta \geq 0$
 $  e_1 \cdot e_3 =\mp \delta e_1 -pe_2 -e_3, \; (e_1,e_2,e_3)=e_1$, \;$ (e_1'e_3,e_2)=e_1, \delta \geq 0$ where $p$ is any nomber,
\item $  e_1 \cdot e_2 =\pm \delta e_1 +e_2, \; (e_1,e_2,e_2)=-e_1$, \;$ (e_1'e_3,e_2)=-e_1, \delta \geq 0$
 $  e_1 \cdot e_3 =\mp \delta e_1 +e_2, \; (e_1,e_2,e_3)=e_1$, \;$ (e_1'e_3,e_2)=e_1$.
\end{itemize}

 Below we reduce to description of 3-Webs correspondimg to the isolated
Bol algebras of Type$ IV^-$ and Type$ IV^+$. The composition law $(\triangle)$,
corresponding to the Lie group $G$ of Lie algebra enveloping Bol algebra is
defined as follows:

$$ \begin{bmatrix} x_1\\x_2\\x_3\\x_4\end{bmatrix} \triangle
   \begin{bmatrix} y_1\\y_2\\y_3\\y_4\end{bmatrix}
   =
   \begin{bmatrix}
     x_{1}+y_{1}+ \cosh (x_{2}+x_3) -y_{4}sinh (x_{2}+x_3)\\x_{2}+y_{2}\\x_{3}+y_{3}\\x_{4}-y_{1}\sinh (x_{2}+x_3) +y_{4}\cosh (x_{2}+x_3)\end{bmatrix}.
$$

 In case $ IV^-.1$ the subgroup $ H=\exp \mathfrak{h} $ can be realized as the
collection of elements

$$
H=\exp \mathfrak{h}=\{\exp \alpha (e_4 +xe_1 +pe_2 +e_3)\}_{ \alpha \in \mathbb{R}}=\{x \alpha,p \alpha ,\alpha, \alpha \}_{ \alpha \in \mathbb{R}}.
$$
The collection of elements
$$
B= \exp \mathfrak{B}=\{\frac{t}{u+v}\sinh (u+v),u,v,\frac{t}{u+v}(1-cosh (u+v))\}_{t,u.v \in \mathbb{R}}
$$
form a local section of left space coset $G \bmod H $.
$$
\exp^{-1}  \begin{bmatrix} x_1\\x_2\\x_3\\x_4\end{bmatrix}= \begin{bmatrix} \frac{x_{1}(x_{2}+x_3)}{\sinh (x_{2}+x_3)}\\x_2\\x_3\end{bmatrix}
$$
$  x_1,x_2,x_3,x_4 \in \mathbb{R}$.

Any element $ (x_1,x_2, x_3, x_4) $ from $ G $, in the neighbourhood  $e$, can be
uniquely represented as follows:

\begin{multline}
 \begin{pmatrix} x_1\\x_2\\x_3\\x_4\end{pmatrix} =
   \begin{bmatrix}\frac{\left[x_1 (\cosh (u+v) -x \sinh (u+v))-x_4 (x\cosh (u+v)-\sinh (u+v)) \right]\sinh (u+v)}{x(1-\cosh (u+v))+\sinh (u+v)}\\u\\ v\\\frac{\left[x_1 (\cosh (u+v) -x \sinh (u+v))-x_4 (x\cosh (u+v)-\sinh (u+v)) \right](\cosh (u+v)-1)}{x(1-\cosh (u+v))+\sinh (u+v)}\end{bmatrix}
   \triangle\\
   \begin{bmatrix} x\frac{x_{4}\sinh (u+v)+x_{1}(\cosh (u+v)-1)}{x(1-\cosh (u+v))+\sinh (u+v)}\\ p\frac{x_{4}\sinh (u+v)+x_{1}(\cosh (u+v)-1)}{x(1-\cosh (u+v))+\sinh (u+v)}\\ \frac{x_{4}\sinh (u+v)+x_{1}(\cosh (u+v)-1)}{x(1-\cosh (u+v))+\sinh (u+v)}\ \\ \frac{x_{4}\sinh (u+v)+x_{1}(\cosh (u+v)-1)}{x(1-\cosh (u+v))+\sinh (u+v)}\end{bmatrix},
\end{multline}

where $u+v$ are any numbers defined from the relation
$$
u+v+(p+1) \frac{x_{1}( \cosh (u+v)-1)+x_{4}\sinh (u+v)}{x(1-\cosh (u+v))+\sinh (u+v)}=x_{2}+x_{3}.
$$

The composition law ($ \star$ ),  corresponding to the local analytical Bol
loop $B( \star )$, is defined as follows:

$$ \begin{pmatrix} t\\u\\v\end{pmatrix} \star
   \begin{pmatrix}t'\\u'\\v'\end{pmatrix}
   = \exp^{-1} \left(\prod_{B}\left(
    \begin{bmatrix}t\\u\\v\\0\end{bmatrix}\triangle \begin{bmatrix}
     t'\\u'\\v'\\0\end{bmatrix}\right) \right)
$$
$$ \begin{pmatrix} t\\u\\v\end{pmatrix} \star
   \begin{pmatrix}t'\\u'\\v'\end{pmatrix}
   = exp^{-1} \left(\prod_{B}\left(
    \begin{bmatrix}t+t'\cosh (u+v)\\u+u'\\v+v'\\-t\sinh (u+v)\end{bmatrix}\right) \right)
$$
$$
= \exp^{-1} \left( \begin{bmatrix}T \sinh (m+n)\\m\\n\\T(\cosh (m+n) -1)\end{bmatrix}\right)
$$.

$$
= \begin{bmatrix}T(m+n)\\m\\n\end{bmatrix},
$$

where $T$ is defined from the relation:
$$
T=\frac{D}{x(1-\cosh(m+n) +\sinh(m+n))},
$$
with
\begin{multline}
D=(t+t' \cosh (u+v))(\cosh (m+n)-x\sinh (m+n))+\\+t\sinh (u+v) (x\cosh (m+n) -\sinh (m+n)),
\end{multline}

and $m+n$ from the relation

$$
m+n+(p+1)\frac{(t+t' \cosh (u+v))(\cosh (m+n) -1)-t\sinh (u+v) \sinh (m+n)}{x(1-\cosh (m+n))+\sinh (m+n)}=\Omega
$$
and $\Omega=u+u'+v+v'$.

 In case $ IV^- .2 $ the subgroup $ H=\exp \mathfrak{h} $, can be realized as the
collection of elements

$$
H=\exp \mathfrak{h}=\{exp \alpha (e_4 +xe_1 +e_2 )\}_{ \alpha \in \mathbb{R}}=\{x \alpha, \alpha ,0, \alpha \}_{ \alpha \in \mathbb{R}}.
$$
The collection of elements
$$
B= exp \mathfrak{B}=\{\frac{t}{u+v}\sinh (u+v),u,v,\frac{t}{u+v}(1-cosh (u+v))\}_{t,u.v \in \mathbb{R}}
$$
form a local section of left space coset $G \bmod H $.
$$
\exp^{-1}  \begin{bmatrix} x_1\\x_2\\x_3\\x_4\end{bmatrix}= \begin{bmatrix} \frac{x_{1}(x_{2}+x_3)}{\sinh (x_{2}+x_3)}\\x_2\\x_3\end{bmatrix}
$$
$  x_1,x_2,x_3,x_4 \in \mathbb{R}$.

Any element $ (x_1,x_2, x_3, x_4) $ from $ G $, in the neighbourhood  $e$, can be
uniquely represented as follows:

\begin{multline}
 \begin{pmatrix} x_1\\x_2\\x_3\\x_4\end{pmatrix} =
   \begin{bmatrix}\frac{[x_1 (\cosh (u+x_3) -x \sinh (u+x_3))-x_4 (x\cosh (u+x_3)-\sinh (u+x_3)) ]\sinh (u+x_3)}{x(1-\cosh (u+x_3))+\sinh (u+x_3)}\\u\\ x_3\\\frac{[x_1 (\cosh (u+x_3) -x \sinh (u+x_3))-x_4 (x\cosh (u+x_3)-\sinh (u+x_3)) ](\cosh (u+x_3)-1)}{x(1-\cosh (u+x_3))+\sinh (u+x_3)}\end{bmatrix}
   \triangle\\
   \begin{bmatrix} x\frac{x_{4}\sinh (u+x_3)+x_{1}(\cosh (u+x_3)-1)}{x(1-\cosh (u+x_3))+\sinh (u+x_3)}\\\frac{x_{4}\sinh (u+x_3)+x_{1}(\cosh (u+x_3)-1)}{x(1-\cosh (u+x_3))+\sinh (u+x_3)}\\0\\ \frac{x_{4}\sinh (u+x_3)+x_{1}(\cosh (u+x_3)-1)}{x(1-\cosh (u+x_3))+\sinh (u+x_3)}\end{bmatrix},
\end{multline}

where $u$ are any numbers defined from the relation
$$
u+x_3 +(p+1) \frac{x_{1}( \cosh (u+x_3)-1)+x_{4}\sinh (u+x_3)}{x(1-\cosh (u+x_3))+\sinh (u+x_3)}=x_{2}+x_{3}.
$$

\; \; The composition law $( \star )$, corresponding to the local analytical Bol
loop $B( \star )$, is defined as follows:

$$ \begin{pmatrix} t\\u\\v\end{pmatrix} \star
   \begin{pmatrix}t'\\u'\\v'\end{pmatrix}
   = \exp^-1 \left(\prod_{B}\left(
    \begin{bmatrix}t\\u\\v\\0\end{bmatrix}\triangle \begin{bmatrix}
     t'\\u'\\v'\\0\end{bmatrix}\right) \right)
$$
$$ \begin{pmatrix} t\\u\\v\end{pmatrix} \star
   \begin{pmatrix}t'\\u'\\v'\end{pmatrix}
   = \exp^{-1} \left(\prod_{B}\left(
    \begin{bmatrix}t+t'\cosh (u+v)\\u+u'\\v+v'\\-t\sinh (u+v)\end{bmatrix}\right) \right)
$$
$$
= \exp^{-1} \left( \begin{bmatrix}\left[H \sinh (m+V+V')\right]\\m\\V+V'\\H(\cosh (m+v+v') -1)\end{bmatrix}\right)
$$.

$$
= \begin{bmatrix}H(m+v+v')\\m\\v+v'\end{bmatrix},
$$

where $H$ is defined from the relation:
$$
H=\frac{A}{x(1-\cosh(m+v+v') +\sinh(m+v+v'))},
$$
\begin{multline}
A=(t+t' \cosh (u+v))(\cosh (m+v+v')-x\sinh (m+v+v'))+\\t\sinh (u+v+v') (x\cosh (m+v+v') -\sinh (m+v+v')),
\end{multline}
and by denoting $X=m+v+v' $, $m$ will be defined from the relation
$$
X+\frac{\Gamma}{x(1-\cosh (m+v+v'))+\sinh (m+v+v')}=u+u',
$$

$\Gamma=(t+t' \cosh (u+v+v'))(\cosh (m+v+v') -1)-t\sinh (u+v+v') \sinh (m+v+v')$.

 We pass to the description of Bol 3-Webs corresponding to the Type $IV^+$.

 The composition law $(\triangle)$, corresponding to Lie group $G$, with
enveloping Lie algebra of Bol algebra of Type $IV^+$ is defined as:

$$ \begin{bmatrix} x_1\\x_2\\x_3\\x_4\end{bmatrix} \triangle
   \begin{bmatrix} y_1\\y_2\\y_3\\y_4\end{bmatrix}
   =
   \begin{bmatrix}
     x_{1}+y_{1}+ \cos (x_{2}-x_3) -y_{4}sin (x_{2}-x_3)\\x_{2}+y_{2}\\x_{3}+y_{3}\\x_{4}-y_{1}\sin (x_{2}-x_3) +y_{4}\cos (x_{2}-x_3)\end{bmatrix}.
$$

 In case $ IV^+.1$ the subgroup $ H=\exp \mathfrak{h} $, can be realized as the
collection of elements

$$
H=\exp \mathfrak{h}=\{exp \alpha (e_4 +xe_1 +pe_2 +e_3)\}_{ \alpha \in \mathbb{R}}=\{x \alpha,p \alpha ,\alpha, \alpha \}_{ \alpha \in \mathbb{R}}.
$$
The collection of elements
$$
B= exp \mathfrak{B}=\{\frac{t}{u+v}\sin (u+v),u,v,\frac{t}{u+v}(1-cos (u+v))\}_{t,u.v \in \mathbb{R}}
$$

form a local section of left space coset $G \bmod H $.

$$
\exp^{-1}  \begin{bmatrix} x_1\\x_2\\x_3\\x_4\end{bmatrix}= \begin{bmatrix} \frac{x_{1}(x_{2}-x_3)}{\sin (x_{2}-x_3)}\\x_2\\x_3\end{bmatrix}.
$$
$  x_1,x_2,x_3,x_4 \in \mathbb{R}$

Any element $ (x_1,x_2, x_3, x_4) $ from $ G $, in the neighbourhood  $e$, can be
uniquely represented as follows:

\begin{multline}
 \begin{pmatrix} x_1\\x_2\\x_3\\x_4\end{pmatrix} =
   \begin{bmatrix}\frac{\left[x_1 (\cos (u-v) -x \sin (u-v))-x_4 (x\cos (u-v)-\sin (u-v)) \right]\sin (u-v)}{x(1-\cos (u-v))+\sin (u-v)}\\u\\ v\\\frac{\left[x_1 (\cos (u-v) -x \sin (u-v))-x_4 (x\cos (u-v)-\sin (u-v)) \right](\cos (u-v)-1)}{x(1-\cos (u-v))+\sin (u-v)}\end{bmatrix}\\
   \triangle
   \begin{bmatrix} x\frac{x_{4}\sin (u-v)+x_{1}(\cos (u-v)-1)}{x(1-\cos (u-v))+\sin (u-v)}\\ p\frac{x_{4}\sin (u-v)+x_{1}(\cos (u-v)-1)}{x(1-\cos (u-v))+\sin (u-v)}\\ \frac{x_{4}\sin (u-v)+x_{1}(\cos (u-v)-1)}{x(1-\cos (u-v))+\sin (u-v)}\ \\ \frac{x_{4}\sin (u-v)+x_{1}(\cos (u-v)-1)}{x(1-\cos (u-v))+\sin (u-v)}\end{bmatrix},
\end{multline}

where $u-v$ are any numbers defined from the relation
$$
u-v+(p+1) \frac{x_{1}( \cos (u-v)-1)+x_{4}\sin (u-v)}{x(1-\cos (u-v))+\sin (u-v)}=x_{2}+x_{3}.
$$

\; \; The composition law $( \star )$, corresponding to the local analytical Bol
loop $B( \star )$, is defined as follows:

$$ \begin{pmatrix} t\\u\\v\end{pmatrix} \star
   \begin{pmatrix}t'\\u'\\v'\end{pmatrix}
   = exp^-1 \left(\prod_{B}\left(
    \begin{bmatrix}t\\u\\v\\0\end{bmatrix}\triangle \begin{bmatrix}
     t'\\u'\\v'\\0\end{bmatrix}\right) \right)
$$
$$ \begin{pmatrix} t\\u\\v\end{pmatrix} \star
   \begin{pmatrix}t'\\u'\\v'\end{pmatrix}
   = \exp^{-1} \left(\prod_{B}\left(
    \begin{bmatrix}t+t'\cos (u-v)\\u+u'\\v+v'\\-t\sin (u-v)\end{bmatrix}\right) \right)
$$
$$
= \exp^{-1} \left( \begin{bmatrix}\left[K \sin (m-n)\right]\\m\\n\\K(\cos (m-n) -1)\end{bmatrix}\right)
$$.

$$
= \begin{bmatrix}K(m-n)\\m\\n\end{bmatrix},
$$

where $K$ is defined from the relation:

$$
K=\frac{\Delta}{x(1-\cos(m-n) +\sin (m-n))},
$$

\begin{multline}
\Delta=[(t+t' \cos (u-v))(\cos (m-n)-x\sin (m-n))+\\+t\sin (u-v) (x\cos (m-n) -\sin (m-n))],
\end{multline}

and $m-n$ from the relation

$$
m-n+(p+1)\frac{Y}{x(1-\cos (m-n))+\sin (m-n)}=u+u'+v+v',
$$

$Y=(t+t' \cos (u-v))(\cos (m-n) -1)-t\sin (u-v) \sin (m-n)$.

 In case $ IV^+.2.$ The subgroup $ H=\exp \mathfrak{h} $, can be realized as the
collection of elements

$$
H=\exp \mathfrak{h}=\{exp \alpha (e_4 +xe_1 +e_2 )\}_{ \alpha \in \mathbb{R}}=\{x \alpha, \alpha ,0, \alpha \}_{ \alpha \in \mathbb{R}}.
$$
The collection of elements
$$
B= \exp \mathfrak{B}=\{\frac{t}{u-v}\sin (u-v),u,v,\frac{t}{u-v}(1-cos (u-v))\}_{t,u.v \in \mathbb{R}}
$$
form a local section of left space coset $G \bmod H $.
$$
\exp^{-1}  \begin{bmatrix} x_1\\x_2\\x_3\\x_4\end{bmatrix}= \begin{bmatrix} \frac{x_{1}(x_{2}-x_3)}{\sinh (x_{2}-x_3)}\\x_2\\x_3\end{bmatrix}
$$
$  x_1,x_2,x_3,x_4 \in \mathbb{R}$.

Any element $ (x_1,x_2, x_3, x_4) $ from $ G $, in the neighbourhood  $e$, can be
uniquely represented as follows:

\begin{multline}
 \begin{pmatrix} x_1\\x_2\\x_3\\x_4\end{pmatrix} =
   \begin{bmatrix}\frac{\left[x_1 (\cos (u-x_3) -x \sin (u-x_3))-x_4 (x\cos (u-x_3)-\sin (u-x_3)) \right]\sin (u-x_3)}{x(1-\cos (u-x_3))+\sin (u-x_3)}\\u\\ x_3\\\frac{\left[x_1 (\cos (u-x_3) -x \sin (u-x_3))-x_4 (x\cos (u-x_3)-\sin (u-x_3)) \right](\cos (u-x_3)-1)}{x(1-\cos (u-x_3))+\sin (u-x_3)}\end{bmatrix}
   \triangle\\
   \begin{bmatrix} x\frac{x_{4}\sin (u-x_3)+x_{1}(\cos (u-x_3)-1)}{x(1-\cos (u-x_3))+\sin (u-x_3)}\\\frac{x_{4}\sin (u-x_3)+x_{1}(\cos (u-x_3)-1)}{x(1-\cos (u-x_3))+\sin (u-x_3)}\\0\\ \frac{x_{4}\sin (u-x_3)+x_{1}(\cos (u-x_3)-1)}{x(1-\cos (u-x_3))+\sin (u-x_3)}\end{bmatrix},
\end{multline}

where $u$ are any numbers defined from the relation
$$
u+x_3 +(p+1) \frac{x_{1}( \cos (u-x_3)-1)+x_{4}\sin (u-x_3)}{x(1-\cos (u-x_3))+\sin (u-x_3)}=x_{2}+x_{3}.
$$

\; \; The composition law $( \star )$ corresponding to the local analytical Bol
loop $B( \star )$ is defined as follows:

$$ \begin{pmatrix} t\\u\\v\end{pmatrix} \star
   \begin{pmatrix}t'\\u'\\v'\end{pmatrix}
   = \exp^{-1} \left(\prod_{B}\left(
    \begin{bmatrix}t\\u\\v\\0\end{bmatrix}\triangle \begin{bmatrix}
     t'\\u'\\v'\\0\end{bmatrix}\right) \right)
$$
$$ \begin{pmatrix} t\\u\\v\end{pmatrix} \star
   \begin{pmatrix}t'\\u'\\v'\end{pmatrix}
   = \exp^{-1} \left(\prod_{B}\left(
    \begin{bmatrix}t+t'\cos (u-v)\\u+u'\\v+v'\\-t\sinh (u-v)\end{bmatrix}\right) \right)
$$
$$
= \exp^{-1} \left( \begin{bmatrix}\left[S \sin (l-V-V')\right]\\l\\v+v'\\S(\cos (l-v-v') -1)\end{bmatrix}\right)
$$.

$$
= \begin{bmatrix}S(l-v-v')\\l\\v+v'\end{bmatrix},
$$

where $S$ is defined from the relation:
$$
S=\frac{W}{x(1-\cos (l-v-v') +\sin (l-v-v'))}
$$
with
\begin{multline}
W=(t+t' \cos (u-v))(\cos (l-v-v')-x\sin (l-v-v'))+\\
+t\sin (u-v-v') (x\cos (l-v-v') -\sin (l-v-v')),
\end{multline}
and $l$ from the relation
$$
l+v+v'+\frac{Z}{x(1-\cos (l-v-v'))+\sin (l-v-v')}=u+u',
$$

$Z=(t+t' \cos (u-v-v'))(\cos (l-v-v') -1)-t\sin (u-v-v') \sin (l-v-v').$

\newpage

\subsection{BOL ALGEBRA WITH TRILINEAR OPERATION OF TYPE V}
 In what follows we consider Bol algebras of dimension 3, from their
construction see [30,34] in this paragraph we base our investigation of
3-dimensional Bol algebras, on the examination of their canonical enveloping
Lie algebras. As we already state it follows that the dimension of their
canonical enveloping Lie algebras can not be more than 6. Below we limit
ourselves to the classification of Bol algebras ( and their corresponding
3-Webs), with canonical enveloping Lie algebras of dimension $ \leq 4$.

 Let us examine the case   $ dim \mathfrak{G}=4$, the structural constants of
 Lie algebra  $ \mathfrak{G}=<e_1,e_2,e_3,e_4>$, $\mathfrak{B}=<e_1,e_2,e_3>$
are defined as follows:
$$
[e_2,e_3]=e_4, \; \; [e_2,e_4]=- e_1
$$
\; \; \; \; \; \; \; \; \; \; \; \; \; \; \; \; \; \; \; \; \; \; \; \; \; \; (1)
$$
[e_3,e_4]=\mp e_4.
$$
in addition $ \mathfrak{G}=\mathfrak{B} \dotplus [\mathfrak{B},\mathfrak{B}]$,
$[\mathfrak{B}, \mathfrak{B}]=<e_4>$.

By introducing in consideration the 3-dimensional subspaces of subalgebras
$$
\mathfrak{h}_{x,y,z}=<e_4 +xe_1 +ye_2 +ze_3>,\; \; x,y,z, \in \mathbb{R}
$$
we obtain a collection of Bol algebras of view:
$$
e_2 \cdot e_3=-xe_1 -ye_2 -ze_3 ,
$$
\; \; \; \; \; \; \; \; \; \; \; \; \; \; \;\;\;\;\;\;\; \; \; \; \; \; \; \; (2)
$$
(e_2 ,e_3 ,e_2)=e_1,\;
$$
$$
 (e_2 ,e_3 ,e_3)=\pm e_1.
$$
\; Our main problem will be to give an isomorphical and isotopical
classification of Bol algebras of view (2).

\; For the full examination of this case, we will split it in cases
Type $ V^-$ and Type $V^+$, corresponding to the upper and the lower signs
of the formulas (1) and (2).

\; The group of automorphisms $F$ of Lie triple system $\mathfrak{B}$
relatively to a fixed base $ e_1, e_2, e_3 $ from Type $IV^-$ is defined as
follows:

\begin{displaymath}
 F=\left\{A=\left(\begin{array}{ccc}
\pm bf^2 & fb & d \\
0 & b & f \\
0 & 0 & \pm 1 \\
\end{array}\right), b \neq 0 \right\}\; \; \;\;\;\;\;\;\;\;\;(3)
\end{displaymath}
\; The extension of automorphism $A$ from $F$ to the automorphism of Lie
algebra $\mathfrak{G}$, transforming the subspace $\mathfrak{B}$ into itself
can be realized as follows:

$$
Ae_4=A[e_2,e_3]=[Ae_2,Ae_3]= \pm be_4 . \; \; \;\;\;\;\;\;\; \; (4)
$$

In addition
$$
A(e_4 +xe_1 +ye_2 +ze_3)=\pm b \left(e_4  \pm \frac{yfb+zd \pm xbf^2}{b}e_1 \pm \frac{yb+zf}{b}e_2 +\frac{z}{b}e_3 \right),
$$
that is
$$
A(\mathfrak{h}_{x,y,z})=\mathfrak{h}_{x',y',z'},
$$
where
$$
x'= \frac{yfb+zd \pm xbf^2}{b},
$$
$$
y'=\frac{yb+zf}{b}, \; \; \; \; \; \;\;\;\;\;\;\;\;\;(5)
$$
$$
z'=\frac{z}{b}.
$$

\begin{itemize}
\item If $ z \neq 0$, then by choosing $b,f$ and $d$ one can make
$x'=0,y'=0,z'=0$;
\item if $z=0$, hence $z'=0$ one can choose $y'= \pm y \geq 0$, and make
$x'=0$.
\end{itemize}
\; \; In this way, we obtain one family and one exceptional Bol algebra.

\underline{{\bf Theorem $ V^-$ .5.1}} Any Bol algebra of dimension 3, with the
trilinear
operation of Type $V^-$ and the canonical enveloping Lie algebra of dimension
4, is isomorphic to one of Bol algebras below:
\begin{itemize}
\item $ V^- .1. e_2 \cdot e_3 =-e_3, \; (e_2,e_3,e_2)=e_1$, $(e_2,e_3,e_3)=e_2, $
\item $ V^- .2. e_2 \cdot e_3 =-ye_2 , \; (e_2,e_3,e_2)=e_1, (e_2,e_3,e_3)=e_2,
 \;  y \geq 0 $.
\end{itemize}

The distinguished Bol algebras are not isomorphic among themselves.

Similarly one can establish the correctness of the following Theorem.

\underline{{\bf Theorem $ V^+$ .3.1}} Any Bol algebra of dimension 3, with the
trilinear
operation of Type $V^+$ and the canonical enveloping Lie algebra of dimension
4, is isomorphic to one of Bol algebras below:
\begin{itemize}
\item $ V^+ .1. e_2 \cdot e_3 =-e_3, \; (e_2,e_3,e_2)=e_1$, $(e_2,e_3,e_3)=e_2 ,$
\item $ V^+ .2. e_2 \cdot e_3 =-ye_2, \; (e_2,e_3,e_2)=e_1, (e_2,e_3,e_3)=-e_2
, \;y \geq 0. $
\end{itemize}

Also this distinguished Bol algebras are not isomorphic among themselves.

Let us  pass to the isotopic classification of Bol algebras given in
Theorems $V^-$.3.1.

We note that for every $ \xi=ue_1 +ve_2 +pe_3 $ from $ \mathfrak{B}$
$ u,v,p, \in \mathbb{R} $

\begin{displaymath}
ad ( \xi )= \left(\begin{array}{cccc}
0 & 0 & 0 & -v \\
0 & 0 & 0 & -p  \\
0 & 0 & 0 & 0 \\
0 & -p & v & 0 \\
\end{array}\right),
\end{displaymath}

\begin{displaymath}
Ad ( \xi )= \left(\begin{array}{cccc}
1 & \frac{v(\cosh p -1)}{p} &- \frac{(\cosh p -1)v^2}{p^2} & -\frac{v \sinh p}{p} \\
0 & \cosh p & \frac{v(\cosh p -1)}{p} & -\sinh p  \\
0 & 0 & 1 & 0 \\
0 & -\sinh p & \frac{v \sinh p}{p} & \cosh p \\
\end{array}\right).
\end{displaymath}
Let us find the the image of $\Phi(\mathfrak{h})$ under the action of
$ \Phi=Ad \xi$ on the one-dimensional subalgebra $\mathfrak{h}$ with a
direction vector
$e_4  +ye_2,\; y \geq 0$;

$$
\Phi (e_4  +ye_2)=(\cosh p -y\sinh p)e_4 +\frac{v}{p}[y(cosh p -1)-\sinh p]e_1 +(ycosh p -\sinh p)e_2,
$$

$$
\Phi (e_4  +ye_2)=e_4 +\frac{v}{p}[\frac{y(cosh p -1)-\sinh p}{\cosh p -y\sinh p}]e_1 +\frac{ycosh p -\sinh p}{\cosh p -y\sinh p}e_2,
$$
$$
x'=\frac{v}{p}[\frac{y(cosh p -1)-\sinh p}{\cosh p -y\sinh p}],
\; \; \; \; \; \; \; \; \; \; \; \; \; \; \; \; \; \; \; \; \; \; \;(5)
$$
$$
y'=+\frac{ycosh p -\sinh p}{\cosh p -y\sinh p}.
$$
\; \; By choosing $p \neq 0$ such that $y=\frac{\sinh p}{\cosh p -1}$,
we obtain $x'=0$. Let us note in addition that $coth p \neq 0$ that is the
map is correctly defined. Applying to the obtain Bol algebra the
automorphism of view 5, one can make $y'=1$.
\begin{itemize}
\item If $p =0$ then $x'=-v$\\
\;\;\;\;\;\;\;\;\; $y'=y$; \\
\item if $v \neq 0$ then one can make $ x'=1, y'=1$;\\
\item if $v=0$ then $x'=0$ and, one can make $y'=1$;\\
\item if $y=0$ then $x'=-\frac{v}{p}tanh p$,\\
\; \; \; \; \; \; \; \; \; \; \; \; \; \; \; \; \; \; \; \;
\; \; \; \; \; \; \; \; (6)

$y'=-tanh p$,
\end{itemize}
\begin{itemize}
\item a) if $p=0$, then $x'=-v, v'=0$
\item b) if $p \neq 0$ then when applying the automorphism of view (6) and with
regards to $v$ we can make $x'=0; 1 , y'=1$.
\end{itemize}
\; \; The selected three cases are isotopic (in the sense of the definition
of isotopy).

 We note the exceptional Bol algebra of Theorem III.5.1 under the
action of isotopic transformation is not changing

\; \; Summarizing the conducted examination one can formulate the Theorem:

\underline{{\bf Theorem $ V^-$ .3.3}} Any Bol algebra of dimension 3, with the
trilinear
operation of Type $V^-$ and the canonical enveloping Lie algebra$ \mathfrak{G}$ of dimension 4, is isotopic to one of the following Bol algebras:
\begin{itemize}
\item $  e_2 \cdot e_3 =-e_3, \; (e_2,e_3,e_2)=e_1, (e_2,e_3,e_3)=e_2$;
\item $ e_2 \cdot e_3 =-e_2, \; (e_2,e_3,e_2)=e_1, \; (e_2,e_3,e_3)=e_2$;
\item $  e_2 \cdot e_3 =-e_1 -e_2, \; (e_2,e_3,e_2)=e_1, \; (e_2,e_3,e_3)=e_2$;
\item trivial bilinear operation, $(e_2,e_3,e_2)=e_1, \; (e_2,e_3,e_3)=e_2$.
\end{itemize}

Analogically one can state the correctness of the Theorem.

\underline{{\bf Theorem $ V^+$ .3.4}} Any Bol algebra of dimension 3, with the
trilinear
operation of Type $V^+$ and the canonical enveloping Lie algebra$ \mathfrak{G}$ of dimension 4, is isotopic to one of the following Bol algebras:
\begin{itemize}
\item $  e_2 \cdot e_3 =-e_3, \; (e_2,e_3,e_2)=e_1, (e_2,e_3,e_3)=-e_2$;
\item $ e_2 \cdot e_3 =e_2, \; (e_2,e_3,e_2)=e_1, \; (e_2,e_3,e_3)=-e_2$;
\item trivial bilinear operation, $(e_2,e_3,e_2)=e_1, \; (e_2,e_3,e_3)=-e_2$.
\end{itemize}

\; \; Below we reduce to description of 3-Webs corresponding to the isolated
Bol algebras of Type$ V^-$ and Type$ V^+$.

The composition law $(\triangle)$,
corresponding to the Lie group $G$ of enveloping Lie algebra for Bol algebra is
defined as follows:

$$ \begin{bmatrix} x_1\\x_2\\x_3\\x_4\end{bmatrix} \triangle
   \begin{bmatrix} y_1\\y_2\\y_3\\y_4\end{bmatrix}
   =
   \begin{bmatrix}
     x_{1}+y_{1}+ \frac{x_{4}y_{2} -y_{4}x_{2}}{2}\\x_{2}+y_{2}\cos x_{3}-y_{4}\sin x{3}\\x_{3}+y_{3}\\x_{4}-y_{2}\sin (x_3) +y_{4}\cos (x_3)\end{bmatrix}.
$$

 In case $ V^-.1$ the subgroup $ H=\exp \mathfrak{h} $, can be realized as the
collection of elements

$$
H=\exp \mathfrak{h}=\{exp \alpha (e_4  +e_3)\}_{ \alpha \in \mathbb{R}}=\{0,0,\alpha, \alpha \}_{ \alpha \in \mathbb{R}}.
$$
The collection of elements
$$
B=\exp \mathfrak{B}=\left\{t+\frac{(v-\sin v)u^2}{2v^2},\frac{u}{v}\sin v,v,\frac{u}{v}(1-cos v)\right\}_{t,u,v \in \mathbb{R}}
$$
form a local section of left space coset $G \bmod H $.
$\exp: \mathfrak{G}\supset\mathfrak{B} \longrightarrow B \subset G$
and
$B =\exp \mathfrak{B}$
$$
\exp^{-1} \begin{bmatrix} x_1\\x_2\\x_3\\x_4\end{bmatrix}= \begin{bmatrix} x_1-\frac{(x_{2})^2 \sin^2 (x_3)}{2(x_3)^3}+\frac{(x_{2})^2
\sin^3 (x_3)}{2(x_3)^4}\\\frac{x_2}{x_3}\sin x_3\\x_3\\x_4\end{bmatrix}
$$
$  x_1,x_2,x_3,x_4 \in \mathbb{R}$.

Any element $ (x_1,x_2, x_3, x_4) $ from $ G $, in the neighborhood  $e$, can be
uniquely represented as follows:

$$ \begin{pmatrix} x_1\\x_2\\x_3\\x_4\end{pmatrix} =
   \begin{bmatrix}x_1 +\frac{(x_3 -v)[x_2 +(x_3 -v)\sin v ]}{2}\\x_2 +(x_3 -v)\sin v\\ v\\x_4 -(x_3 -v)\cos v\end{bmatrix}
   \triangle
   \begin{bmatrix} 0\\0\\x_3 -v\\ x_3 -v\end{bmatrix},
$$

where $v$ are any numbers defined from the relation
$$
\left[ x_4 -(x_3 -v)\cos v \right]\sin v =\left[x_2 +(x_3 -v)\cos v\right](\cos v -1).
$$

\; \; The composition law $( \star )$ corresponding to the local analytical Bol
loop $B( \star )$ is defined as follows:

$$ \begin{pmatrix} t\\u\\v\end{pmatrix} \star
   \begin{pmatrix}t'\\u'\\v'\end{pmatrix}
   = \exp^{-1} \left(\prod_{B}\left(
    \begin{bmatrix}t\\u\\v\\0\end{bmatrix}\triangle \begin{bmatrix}
     t'\\u'\\v'\\0\end{bmatrix}\right) \right)
$$
$$ \begin{pmatrix} t\\u\\v\end{pmatrix} \star
   \begin{pmatrix}t'\\u'\\v'\end{pmatrix}
   = \exp^{-1} \left(\prod_{B}\left(
    \begin{bmatrix}t+t'\\u+u'\cos v\\v+v'\\u \sin (v)\end{bmatrix}\right) \right)
$$
$$
= \exp^{-1} \left( \begin{bmatrix}t+t'+\frac{(v+v'-T)\left[u+u'\cos v +(v+v'-T)\sin T\right]}{2}\\u+u'\cos v +(v+v'-T)\sin T\\T\\u\sin v -(v+v'-T)\cos T\end{bmatrix}\right)
$$.

$$
= \begin{bmatrix}F_1(t,t',u'u',v'v',T)\\\left[u+u'\cos v +(v+v'-T)sin T\right]\frac{\sin T}{T}\\T\end{bmatrix},
$$

where $T$ is defined from the relation:
$$
\left[ u\sin v -(v+v' -T)\cos T \right]\sin T =\left[u+u'\cos v +(v+v' -T)\cos T\right](\cos T -1).
$$

And $F_1(t,t',u'u',v'v',T)$ from the relation
\begin{multline}
F_1(t,t',u'u',v'v',T)=t+t'+\frac{(v+v'-T)(u+u'\cos v)}{2}+\frac{(v+v'-T)^2}{2}\sin T-\\
-\frac{\left[u+u'cos v + (v+v'-T)\sin T \right]^2}{2T^4}(T-\sin T)\sin^2 T.
\end{multline}

In  case $ V^-.2$ the subgroup $ H=\exp \mathfrak{h} $, can be realized as the
collection of elements

$$
H=\exp \mathfrak{h}=\{exp \alpha (e_4 +ye_2 )\}_{ \alpha \in \mathbb{R}}=\{0,0, \alpha , \alpha \}_{ \alpha \in \mathbb{R}}.
$$
The collection of elements
$$
B= exp \mathfrak{B}=\left\{t+\frac{(v-\sin v)u^2}{2v^2},\frac{u}{v}\sin v,v,\frac{u}{v}(1-cos v)\right\}_{t,u.v \in \mathbb{R}}
$$
form a local section of left space coset $G \bmod H $.

Here $\exp^{-1}$  is defined as in the case above.

\; \; Any element $ (x_1,x_2, x_3, x_4) $ from $ G $, in the neighborhood
 $e$, can be uniquely represented as follows:

$$ \begin{pmatrix} x_1\\x_2\\x_3\\x_4\end{pmatrix} =
   \begin{bmatrix}x_1 AB\frac{y-y\cos x_3 -\sin x_3}{ 2x_3}\\A\frac{\sin x_3}{x_3}\\ x_3\\A\frac{1-\cos x_3}{x_3}\end{bmatrix}
   \triangle
   \begin{bmatrix} 0\\y\frac{x_{4}\sin x_3 -x_{2}(1-\cos x_3)}{y-y\cos x_3 +\sin x_3}\\0\\\frac{x_{4}\sin x_3 -x_{2}(1-\cos x_3)}{y-y\cos x_3 +\sin x_3}\end{bmatrix},
$$

where $A,B$ are any numbers defined from the relations
$$
A=x_3 \frac{x_2 (y\sin x_3 +\cos x_3)-x_4 (y \cos x_3 -\sin x_3)}{y-y\cos x_3 + \sin x_3},
$$
$$
B= \frac{x_4 \sin x_3 -x_2 (1- \cos x_3)}{y-y\cos x_3 + \sin x_3}.
$$

\; \; The composition law $( \star )$ corresponding to the local analytical Bol
loop $B( \star )$ is defined as follows:

$$ \begin{pmatrix} t\\u\\v\end{pmatrix} \star
   \begin{pmatrix}t'\\u'\\v'\end{pmatrix}
   = \exp^{-1} \left(\prod_{B}\left(
    \begin{bmatrix}t\\u\\v\\0\end{bmatrix}\triangle \begin{bmatrix}
     t'\\u'\\v'\\0\end{bmatrix}\right) \right)
$$
$$ \begin{pmatrix} t\\u\\v\end{pmatrix} \star
   \begin{pmatrix}t'\\u'\\v'\end{pmatrix}
   = \exp^{-1} \left(\prod_{B}\left(
    \begin{bmatrix}t+t'\\u+u'\cos v\\v+v'\\-u\sin (v)\end{bmatrix}\right) \right)
$$

$$
= \begin{bmatrix}t+t'-C'\\A'\frac{\sin^2 (v+v')}{(v+v')^2}\\v+v'\end{bmatrix},
$$

where $A',C'$ are defined from the relations:
$$
A=(v+v') \frac{\left\{(u+u'\cos v)\left[ (y\sin (v+v') +\cos (v+v')\right]-u\sin v [y \cos (v+v') -\sin (v+v')]\right\}}{y-y\cos (v+v') + \sin (v+v')},
$$

$$
B= \frac{u\sin v \sin (v+v') -(u+u'\cos v) (1- \cos (v+v'))}{y-y\cos (v+v') + \sin (v+v')},
$$
$$
C'=A'B'\frac{(y-y\cos (v+v') -\sin (v+v'))}{v+v'}+(A')^2 \frac{\sin^4 (v+v')}{2(v+v')^4} \cdot \frac{-1+\cos (v+v')}{v+v'}.
$$

\; \; We pass to the description of Bol 3-Webs corresponding to the Type $V^+$.

\; \; In this case the composition law $(\triangle)$, corresponding to Lie group $G$, with
enveloping Lie algebra of Bol algebra of Type $V^+$ is defined as:

$$ \begin{bmatrix} x_1\\x_2\\x_3\\x_4\end{bmatrix} \triangle
   \begin{bmatrix} y_1\\y_2\\y_3\\y_4\end{bmatrix}
   =
   \begin{bmatrix}
     x_{1}+y_{1}+ \frac{x_{4}y_{2} -y_{4}x_{2}}{2}\\x_{2}+y_{2}\cosh x_{3}-y_{4}\sinh x{3}\\x_{3}+y_{3}\\x_{4}-y_{2}\sinh (x_3) +y_{4}\cosh (x_3)\end{bmatrix}.
$$

 In case $ V^+.1$ the subgroup $ H=exp \mathfrak{h} $, can be realized as the
collection of elements

$$
H=\exp \mathfrak{h}=\{exp \alpha (e_4  +e_3)\}_{ \alpha \in \mathbb{R}}=\{0,0,\alpha, \alpha \}_{ \alpha \in \mathbb{R}}.
$$
The collection of elements
$$
B= \exp \mathfrak{B}=\left\{t+\frac{(v-\sinh v)u^2}{2v^2},\frac{u}{v}\sinh v,v,\frac{u}{v}(1-cosh v)\right\}_{t,u.v \in \mathbb{R}}
$$
form a local section of left space coset $G \bmod H $.
$\exp: \mathfrak{G}\supset\mathfrak{B} \longrightarrow B \subset G$
and
$B =\exp \mathfrak{B}$
$$
\exp^{-1}  \begin{bmatrix} x_1\\x_2\\x_3\\x_4\end{bmatrix}= \begin{bmatrix} x_1-\frac{(x_{2})^2 \sinh^2 (x_3)}{2(x_3)^3}+\frac{(x_{2})^2 \sinh^3 (x_3)}{2(x_3)^4}\\\frac{x_2}{x_3}\sinh x_3\\x_3\end{bmatrix}
$$
$  x_1,x_2,x_3,x_4 \in \mathbb{R}$.

Any element $ (x_1,x_2, x_3, x_4) $ from $ G $, in the neighborhood  $e$, can be
uniquely represented as follows:

$$ \begin{pmatrix} x_1\\x_2\\x_3\\x_4\end{pmatrix} =
   \begin{bmatrix}x_1 +\frac{(x_3 -v)[x_2 +(x_3 -v)\sinh v]}{2}\\x_2 +(x_3 -v)\sin v\\ v\\x_4 -(x_3 -v)\cos v\end{bmatrix}
   \triangle
   \begin{bmatrix} 0\\0\\x_3 -v\\ x_3 -v\end{bmatrix}
$$

where $v$ are any numbers defined from the relation
$$
[ x_4 -(x_3 -v)\cosh v ]\sinh v =[x_2 +(x_3 -v)\cosh v](\cosh v -1).
$$

\; \; The composition law $( \star )$, corresponding to the local analytical Bol
loop $B( \star )$, is defined as follows:

$$ \begin{pmatrix} t\\u\\v\end{pmatrix} \star
   \begin{pmatrix}t'\\u'\\v'\end{pmatrix}
   = \exp^{-1} \left(\prod_{B}\left(
    \begin{bmatrix}t\\u\\v\\0\end{bmatrix}\triangle \begin{bmatrix}
     t'\\u'\\v'\\0\end{bmatrix}\right) \right)
$$
$$ \begin{pmatrix} t\\u\\v\end{pmatrix} \star
   \begin{pmatrix}t'\\u'\\v'\end{pmatrix}
   = \exp^{-1} \left(\prod_{B}\left(
    \begin{bmatrix}t+t'\\u+u'\cosh v\\v+v'\\u \sinh (v)\end{bmatrix}\right) \right)
$$
$$
= \exp^{-1} \left( \begin{bmatrix}t+t'+\frac{(v+v'-P)\left[u+u'\cosh v +(v+v'-P)\sinh P\right]}{2}\\u+u'\cosh v +(v+v'-P)\sin P\\P\\u\sinh v -(v+v'-P)\cosh P\end{bmatrix}\right)
$$.

$$
= \begin{bmatrix}F_2(t,t',u'u',v'v',P)\\\left[u+u'\cosh v +(v+v'-P)sinh P\right]\frac{\sinh P}{P}\\P\end{bmatrix}
$$

where $P$ is defined from the relation:
$$
\left[ u\sinh v -(v+v' -P)\cosh P \right]\sinh P =\left[u+u'\cosh v +(v+v' -P)\cosh P\right](\cosh P -1)
$$

and $F_1(t,t',u'u',v'v',P)_1$ from the relation
\begin{multline}
F_1(t,t',u'u',v'v',P)=t+t'+\frac{(v+v'-P)(u+u'\cosh v)}{2}+\frac{(v+v'-P)^2}{2}\sinh P-\\
-\frac{\left[u+u'cosh v + (v+v'-P)\sinh P \right]^2}{2P^4}(P-\sinh P)\sinh^2 P.
\end{multline}

In case $ V^+.2$ the subgroup $ H=\exp \mathfrak{h} $, can be realized as the
collection of elements

$$
H=\exp \mathfrak{h}=\{\exp \alpha (e_4 +ye_2 )\}_{ \alpha \in \mathbb{R}}=\{0,0, \alpha , \alpha \}_{ \alpha \in \mathbb{R}}.
$$
The collection of elements
$$
B= \exp \mathfrak{B}=\left\{t+\frac{(v-\sinh v)u^2}{2v^2},\frac{u}{v}\sinh v,v,\frac{u}{v}(1-cosh v)\right\}_{t,u.v \in \mathbb{R}}
$$
form a local section of left space coset $G mod H $.

Here $\exp^{-1}$  is defined as in the case above.

\; \; Any element $ (x_1,x_2, x_3, x_4) $ from $ G $, in the neighborhood
 $e$, can be uniquely represented as follows:

$$ \begin{pmatrix} x_1\\x_2\\x_3\\x_4\end{pmatrix} =
   \begin{bmatrix}x_1 DE\frac{y-y\cosh x_3 -\sinh x_3}{ 2x_3}\\D\frac{\sinh x_3}{x_3}\\ x_3\\D\frac{1-\cosh x_3}{x_3}\end{bmatrix}
   \triangle
   \begin{bmatrix} 0\\y\frac{x_{4}\sin x_3 -x_{2}(1-\cosh x_3)}{y-y\cosh x_3 +\sinh x_3}\\0\\\frac{x_{4}\sinh x_3 -x_{2}(1-\cosh x_3)}{y-y\cosh x_3 +\sinh x_3}\end{bmatrix},
$$

where $D,E$ are any numbers defined from the relations
$$
D=x_3 \frac{[x_2 (y\sinh x_3 +\cosh x_3)-x_4 (y \cosh x_3 -\sinh x_3)]}{y-y\cosh x_3 + \sinh x_3},
$$
$$
E= \frac{x_4 \sinh x_3 -x_2 (1- \cosh x_3)}{y-y\cosh x_3 + \sinh x_3}.
$$

\; \; The composition law $( \star )$ corresponding to the local analytical Bol
loop $B( \star )$ is defined as follows:

$$ \begin{pmatrix} t\\u\\v\end{pmatrix} \star
   \begin{pmatrix}t'\\u'\\v'\end{pmatrix}
   = \exp^{-1} \left(\prod_{B}\left(
    \begin{bmatrix}t\\u\\v\\0\end{bmatrix}\triangle \begin{bmatrix}
     t'\\u'\\v'\\0\end{bmatrix}\right) \right)
$$
$$ \begin{pmatrix} t\\u\\v\end{pmatrix} \star
   \begin{pmatrix}t'\\u'\\v'\end{pmatrix}
   = \exp^{-1} \left(\prod_{B}\left(
    \begin{bmatrix}t+t'\\u+u'\cosh v\\v+v'\\-u\sinh (v)\end{bmatrix}\right) \right)
$$

$$
= \begin{bmatrix}t+t'-F'\\D'\frac{\sinh^2 (v+v')}{(v+v')^2}\\v+v'\end{bmatrix}
$$

where $D',F'$ are defined from the relations:
$$
D'=(v+v') \frac{\Lambda}{y-y\cosh (v+v') + \sinh (v+v')},
$$

\begin{multline}
\Lambda=\{(u+u'\cosh v)\left[ (y\sinh (v+v') +\cosh (v+v')\right]-\\-u\sinh v \left[y \cosh (v+v') -\sinh (v+v')\right]\},
\end{multline}
$$
E'= \frac{u\sinh v \sinh (v+v') -(u+u'\cosh v) (1- \cosh (v+v'))}{y-y\cosh (v+v') + \sinh (v+v')},
$$
$$
F'=D'E'\frac{(y-y\cosh (v+v') -\sinh (v+v'))}{v+v'}+(A')^2 \frac{\sinh^4 (v+v')}{2(v+v')^4} \cdot \frac{-1+\cosh (v+v')}{v+v'}.
$$
\newpage

\subsection{BOL ALGEBRAS WITH TRILINEAR OPERATION OF TYPE VI}

\; \; As in the previous chapter, we will base our investigation of
3-dimensional Bol algebras, on the examination of their canonical enveloping
Lie algebras. In what follows, we consider Bol algebras of dimension 3, from
their construction see [30,34]; it follows that the dimension of their
canonical enveloping Lie algebras, can not be more than 6. Below we limit
ourselves to the classification of Bol algebras (and their corresponding
3-webs), with canonical enveloping Lie algebras of dimension$ \leq 5$.

Let $ \mathfrak{B} $ be a 3-dimensional Bol algebra with a trilinear
operation Type VI see chapter II \S 3, and
$\mathfrak{G}=\mathfrak{B} \dotplus \mathfrak{h}$-its canonical enveloping Lie
 algebra according to the table given in chapter II (case 4). We note that the
 situation
$dim \mathfrak{G}=3$ is possible, that means we obtain a total grouped 3-Web the
corresponding Lie group $G$, is isomorphic to the matrix of the view:

\begin{displaymath}
 \left(\begin{array}{ccc}
\cosh x & -\sinh x & \frac{z(1-\cosh x)+y\sinh x}{x} \\
\sinh x & -\cosh x & \frac{z(1-\cosh x)+z\sinh x}{x} \\
0 & 0 & 1 \\
\end{array}\right),
\end{displaymath}
with $ x,y,z \in \mathbb{R}$.

\; \; Let us examine the case $dim \mathfrak{G}=4$, the structural constants of Lie algebra $ \mathfrak{G}=<e_1,e_2,e_3,e_4>$, $\mathfrak{B}=<e_1,e_2,e_3>$
are defined as follows:
$$
[e_2,e_3]=e_4, \; \; [e_3,e_4]=- e_1, \; \; \; \; \; (1)
$$
$$
[e_1,e_3]=-e_5, [e_3,e_5]=-e_2,
$$
in addition $ \mathfrak{G}=\mathfrak{B} \dotplus [\mathfrak{B},\mathfrak{B}]$,
$[\mathfrak{B}, \mathfrak{B}]=<e_4>$.

By introducing in consideration the 3-dimensional subspaces of subalgebras
$$
\mathfrak{h}_{x,y,z,x_1,y_1,z_1}=<e_4 +xe_1 +ye_2 +ze_3,e_5 +x_1 e_1 +y_1 e_2 +z_1 e_3>,
$$
where $ x,y,z, x_1,y_1,z_1 \in \mathbb{R}$
\; \;  Let us note
$$
e'_4=e_4 +xe_1 +ye_2 +ze_3,
$$
$$
e'_5=e_5 +x_1 e_1 +y_1 e_2 +z_1 e_3. \; \; \; \; \; \; (2)
$$
We are interested in those spaces which are Lie subalgebras in $\mathfrak{G}$,
that is $[e'_4,e'_5] \in <e'_4,e'_5>$.

\; \; \; The condition of algebraical closure of the subspaces $\mathfrak{h}_{x,y,z,x_1,y_1,z_1}$
consists that there exist $\alpha, \beta \in \mathbb{R}$ such that:

$$
[e'_4,e'_5]=\alpha e'_4+\beta e'_5,
$$
that is
$$
\beta= x_1 z-z_1 x, \alpha=yz_1 -zy_1, \; \; (\ast)
$$
$$
z_1=\alpha x+\beta x_1, \alpha z+ \beta z_1=0, \; \; \; \; (3)
$$
$$
-z=\alpha y + \beta y_1.
$$
\begin{enumerate}
\item 1. If $\beta =0$ then, $\alpha \neq 0$ or $\alpha=0$.
\begin{itemize}
\item If $\alpha \neq 0$, then $z=y=0$ and we come to the contradiction with the
condition $(\ast)$.
\item $\alpha=0$, then $z=z_1=0$ and, the subspace will be defined as follows:
$$
<e_4 +xe_1 +ye_2, e_5 +x_1 e_1 +y_1 e_2>,
$$
it is a subalgebra we note it $\mathfrak{h}_{x,y,0,x_1,y_1,0}$.
\end{itemize}
\item 2. If $\beta \neq 0$ assume that:
\begin{itemize}
\item 2.1.  $x_1=y_1=0, z_1=t>0$, then $\beta =tx, \alpha =yt$. And, associating it
with (3), we obtain:
\begin{displaymath}
\left\{\begin{array}{l}
t=tyx,\\-z=ty^2\\-yt-xt^2 \end{array}\right.
\end{displaymath}
then $xy=1$, hence $x \neq 0, y \neq 0, z neq 0$;
\begin{displaymath}
\left\{\begin{array}{l}
yx=1,\\x=-y^3 \end{array}\right.
\end{displaymath},
 hence $y^4=-1$, and we obtain a contradiction
\item  2.2 $x=y=0$, thus $\beta =tx_1, \alpha =-iy_1, z=p>0,z_1=px^2_1$,
$-p=\beta y_1=x_1 y_1 p$ hence,
\begin{displaymath}
\left\{\begin{array}{l}
y_1 \cdot x_1=-1,\\(x_1)^4=-1 \end{array}\right.
\end{displaymath}
and we also come to the contradiction.
\item  2.3 $x_1=z=0$, then $\beta =-xz_1, \alpha =yz_1, \beta z_1=0$, hence $z=0$
we obtain a contradiction.
\item 2.4 $x=z_1=0$, then $\beta =-zx_1, \alpha=-zy_1, \alpha z=0$.
\begin{itemize}
\item 2.4.1 If $\alpha =0$ then $z=0$, and come to a contradiction.
\item 2.4.2 If $\alpha=0$, then $-zy_1=0$, hence $y_1=0$ hence, $z_1=0$ also
we obtain a contradiction.
\end{itemize}
\item 2.5 $y=y_1=0,$ then $ \alpha=z=0$, $\beta=-xz_1, \beta z_1=0$, hence,
$z_1=0$ also a contradiction.
\item 2.6 $y=z=0$, then $\beta =zx_1, \alpha =-zy_1, \beta z_1=0$, hence,
$x_1=0$ we still get a contradiction.
\item 2.7 $y=z=0$, then $\alpha z=0$; here $z \neq 0$, otherwise it will
contradict the condition $\beta \neq 0$; hence $\alpha =0, \beta =zx_1$,
$\alpha y =-z=0$ this is also a contradiction.
\item 2.8 $y_1=z=0$, then $z_1=0$ it's a contradiction.
\item 2.9 $y=z_1=0$, $ \alpha=-zy_1, \beta=zx_1, z \neq 0, \alpha z =0$, hence,
$\alpha=0, y_1=0$, this implies $z=0$ it's a contradiction.
\end{itemize}
\end{enumerate}

Summarizing: the condition $\beta \neq 0$ its impossible that is $\beta=\alpha=0$.

The analogical reasoning can be done if we consider $\alpha=0$. In result we
obtain only the subalgebra:
$$
<e_4 +xe_1 +ye_2, e_5 +x_1 e_1 +y_1 e_2>.
$$
\; \;  The group of automorphisms $F$ of Lie triple system $\mathfrak{B}$
relatively to a fixed base $ e_1, e_2 e_3 $ from Type $IV^-$ is defined as
follows:

\begin{displaymath}
 F=\left\{P=\left(\begin{array}{ccc}
\alpha & -\beta & d \\
\beta & \alpha & f \\
0 & 0 & \pm 1 \\
\end{array}\right), \alpha^2 +\beta^2 \neq 0 \alpha,\beta, d,f \in \mathbb{R}\right\}.
\end{displaymath}
\; \; \; The extension of automorphism $P$ from $F$ to the automorphism of Lie
algebra $\mathfrak{G}$, transforming the subspace $\mathfrak{B}$ into itself
can be realized as follows:

$$
Ae_4=[Ae_2,Ae_3]= \pm \beta e_5 \pm \alpha e_4 ,
$$
$$
Ae_5=[Ae_1,Ae_3]= \pm \alpha e_5 \pm \beta e_4.
$$
In addition for the subspace  $\mathfrak{h}_{x,y,0,x_1,y_1,0}$ we have:
$$
A(e_4 +xe_1 +ye_2 )=e_4  \pm \frac{x\alpha^2 +y \beta^2-y \alpha \beta -\beta \alpha x}{\alpha^2 +\beta^2}e_1 \pm \frac{y\alpha^2 -x \beta^2+x\alpha \beta -y \alpha \beta}{\alpha^2 +\beta^2}e_2,
$$
\; \; \; \; \; \; \; \; \; \; \; \; \; \; \; \; \; \; \; \; \; \; \;\; (4)
$$
A(e_5 +x_1 e_1 +y_1 e_2 )=e_5  \pm \frac{x_1 \alpha^2 +y_1 \beta^2-y_1 \alpha \beta +\beta \alpha x_1}{\alpha^2 +\beta^2}e_1 \pm \frac{y_1\alpha^2 +x_1 \beta^2+x\alpha \beta +y_1 \alpha \beta}{\alpha^2 +\beta^2}e_2,
$$
such that the action of $A$ on the subspace  $\mathfrak{h}_{x,y,0,x_1,y_1,0}$,
can be represented as follows:

$$
\begin{pmatrix} \frac{\alpha}{\sqrt{\alpha^2+\beta^2}}&-\frac{\beta}{\sqrt{\alpha^2+\beta^2}}\\\frac{\beta}{\sqrt{\alpha^2+\beta^2}}&\frac{\alpha}{\sqrt{\alpha^2+\beta^2}} \end{pmatrix}
\begin{pmatrix} x&x_1\\y&y_1 \end{pmatrix}
\begin{pmatrix} \frac{\alpha}{\sqrt{\alpha^2+\beta^2}}&\frac{\beta}{\sqrt{\alpha^2+\beta^2}}\\-\frac{\beta}{\sqrt{\alpha^2+\beta^2}}&\frac{\alpha}{\sqrt{\alpha^2+\beta^2}} \end{pmatrix}.
$$

If the matrix
\begin{displaymath}
 \left(\begin{array}{cc}
\lambda_1 & \delta  \\
0 & \lambda_2 \\
\end{array}\right)
\end{displaymath},
with $ \lambda_1, \delta \in \mathbb{R}$. \; \;  \; \; \;  \; \; \; \; (6)

If the matrix
\begin{displaymath}
 \left(\begin{array}{cc}
x & x_1  \\
y & y_1 \\
\end{array}\right),
\end{displaymath}
has two complex conjugate eigenvalues
$\lambda_1 =\alpha +i\beta (\lambda_2 =\alpha -i\beta)$ correspondently, where
$\beta \neq 0$, and eigenvectors
$\bar{a}=\bar{\xi}+i\bar{\eta}(\bar{a}=\bar{\xi}-i\bar{\eta)}$ correspondently.

\; \; \;  Let $\eta \neq 0$ the first vector base
$(e_1,e_2)\longrightarrow (\bar{\eta}, \bar{\epsilon})$, where $\bar{\epsilon}$
is perpendicular to $\bar{\eta}$
$$
A\bar{a}=A\bar{\xi}+iA\bar{\eta}=(\alpha \bar{\xi}-\beta \bar{\eta})+i(\alpha \bar{\eta}+\beta \bar{\xi}),
$$
$$
A\bar{\eta}=\alpha \bar{\eta}+\beta \bar{\xi}, \; \; \; \; \; \; (7),
$$
where $ \bar{\xi}=\mu \bar{\eta}+\nu \bar{\epsilon}$
$$
A\bar{\xi}=\alpha \bar{\xi}-\beta \bar{\eta}=\mu A\bar{\eta}+\nu A\bar{\epsilon}.
$$

\; \; Thus  we can write out $A\bar{\eta}=(\alpha+\beta \mu)\bar{\eta}+\beta \mu \bar{\epsilon}$,

on the other hand

$$
\nu A\bar{\epsilon}=-\mu A \bar{\eta}+\alpha \bar{\xi}-\beta \bar{\eta}=\bar{\eta}(-\beta-\mu^2 \beta)+\bar{\epsilon}(-\mu \beta \nu+\alpha \nu), \;\nu \neq 0
.$$

\; \; \; Hence the canonical view of the matrix will be represented as follows:

\begin{displaymath}
 \left(\begin{array}{cc}
\alpha+\beta \mu & -\frac{\beta}{\nu}(1+\mu^2)  \\
\beta \nu & \alpha-\mu \beta \\
\end{array}\right)
\nu \neq 0, \beta \neq 0 .  \; \; \; \; (8)
\end{displaymath}

\; \; Alternative approach: if eigenvalues $\lambda_1=\alpha+i\beta (\lambda_2=\alpha-i\beta)$
are two complex conjugate numbers corresponding to two complex eigenvectors
$\bar{a}=\bar{\xi}+i\bar{\eta} (\bar{b}=\bar{\xi}-i\bar{\eta})$. Let's
consider vectors $\bar{\nu}, \bar{\mu}$ defined in terms of $\bar{a}$ and
$ \bar{b}$ as follow:
$$
\bar{\nu}=\frac{\bar{a}+\bar{b}}{2}, \bar{\nu}=\frac{\bar{a}-\bar{b}}{2i}.\;\;(9)
$$
The are reals, moreover we have:
$$
A\bar{\nu}=\alpha \bar{\nu}-\beta \bar{\mu}
$$
\; \; \; \; \; \; \; \; \; \; \; \; \; \; \;  \; \; \; \; \; \; \; \; \;\;\;\;\;\; \; (10)
$$
A\bar{\mu}=\beta \bar{\nu}+\alpha \bar{\mu}
$$

\; \; \;  Therefore the linear span in the space of reals numbers $\mathbb{R}$,
constructed in terms of vectors in (9), is an invariant subspace under the
action of automorphism $A$. Therefore the matrix induced by the automorphism in that
subspace, in the base defined in (9)is:

\begin{displaymath}
 \left(\begin{array}{cc}
\alpha & \beta  \\
-\beta & \alpha \\
\end{array}\right)
\beta \neq 0, \alpha \in \mathbb{R}.
\end{displaymath}

\; \; \;  We then obtain the following theorem:

\underline{{\bf Theorem III .6.1}} Any Bol algebra of dimension 3, with the
trilinear
operation of Type $III^-$, and the canonical enveloping Lie algebra of dimension 5, is isomorphic to one of Bol algebras below:
\begin{itemize}
\item $ VI .1. e_2 \cdot e_3 =-\lambda_1 e_1, \; (e_2,e_3,e_3)=e_1$,
\item $  e_1 \cdot e_3 =-\lambda_2 e_1-\delta e_2, \; (e_3,e_1,e_3)=e_2$ ,$ \lambda_1 \neq 0, \lambda_2 \neq 0 \delta \in \mathbb{R}$.
\end{itemize}
\begin{itemize}
\item $ VI .2. e_2 \cdot e_3 =-(\alpha +\beta \mu)e_1 -\beta \nu e_2, \; (e_2,e_3,e_3)=e_1$ $\nu > 0, \beta > 0 $,
\item $ e_1 \cdot e_3 =\frac{\beta (1+\mu^2)}{\nu}e_1-(\alpha-\mu)e_2, \;(e_3,e_1,e_3)=e_2$ $\mu ,\alpha \geq 0 , \nu neq 0$,
or \;\; \; \; \; \; \; \; \; \; \; \; \; \\; \; \; \; \; \; \; \; \; \; \; \;
$ e_2 \cdot e_3 =-\alpha e_1 +\beta e_2, e_1 \cdot e_3 =-\beta e_1 -\alpha e_2, \alpha \geq 0, \beta >0$.
\end{itemize}

We note that the distinguished Bol algebras are not isomorphic among themselves.

\; \; \;  The composition law $(\triangle)$, of local Lie group $G$, is defined as
follows:

$$ \begin{bmatrix} x_1\\x_2\\x_3\\x_4\\x_5\end{bmatrix} \triangle
   \begin{bmatrix} y_1\\y_2\\y_3\\y_4\\y_5\end{bmatrix}
   =
   \begin{bmatrix}
     F_1(x_{1},x_{3},y_{1},y_{2},y_{4},y_{5})\\F_2(x_{2},x_{3}, y_{1},y_{2},y_{4},y_{5})\\x_{3}+y_{3}\\F_4(x_{4},x_{3},y_1,y_{2},y_4,y_5)\\ F_5(x_5,x_3,y_1,y_2,y_4,y_5)\end{bmatrix}.
$$

 \; \; \;  Where
\begin{multline}
 F_1(x_{1},x_{3},y_{1},y_{2},y_{4},y_{5})=x_1 +(\cosh (x_3 \frac{\sqrt{2}}{2})\cos (x_3 \frac{\sqrt{2}}{2}))y_1 +(\sinh (x_3 \frac{\sqrt{2}}{2})\sin (x_3 \frac{\sqrt{2}}{2}))y_2\\
-\frac{\sqrt{2}}{2}\left[\sinh (x_3 \frac{\sqrt{2}}{2})\cos (x_3 \frac{\sqrt{2}}{2}) +\cosh (x_3 \frac{\sqrt{2}}{2})\sin (x_3 \frac{\sqrt{2}}{2})\right]y_4\\\
-\frac{\sqrt{2}}{2}\left[\sinh (x_3 \frac{\sqrt{2}}{2})\cos (x_3 \frac{\sqrt{2}}{2}) +\cosh (x_3 \frac{\sqrt{2}}{2})\sin (x_3 \frac{\sqrt{2}}{2})\right]y_5,
\end{multline}
\begin{multline}
F_2(x_{2},x_{3}, y_{1},y_{2},y_{4},y_{5})=x_2 +\frac{\sqrt{2}}{2}[\cos (x_3 \frac{\sqrt{2}}{2})\sinh (x_3 \frac{\sqrt{2}}{2}) -\sin (x_3 \frac{\sqrt{2}}{2})\cosh (x_3 \frac{\sqrt{2}}{2})]y_1\\
+(\cosh (x_3 \frac{\sqrt{2}}{2})\cos (x_3 \frac{\sqrt{2}}{2}))y_2 -(\sinh (x_3 \frac{\sqrt{2}}{2})\sin (x_3 \frac{\sqrt{2}}{2}))y_4\\
+\frac{\sqrt{2}}{2}[\cosh (x_3 \frac{\sqrt{2}}{2})\sin (x_3 \frac{\sqrt{2}}{2}) +\sinh (x_3 \frac{\sqrt{2}}{2})\cos (x_3 \frac{\sqrt{2}}{2})]y_5,
\end{multline}
\begin{multline}
F_4(x_{4},x_{3},y_1,y_{2},y_4,y_5)=x_4 -(\sinh (x_3 \frac{\sqrt{2}}{2})\sin (x_3 \frac{\sqrt{2}}{2}))y_1 \\
-\frac{\sqrt{2}}{2}\left[\sinh (x_3 \frac{\sqrt{2}}{2})\cos (x_3 \frac{\sqrt{2}}{2}) +\cosh (x_3 \frac{\sqrt{2}}{2})\sin (x_3 \frac{\sqrt{2}}{2})\right]y_2+\\+
(\cosh (x_3 \frac{\sqrt{2}}{2})\cos (x_3 \frac{\sqrt{2}}{2}))y_4-\\
-\frac{\sqrt{2}}{2}\left[\sinh (x_3 \frac{\sqrt{2}}{2})\cos (x_3 \frac{\sqrt{2}}{2}) +\cosh (x_3 \frac{\sqrt{2}}{2})\sin (x_3 \frac{\sqrt{2}}{2})\right]y_5,
\end{multline}
\begin{multline}
 F_5(x_5,x_3,y_1,y_2,y_4,y_5)= x_5 +\frac{\sqrt{2}}{2}\left[\sinh (x_3 \frac{\sqrt{2}}{2})\cos (x_3 \frac{\sqrt{2}}{2}) +\sin (x_3 \frac{\sqrt{2}}{2})\cosh (x_3 \frac{\sqrt{2}}{2})\right]y_1+\\+
(\sinh (x_3 \frac{\sqrt{2}}{2})\sin (x_3 \frac{\sqrt{2}}{2}))y_2-\\-
-\frac{\sqrt{2}}{2}\left[\sinh (x_3 \frac{\sqrt{2}}{2})\cos (x_3 \frac{\sqrt{2}}{2}) +\sin (x_3 \frac{\sqrt{2}}{2})\cosh (x_3 \frac{\sqrt{2}}{2})\right]y_4+\\+
(\cosh (x_3 \frac{\sqrt{2}}{2})\cos (x_3 \frac{\sqrt{2}}{2}))y_5
\end{multline}

\; \; \;  For this particular type of Bol algebras the classification with
accuracy to isotopy and the description of the corresponding 3-Webs are not
given for reason of awkwardness.

\newpage
\subsection{BOL ALGEBRAS WITH TRILINEAR OPERATION OF TYPE VII}

\; \; As in the last chapter we will base our investigation of
3-dimensional Bol algebras on the examination of their canonical enveloping
Lie algebras. Here we limit ourselves to the classification of Bol algebras
(and their corresponding
3-webs), with canonical enveloping Lie algebras of dimension$ \leq 5$.

Let $ \mathfrak{B} $ be a 3-dimensional Bol algebra with a trilinear
operation Type VII see chapter II \S 3, and
$\mathfrak{G}=\mathfrak{B} \dotplus \mathfrak{h}$-its canonical enveloping Lie
 algebra, according to the table given in chapter II (case 4). We note that the
 situation
$dim \mathfrak{G}=3$ is impossible that mean we obtain a total grouped 3-Web
 is excluded.

\; \;  The structural constants of Lie algebra $ \mathfrak{G}=<e_1,e_2,e_3,e_4>$, $\mathfrak{B}=<e_1,e_2,e_3>$
are defined as follows:
$$
[e_2,e_3]=e_4, \; \; [e_1,e_4]=- e_1,
$$
$$
[e_1,e_3]=-e_5, [e_2,e_5]=e_1, [e_4,e_5]=e_5,
$$
in addition $ \mathfrak{G}=\mathfrak{B} \dotplus [\mathfrak{B},\mathfrak{B}]$,
$[\mathfrak{B}, \mathfrak{B}]=<e_4>$.

By introducing in consideration the 3-dimensional subspaces of subalgebras
$\mathfrak{h}_{x,y,z,x',y',z'>}$ such that:
$$
\mathfrak{h}_{x,y,z,x',y',z'}=<e_4 +xe_1 +ye_2 +ze_3,e_5 +x' e_1 +y' e_2 +z' e_3>_{ x,y,z, x',y',z' \in \mathbb{R}}.
$$
\; \; \;  Let us note
$$
e'_4=e_4 +xe_1 +ye_2 +ze_3,
$$
$$
e'_5=e_5 +x_1 e_1 +y_1 e_2 +z_1 e_3.
$$
Our interest  goes to those spaces which are Lie subalgebras in $\mathfrak{G}$,
that is $[e'_4,e'_5] \in <e'_4,e'_5>$.

where
$$
[e'_4,e'_5]=(x'+y)e_1 +(yz'-zy')e_4 +(1+xz'-zy')e_5,
$$
so that the matrix

\begin{displaymath}
 \left(\begin{array}{ccccc}
x & y & z & 1 & 0 \\
x' & y' & z' & 0 & 1 \\
x'+y & 0 & 0 & yz'-zy' & 1+xz'-zx' \\
\end{array}\right),
\end{displaymath}

must be linearly dependent or what equivalent to:

\begin{enumerate}
\item (1)  $(x'+y)(yz'-y'z)=0$,
\item (2)  $(yz'-y'z)^2 -y'(x'+y)=0$,
\item (3)  $ y'(x'+y)+(1+xz'-zx')(xy'-x'y)=0$,
\item (4)  $(yz'-y'z)^2 =0$,
\item (5)  $(yz'-y'z)(1+xz'-zx')=0$,
\item (6)  $(yz'-y'z)(xz'-x'z) -z'(x'+y)=0$,
\item (7)  $z(x'+y)+(1+xz'-zx')(xz'-x'z)=0$,
\item (8)  $(x'+y)-x'(1+xz'-zx')-x(yz'-y'z)=0$,
\item (9)  $y(yz'-zy')-y'(1+xz'-zx')=0$,
\item (10) $z(yz'-zy')-z'(1+xz'-zx')=0$.
\end{enumerate}

The system of relation (1)--(10) equivalent to the following:
\begin{itemize}
\item (2)'  $y'(x'+y)=0$,
\item (3)'  $y(x'+y)+(1+xz'-zx')(xy'-x'y)=0$,
\item (6)'  $z'(x'+y)=0$,
\item (8)'  $(x'+y)-x'(1+xz'-zx')=0$,
\item (9)'  $y'(1+xz'-zx')=0$,
\item (10)' $z'(1+xz'-zx')=0$,
\end{itemize}
\begin{enumerate}
\item I. case 1. Let $x'+y \neq 0$, then $y'=z'=0$ and $y=-(x')^2 z$. in result
we obtain the family of subalgebras $\mathfrak{h}$:
$$
\mathfrak{h}=<e_4 +xe_1 -(x')^2 ze_2, e_5 +x'e_1>.
$$
\item II case 2. Let $ x'+y=0$, then
\begin{itemize}
\item  (3)'  $ (1+xz'+zy')(xy'+y^2)=0$,
\item  (7)'  $ (1+xz+zy)(xz'+yz)=0$,
\item  (8)'  $ y(1+xz'+zy)=0$,
\item  (9)'  $y'(1+xz'+zy)=0$,
\item  (10)' $z'(1+xz'+zy)=0$;
\end{itemize}
\end{enumerate}
\begin{itemize}
\item a) if $1+xz'+zy=0$ then
$$
-x'=y, \; yz'-zy'=0,
$$
and we obtain the family of subalgebras $\mathfrak{h}$
$$
\mathfrak{h}=<e_4 +xe_1 +ye_2 +ze_3,e_5 +ye_1 +y'e_2 +z'e_3>
$$
where $1+xz'+zy=0$ and $yz'-zy'=0$
\item b) if $1+xz'+zy \neq 0$, then $x'=y'=z'=y=0$, and we obtain the family of
subalgebras $\mathfrak{h}$:
$$
\mathfrak{h}=<e_4 +xe_1 +ze_3,e_5>.
$$
\end{itemize}

\; \; \; In result we obtain the collection of Bol algebras with structural
equations given as follows:

\begin{enumerate}
\item $ e_2 \cdot e_3 =-xe_1 +(x')^2 ze_2 -ze_3, (e_2,e_3,e_1)=e_1$;\\
$ e_1 \cdot e_3 =-x'e_1, \; (e_3,e_1,e_2)=e_1$,
\item $ e_2 \cdot e_3 =-xe_1 -ye_2 -ze_3, (e_2,e_3,e_1)=e_1$;\\
$ e_1 \cdot e_3 =ye_1 -y'e_2 -z'e_3, \; (e_3,e_1,e_2)=e_1$,
with condition
$$
1+xz'+zy=0, yz'-zy'=0, \; \; \; (11)
$$
\item $ e_2 \cdot e_3 =-xe_1  -ze_3, (e_2,e_3,e_1)=e_1$;\\
$  \; (e_3,e_1,e_2)=e_1$.
\end{enumerate}

\; \; \;  The group of automorphisms $F$, of Lie triple system $\mathfrak{B}$,
relatively to a fixed base $ e_1, e_2, e_3 $  is defined as follows:

\begin{displaymath}
 P=\left\{P=\left(\begin{array}{ccc}
\alpha & a & d \\
0 & b & f \\
0 & 0 & g \\
\end{array}\right), \alpha \neq 0 ,bg=1 bd=af\right\}.
\end{displaymath}
\; \; \;  The extension of automorphism $P$ from $F$ to the automorphism of Lie
algebra $\mathfrak{G}$, transforming the subspace $\mathfrak{B}$ into itself
can be realized as follows:

$$
Pe_4=[Pe_2,Pe_3]= e_4 -age_5 ,
$$
$$
Pe_5=[Pe_1,Pe_3]= g\alpha e_5.
$$
In addition the examination of Lie subalgebras in (11) gives:
for 1.
$$
P(e_4 +xe_1 -(x')^2 ze_2 +ze_3)=e_4 +(x\alpha -(x')^2 za+zd+ax')e_1 +(zf-(x')^2zb)e_2 +gze_3,
$$
$$
P(e_5 +x'e_1)=e_5 +\frac{x'}{g}e_1.
$$
We will denote:
$$
x_1=x\alpha -(x)^2 za+zd+ax'
$$
$$
y_1=-(x')^2zb+zf,
$$
$$
z_1=gz,
$$
$$
x'_1=\frac{x'}{g},
$$
$$
y'_1=z'_1=0.
$$
\begin{enumerate}
\item I. If $x' \neq 0$, by acting on $g$ we can make $x'_1=1$and here two
cases must be investigate:
\begin{itemize}
\item I.a) if $z=0$, then $z_1=y_1=0$. By acting on $a$ and $\alpha$ we can
make $x_1=0$;
\item I.b) if $z \neq 0$ then $z_1=\beta \neq 0$, hence by acting on $a,f$ and
$\alpha$ we can make $y_1=x_1=0$.
\end{itemize}
\item II. If $x'=0$ then $x'_1=0$;
\begin{itemize}
\item II.a) if $z=0$, then $z_1=y_1=0$. By acting on $\alpha$ and $a$ we can
make $x_1=0$
\item II.b) if $z \neq 0$, then by acting on $g$ we can make $z_1=1$, while
acting on $a,f$ and $\alpha$ we make $y_1=x_1=0$.
\end{itemize}
\end{enumerate}

\; \; \;  In result we obtain for $(x,y,z,x',y',z')$ one family of values and 3
exceptional six-uplet values:
$$
(0,0,0,1,0,0),
$$
$$
(0,0,\beta,1,0,0), \; \beta \neq 0
$$
$$
(0,0,0,0,0,0),
$$
$$
(0,0,1,0,0,0),
$$

\; \; Now we pass to the examination of Bol algebras from the family of Lie subalgebras 2 in (11).
\begin{multline}
P(e_4 +xe_1+ye_2+ze_3)=e_4 +\left[(x+x'\frac{a}{\alpha})\alpha+(y+y'\frac{a}{\alpha})a+(z+z'\frac{a}{\alpha})d\right]e_1 +\\+\left[(y+y'\frac{a}{\alpha})b+(z+z'\frac{a}{\alpha}f)\right]e_2 +(z+z'\frac{a}{\alpha})ge_3,
\end{multline}
$$
P(e_5 -ye_1 +y'e_2 +z'e_3)=e_5 +\frac{-y\alpha +y'a+z'd}{g\alpha}e_1 +\frac{y'b+z'f}{g\alpha}e_2 +\frac{z'}{\alpha}e_3,
$$

where we will denote

$$
x_1=(x\alpha -ya)+(y+y'\frac{a}{\alpha})a+(z+z'\frac{a}{\alpha})d,
$$
$$
y_1=(y+y'\frac{a}{\alpha})b+(z+z'\frac{a}{\alpha}f),
$$
$$
z_1=(z+z'\frac{a}{\alpha})g,
$$
$$
x'_1=\frac{-y\alpha +y'a+z'd}{g\alpha},
$$
$$
y'_1=\frac{y'b+z'f}{g\alpha},
$$
$$
z'_1=\frac{z'}{\alpha},
$$

with the conditions
$$
1+xz'+zy=0 \; \; \; \; \; \; \; \; \; \; (14)
$$
$$
yz'-zy'=0.
$$
\begin{enumerate}
\item I. If $z'=0$, then $z'_1=0$ and from (14) follow that, $zy=-1$ and
make $zy'=0$ then $y'=0, z \neq 0$ hence $y'_1=0$, by acting on $ g$ we can
make
$x'_1=1, z_1=\omega \neq 0$, while acting on $a,f$ and $ \alpha$ we can make
$x_1=y_1=0$.
\item II. If $z' \neq 0$, then by acting on $\alpha$ we can in additional $z'_1=1$
\begin{itemize}
\item II.a) if $y'=0$, the from (14)it follow that, $yz'=0$ hence $y=0$ therefore
$xz'=-1$; hence  acting on $g,a,f$ we can make $y'_1=0;1$, $x'_1=\lambda \neq 0$,
$z_1=0, y_1=0, x_1=s$ (where $s$ any number)
\item II.b) if $y' \neq 0$ and $z=0$ from (14) it follow that: $ xz'=-1, yz'=0$, hence
$y=0$, then  acting on $f$ we can make $y'_1=0$ and $x'_1=t$ (where $t$ is any number)
acting on $g$ we can make $z_1=1, y_1=0, x_1=\gamma$- (any number).
\item II.b.2) If $y'\neq 0$ and $z\neq 0$, then by acting on $f$ and $b$ we make
$y'_1=0, x'_1=\tau$ - any number, $z_1=\kappa \neq 1, y_1=0, x_1=\eta$ -any number.
\end{itemize}
\end{enumerate}

\; \; \;  In result we obtain for $(x,y,z,x',y'z')$ six families of values:

$$
(0,0,\omega,1,0,0) \omega \neq 0,
$$
$$
(s,0,0,\lambda,1,1), \forall s
$$
$$
(s,0,0,\lambda,0,1), \lambda \neq 0, \forall s
$$
$$
(\gamma, 0,1,t,0,1), \forall \gamma, t
$$
$$
(\eta,0,\kappa,\tau,0,1), \kappa \neq 1.
$$

We pass to the examination of Bol algebras from the family 3 of (11)

$$
P(e_4 +xe_1 +xe_3)=e_4+(x\alpha+zd)e_1 +zfe_2 +zye_3,
$$
$$
P(e_5)=e_5,
$$
where we can denote
$$
x_1=x\alpha +zd,
$$
$$
y_1=zf,
$$
$$
z_1=zy,
$$
$$
x'_1=y'_1=z'_1=0.
$$
\begin{itemize}
\item a) If$ z=0$, then $z_1=y_1=0$, therefore $x_1=x\alpha$.
\begin{itemize}
\item a.1. If $x \neq 0$, by choosing $\alpha=\frac{1}{x}$ we can make $x_1=1$.
\item a.2. If $x=0$, then $z_1=y_1=x_1=0$.
\end{itemize}
\item b) If $z \neq 0$, by choosing $g=\frac{1}{z}$ we can make $z_1=1, z=b$, but if
$f=0$ we make $y_1=0$ by acting on $\alpha$ and $d$ we make $x_1=0;1$.
\end{itemize}

\; \;  In result we obtain the following values $(x,y,z,x',y',z')$:
$$
(1,0,1,0,0,0),
$$
$$
(0,0,1,0,0,0),
$$
$$
(1,0,0,0,0,0),
$$
$$
(0,0,0,0,0,0).
$$

\; \;  We require to note that for all the isolated values of the form
$(x,y,z,x',y',z')$, the subalgebras so obtained are no more changing, under the
action of the automorphism of Lie algebras.

\; \;  Summarizing the conducted examination, one can formulate the theorem:

\underline{{\bf Theorem III.7.1}} Any Bol algebra of dimension 3, with the
trilinear
operation of Type VII and the canonical enveloping Lie algebra $\mathfrak{G}$
of dimension 5, is isomorphic to one of Bol algebras below:

\begin{tabular}{rll}
1. &Trivial & $(e_2,e_3,e_1)=e_1$\\
   &bilinear & $(e_3,e_1,e_2)=e_1$\\
   &operation &                   \\
2. &$e_2 \cdot e_3=-e_1$ & $(e_2,e_3,e_1)=e_1$\\
   &                     & $(e_3,e_1,e_2)=e_1$\\
3. &$e_2 \cdot e_3=-e_3$ & $(e_2,e_3,e_1)=e_1$\\
   &                     & $(e_3,e_1,e_2)=e_1$\\
4. &$e_2 \cdot e_3=-e_1$ & $(e_2,e_3,e_1)=e_1$\\
   &                     & $(e_3,e_1,e_2)=e_1$\\
5. &$e_2 \cdot e_3=e_1 -e_3$ & $(e_2,e_3,e_1)=e_1$\\
   &                     & $(e_3,e_1,e_2)=e_1$\\
6. &$e_2 \cdot e_3=-\omega e_1$ & $(e_2,e_3,e_1)=e_1, \omega >0$\\
   &$ e_1 \cdot e_3=-e_1$  & $(e_3,e_1,e_2)=e_1$\\
7. &$e_2 \cdot e_3=-se_1$ & $(e_2,e_3,e_1)=e_1, s \ge 0$\\
   &$ e_1 \cdot e_3=-\lambda e_1 -e_2 -e_3$  & $(e_3,e_1,e_2)=e_1,\lambda >0$\\
8. &$e_2 \cdot e_3=-se_1$ & $(e_2,e_3,e_1)=e_1, s \ge >0$\\
   &$ e_1 \cdot e_3=-\lambda e_1 -e_2$  & $(e_3,e_1,e_2)=e_1, \lambda >0$\\
9. &$e_2 \cdot e_3=-\gamma e_1 -e_3$ & $(e_2,e_3,e_1)=e_1, \gamma >0$\\
   &$ e_1 \cdot e_3=-te_1 -e_3$  & $(e_3,e_1,e_2)=e_1, t>0$\\
10. &$e_2 \cdot e_3=-\eta e_1 -\kappa e_3$ & $(e_2,e_3,e_1)=e_1, \eta,\tau \ge 0$\\
   &$ e_1 \cdot e_3=-\tau e_1 -e_3,$  & $(e_3,e_1,e_2)=e_1 \kappa \ge 0$ but $ \kappa \neq 1$.\\
\end{tabular}

\vspace{12pt}
\; \;  In addition the distinguished Bol algebras are, not isomorphic among
themselves.

\; \;  We pass to the description of Bol 3-Webs corresponding to the isolated
Bol algebras of Type VII. We note for the selected Bol algebras above, we will
limit ourselves to the description of Bol 3-Webs from Bol algebra number 6.

 The composition law  $(\triangle)$, corresponding to the Lie group $G$, of Lie
algebra enveloping Bol algebra is defined as follows:

$$ \begin{bmatrix} x_1\\x_2\\x_3\\x_4\\x_5\end{bmatrix} \triangle
   \begin{bmatrix} y_1\\y_2\\y_3\\y_4\\y_5\end{bmatrix}
   =
   \begin{bmatrix}
     x_{1}+(1+x_2 x_3)y_{1}\exp (x_4 -\frac{x_2 x_3}{2})\\x_{2}+y_{2}\\x_{3}+y_{3}\\x_{4}+y_{4}+\frac{x_{2}y_3 +y_{2}x_{3}}{2}\\ x_5 -x_3y_1 \exp (x_4 -\frac{x_2 x_3}{2})+y_5 \exp (x_4 -\frac{x_2 x_3}{2})\end{bmatrix}.
$$

The collection of elements in $B$ is define as:
$$
B= \exp \mathfrak{B}=\{t(2\int_{0}^{1}\exp (-\frac{uv\alpha^2}{2})\ud \alpha -\exp \frac{-uv}{2}),u,v,\frac{t}{u}(\exp \frac{-uv}{2} -1)\}_{t,u.v \in \mathbb{R}}
$$

For a convenience of examination we will divide it into subcases:
\begin{enumerate}
\item 1. If u=0 and v-any number,
$$
B= \exp \mathfrak{B}=\{t,0,v,\frac{-vt}{2}\}
$$
\item 2. If v=0 and u- any number,

$$
B= \exp \mathfrak{B}=\{t,u,0,0\}.
$$
\item 3. If $u \neq 0$, $v \neq 0$ and have same sign then,
$$
B= \exp \mathfrak{B}=\{t(2\beta \sqrt{\frac{2}{uv}} -\exp \frac{-uv}{2}),u,v,\frac{t}{u}(\exp \frac{-uv}{2} -1)\}_{t,u.v \in \mathbb{R}},
$$
where $\beta=\int_{0}^{1} \exp -p^2 \ud p$ and $p=\frac{uv}{2}\alpha^2$.
\end{enumerate}

\; \; \;  Let Bol algebra be defined as follows:

\begin{tabular}{rll}
   &$e_2 \cdot e_3=-\omega e_1$ & $(e_2,e_3,e_1)=e_1, \omega >0$\\
   &$ e_1 \cdot e_3=-e_1$  & $(e_3,e_1,e_2)=e_1$\\
\end{tabular}

 the subalgebra $ H=\exp \mathfrak{h} $ can be defined as:

$$
H=\exp \mathfrak{h}=\{\exp \alpha (e_4 +\omega e_1),\exp l(e_5 +e_1)\}_{ \alpha \in \mathbb{R}}=\{l+\alpha \omega,0 ,0, \alpha, l \}_{ \alpha,l \in \mathbb{R}}.
$$

The collection of elements $B=\exp \mathfrak{B}$ is:

\begin{enumerate}
\item 1. If u=0 and v-any number'
$$
B= \exp \mathfrak{B}=\{t,0,v,\frac{-vt}{2}\}.
$$
Hence, any element $(x_1,x_2,x_3,x_4,x_5)$ from $G$ can not be always represented as an element of
$B$ and $H$.
\item 2. If v=0 and u- any number,

$$
B= \exp \mathfrak{B}=\{t,u,0,0\}.
$$
In this case too, any element from $G$, can not be always represented as an
element of $B$ and $H$.
\item 3. If $u \neq 0$, $v \neq 0$ and have same sign then,
$$
B= \exp \mathfrak{B}=\{t(2\beta \sqrt{\frac{2}{uv}} -\exp \frac{-uv}{2}),u,v,\frac{t}{u}(\exp \frac{-uv}{2} -1)\}_{t,u.v \in \mathbb{R}}
$$
where $\beta=\int_{0}^{1} \exp -p^2 \ud p$ and $p=\frac{uv}{2}\alpha^2$

Here, any element $(x_1,x_2,x_3,x_4,x_5)$ from $G$ sufficiently in the
neighborhood of $e$ can be uniquely represented as follows:
\end{enumerate}

$$ \begin{pmatrix} x_1\\x_2\\x_3\\x_4\\x_5\end{pmatrix} =
   \begin{bmatrix}P(2\beta\sqrt{\frac{2}{x_2 x_3}}-e^{-\frac{x_2 x_3}{2}})\\x_2\\ x_3\\0\\P\frac{e^{-\frac{x_2 x_3}{2}-1}}{x_2}\end{bmatrix}
   \triangle
   \begin{bmatrix} L+\omega x_4\\0\\0\\x_4\\ L\end{bmatrix}
$$

where $P, L$ are defined from the relations
$$
P=\frac{We^{-\frac{x_2 x_3}{2}}}{S},
$$
$$
L=\frac{Q}{Z},
$$
$$
K=\frac{e^{-\frac{x_2 x_3}{2}}-1}{x_2},
$$
$$
W=[[x_1 -\omega x_4 (1+x_2 x_3)](1-x_3)-(1-x_2 +x_2 x_3)(x_5+ \omega x_3 x_4)]
$$
$$
S=(2\beta \sqrt{\frac{2}{x_2 x_3}}-e^{-\frac{x_2 x_3}{2}})(1-x_3)e^{-\frac{x_2 x_3}{2}}-(1-x_2 +x_2 x_3)e^{-\frac{x_2 x_3}{2}}K,
$$
$$
Q=(2\beta \sqrt{\frac{2}{x_2 x_3}}-e^{-\frac{x_2 x_3}{2}})(x_5+\omega x_3 x_4 e^{-\frac{x_2 x_3}{2}})-(x_1-\omega x_4(1+x_2 x_3)e^{-\frac{x_2 x_3}{2}})K,
$$
$$
Z=(2\beta \sqrt{\frac{2}{x_2 x_3}}-e^{-\frac{x_2 x_3}{2}})(1-x_3)e^{-\frac{x_2 x_3}{2}}-(1-x_2 +x_2 x_3)e^{-\frac{x_2 x_3}{2}}K.
$$

\; \; The composition law $( \star )$, corresponding to the local analytical Bol
loop $B( \star )$ is defined as follows:

$$
= \exp^{-1}
  \left(\prod_B \begin{bmatrix}
           t+t'(1+uv)e^{-\frac{uv}{2}}\\
           u+u'\\
           v+v'\\
           \frac{uv'-vu'}{2}\\
           vt'e^{-\frac{uv}{2}}
  \end{bmatrix} \right)
$$
$$
= \exp^{-1}
  \left( \begin{bmatrix}
   P_1 (2\beta \sqrt{\frac{2}{(u+u')(v+v')}}-e^{-\frac{(u+u')(v+v')}{2}})\\
           u+u'\\
           v+v'\\
           0\\
           T_1 \frac{e^{-\frac{(u+u')(v+v')}{2}}-1}{u+u'}
  \end{bmatrix} \right)
$$
$$
= \begin{bmatrix} P_1\\u+u'\\v+v'\end{bmatrix}.
$$

Where:
$$
T_2 =\left\{t+t'(1+uv)e^{-\frac{uv}{2}}-\omega \frac{uv'-vu'}{2}[1+(u+u')(v+v')]\right\}[1-(u+u')],
$$
$$
R_2 =[1-u-u'+(u+u')(v+v')][vt'e^{-\frac{uv}{2}}+\omega (v+v')\frac{(uv'-vu')}{2}]
$$
$$
M_2 = \left(2\beta \sqrt{\frac{2}{(u+u')(v+v')}}-e^{-\frac{(u+u')(v+v')}{2}}\right)[1-(v+v')]e^{-\frac{(u+u')(v+v')}{2}},
$$
$$
N_2 =[1-u-u'+(u+u')(v+v')]e^{-\frac{(u+u')(v+v')}{2}}\frac{(e^{-\frac{(u+u')(v+v')}{2}}-1)}{u+u'}
$$
$$
P_1 =\frac{T_2- R_2}{M_2 -N_2}.
$$
\newpage

\end{document}